\documentclass[leqno, 10pt]{amsart}
\usepackage{amsmath, amsthm, amssymb, latexsym}
 \newtheorem{definition}{Definition}[section]
 \newtheorem{theorem}[definition]{Theorem}
 \newtheorem{lemma}[definition]{Lemma}
 \newtheorem{proposition}[definition]{Proposition}
 \newtheorem{corollary}[definition]{Corollary}

 \newtheorem*{theorem*}{Theorem}
\newtheorem*{proposition*}{Proposition}
\newtheorem*{lemma*}{Lemma}

 \theoremstyle{remark}
 
 \newtheorem{remark}[definition]{Remark}
 
  \newtheorem*{claim}{Claim}
  \newtheorem*{acknowledgements}{Acknowledgements}

\newcommand{\op}[1]{\operatorname{#1}}

\newcommand{\norminf}[1]{\ensuremath{\|{#1}\|_{(1,\infty)}}}

\newcommand{\acou}[2]{\ensuremath{\langle #1 , #2 \rangle}}

\newcommand{\brak}[1]{\ensuremath{\langle\! #1\!\rangle}}

\newcommand{\Tr}{\ensuremath{\op{Tr}}}
\newcommand{\tr}{\op{tr}}
\newcommand{\Tra}{\ensuremath{\op{Trace}}}
\newcommand{\Trace}{\ensuremath{\op{Trace}}}
\newcommand{\TR}{\ensuremath{\op{TR}}}

\newcommand{\Res}{\ensuremath{\op{Res}}}
\newcommand{\res}{\ensuremath{\op{Res}}}

\newcommand{\Trw}{\ensuremath{\op{Tr}_{\omega}}}
\newcommand{\bint}{\ensuremath{-\hspace{-2,4ex}\int}}

\newcommand{\str}{\op{str}}
\newcommand{\ind}{\op{ind}}

\newcommand{\Hol}{\op{Hol}}

\newcommand{\C}{\ensuremath{\mathbb{C}}} 
\newcommand{\bH}{\ensuremath{\mathbb{H}}} 
\newcommand{\N}{\ensuremath{\mathbb{N}}} 
\newcommand{\R}{\ensuremath{\mathbb{R}}} 
\newcommand{\Z}{\ensuremath{\mathbb{Z}}} 
\newcommand{\CZ}{\ensuremath{\mathbb{C}\!\setminus\!\mathbb{Z}}} 

\newcommand{\Rd}{\ensuremath{\R^{d+1}}}

\newcommand{\Rdo}{\R^{d+1}\!\setminus\! 0}

\newcommand{\URd}{U\times\R^{d+1}}

\newcommand{\URdo}{U\times(\R^{d+1}\!\setminus\! 0)}

\newcommand{\fg}{\ensuremath{\mathfrak{g}}}
\newcommand{\fh}{\ensuremath{\mathfrak{h}}}

\newcommand{\cD}{\ensuremath{\mathcal{D}}}
\newcommand{\cE}{\ensuremath{\mathcal{E}}}
\newcommand{\cF}{\ensuremath{\mathcal{F}}}

\newcommand{\cH}{\ensuremath{\mathcal{H}}}

\newcommand{\cK}{\ensuremath{\mathcal{K}}}
\newcommand{\cL}{\ensuremath{\mathcal{L}}}
\newcommand{\cM}{\ensuremath{\mathcal{M}}}
\newcommand{\cN}{\ensuremath{\mathcal{N}}}
\newcommand{\cS}{\ensuremath{\mathcal{S}}}
\newcommand{\cU}{\ensuremath{\mathcal{U}}}

\newcommand{\lunf}{\ensuremath{\mathcal{L}^{(1,\infty)}}}

\newcommand{\psivdo}{$\Psi_{H}$DO}
\newcommand{\psivdos}{$\Psi_{H}$DOs}
 
\newcommand{\pvdo}{\ensuremath{\Psi_{H}}}

\newcommand{\pvdoi}{\ensuremath{\Psi_{H}^{\text{int}}}} 
\newcommand{\pvdoz}{\ensuremath{\Psi_{H}^{\Z}}} 
\newcommand{\pvdocz}{\ensuremath{\Psi_{H}^{\CZ}}} 
\newcommand{\pvdoc}{\ensuremath{\Psi_{H,\op{c}}}}

\newcommand{\psido}{$\Psi$DO} 
\newcommand{\psidos}{$\Psi$DOs} 
\newcommand{\psinf}{\ensuremath{\Psi^{-\infty}}} 
\newcommand{\psinfc}{\ensuremath{\Psi^{-\infty}_{c}}}

\newcommand{\Shi}{\ensuremath{S^{\text{int}}}}

\newcommand{\reg}{{\text{reg}}}
\newcommand{\ord}{{\op{ord}}}

\newcommand{\xiy}{{\xi\rightarrow y}}
\newcommand{\yxi}{{y\rightarrow\xi}}

\newcommand{\supp}{\op{supp}}

\newcommand{\im}{\op{im}}

\newcommand{\End}{\ensuremath{\op{End}}}
\newcommand{\END}{\ensuremath{\op{END}}}

\newcommand{\updown}{\uparrow \downarrow}

\newcommand{\hotimes}{\hat\otimes}

\renewcommand{\Box}{\square}

\newcommand{\dbarb}{\overline{\partial}_{b}}
\newcommand{\dbarbpq}{\overline{\partial}_{b;p,q}}

\newcommand{\Boxbpq}{\Box_{b;p,q}}

\newcommand{\Vol}{\op{Vol}}

\begin{document}
\title{Noncommutative residue for Heisenberg manifolds. Applications in CR and contact geometry} 

\author{Rapha\"el Ponge}

\address{Department of Mathematics, University of Toronto, Canada.}
\email{ponge@math.toronto.edu}


\begin{abstract}
    This paper has four main parts. In the first part, we construct a noncommutative residue for the hypoelliptic calculus on Heisenberg manifolds,  
    that is, for the class of \psivdo\ operators introduced by Beals-Greiner and Taylor. This noncommutative residue appears as the residual trace 
     on integer order \psivdos\  induced by the analytic extension of the usual 
    trace to non-integer order \psivdos.  Moreover, it agrees with the integral of the density defined by the logarithmic singularity of the Schwartz kernel of 
    the corresponding \psivdo. In addition, we show that this noncommutative residue provides us with the unique trace up to constant multiple on 
    the algebra of integer order \psivdos. In the second part, we give some analytic applications of this construction concerning zeta functions of 
    hypoelliptic operators, logarithmic metric estimates for Green kernels of hypoelliptic operators, and the extension of the 
    Dixmier trace to the whole algebra of integer order \psivdos. In the third part, we present examples of 
computations of noncommutative residues of some powers of the horizontal sublaplacian and the contact Laplacian on contact manifolds. In the fourth part, 
we present two applications in CR geometry. First, we give 
some examples of geometric computations of noncommutative residues of some powers of the horizontal sublaplacian and of the Kohn 
Laplacian. Second, we 
make use of the framework of noncommutative geometry and of our noncommutative residue to define lower dimensional volumes in pseudohermitian 
geometry, e.g., we can give sense to the area of any 3-dimensional CR manifold. On the way we obtain a spectral interpretation of the 
Einstein-Hilbert action in pseudohermitian geometry.
\end{abstract}

\maketitle 
\numberwithin{equation}{section}

    
\section{Introduction}
The aim of this paper is to construct a noncommutative residue trace for the Heisenberg calculus and to present several of its applications, 
in particular in CR and contact geometry.
The Heisenberg calculus was built independently by Beals-Greiner~\cite{BG:CHM} and 
Taylor~\cite{Ta:NCMA}  as the relevant pseudodifferential tool to study the main geometric operators on contact and CR manifolds, which fail to 
be elliptic, but may be hypoelliptic (see also~\cite{BdM:HODCRPDO}, \cite{EM:HAITH}, 
\cite{FS:EDdbarbCAHG}, \cite{Po:MAMS1}). This calculus holds in the general 
setting of a Heisenberg manifold, that is, a manifold $M$ together with a distinguished hyperplane bundle $H\subset TM$, and we construct a 
noncommutative residue trace in this general context. 

The noncommutative residue trace of Wodzicki~(\cite{Wo:LISA}, \cite{Wo:NCRF}) and Guillemin~\cite{Gu:NPWF} was originally constructed for classical 
\psidos\ and it appears as the residual trace on integer order \psidos\ induced by analytic extension of the operator trace to \psidos\ of non-integer 
order. It has numerous applications and generalizations 
(see, e.g., \cite{Co:AFNG}, \cite{Co:GCMFNCG}, \cite{CM:LIFNCG}, \cite{FGLS:NRMB}, \cite{Gu:RTCAFIO}, \cite{Ka:RNC}, \cite{Le:NCRPDOLPS}, 
\cite{MMS:FIT}, \cite{MN:HPDO1}, \cite{PR:CDBFCF}, \cite{Po:IJM1}, \cite{Sc:NCRMCS}, \cite{Vas:PhD}). 
In particular, the existence of a residual trace is an essential ingredient in the framework for the local 
index formula in noncommutative geometry of Connes-Moscovici~\cite{CM:LIFNCG}. 

Accordingly, the noncommutative residue for the Heisenberg calculus has various applications and several of them are presented in this paper.
Further geometric applications can be found in~\cite{Po:Crelle2}.

%
\subsection{Noncommutative residue for Heisenberg manifolds}
Our construction of a noncommutative residue trace for \psivdos, i.e., for the pseudodifferential operators in the Heisenberg calculus,  
follows the approach of~\cite{CM:LIFNCG}. It has two main ingredients:\smallskip 

(i) The observation that the coefficient of 
the logarithmic singularity of the Schwartz kernel of a \psivdo\ operator $P$ can be defined globally as a density $c_{P}(x)$ functorial with respect to the 
action of Heisenberg diffeomorphisms, i.e., diffeomorphisms preserving the Heisenberg structure  
(see Proposition~\ref{thm:NCR.log-singularity}).\smallskip

(ii) The analytic extension of the operator trace to \psivdos\ of complex non-integer order (Proposition~\ref{thm:NCR.TR.global}).\smallskip 

The analytic extension of the trace from (ii) is obtained by working directly at the level of densities and induces on \psivdos\ of integer order a 
residual trace given by (minus) the integral of the density from (i) (Proposition~\ref{thm:NCR.TR.local}). 
This residual trace is our noncommutative residue for the 
Heisenberg calculus. 

In particular, as an immediate byproduct of this construction the noncommutative residue is invariant under the action of 
Heisenberg diffeomorphisms. Moreover, in the foliated case our noncommutative residue agrees with that of~\cite{CM:LIFNCG}, and on the algebra of Toeplitz 
pseudodifferential operators on a contact manifold of Boutet de Monvel-Guillemin~\cite{BG:STTO} we recover the noncommutative residue of 
Guillemin~\cite{Gu:RTCAFIO}. 
%
%

As a first application of this construction we show that when the Heisenberg manifold is connected the noncommutative residue is the unique trace up to constant 
multiple  on the algebra of integer order \psivdos~(Theorem~\ref{thm:Traces.traces}).  As a 
consequence we get a characterization sums of \psivdo\ commutators and we obtain that any smoothing operator can be written as a sum of \psivdo\ 
commutators. 

These results are the analogues for \psivdos\ of well known results of Wodzicki~(\cite{Wo:PhD}; see also~\cite{Gu:RTCAFIO}) 
for classical \psidos. Our arguments are somewhat elementary and partly rely on the characterization of the Schwartz kernels of 
\psivdos\ that was used in the analysis of their logarithmic singularities near the diagonal.
%

\subsection{Analytic applications on general Heisenberg manifolds}
The analytic extension of the trace allows us to directly define the zeta function 
  $\zeta_{\theta}(P;s)$ of a hypoelliptic \psivdo\ operator $P$ as a meromorphic functions on $\C$. The definition depends on the choice of a ray $L_{\theta}=\{\arg 
  \lambda =\theta\}$, $0\leq \theta <2\pi$, which is a ray of principal values for the principal 
symbol of $P$ in the sense of~\cite{Po:CPDE1} and is not through an eigenvalue of $P$, so that $L_{\theta}$ is a ray of minimal growth for $P$. 
Moreover, the residues at the potential singularity points of  $\zeta_{\theta}(P;s)$ can be expressed as noncommutative residues.  

When the set of principal values of the principal symbol of $P$ contains the left half-plane $\Re \lambda\leq 0$ we further can relate the residues 
and regular values of $\zeta_{\theta}(P;s)$ to the coefficients in the heat kernel asymptotics for $P$ (see 
Proposition~\ref{prop:Zeta.heat-zeta-global} for the precise statement). 
We then use this to derive a local formula for the index of a hypoelliptic \psivdo\ and to rephrase in terms of noncommutative residues 
the Weyl asymptotics for hypoelliptic \psidos\ from~\cite{Po:MAMS1} and~\cite{Po:CPDE1}. 

An interesting application concerns logarithmic metric estimates for Green kernels of hypoelliptic \psivdos. 
It is not true that a positive hypoelliptic \psivdo\ has a Green kernel positive near the diagonal. Nevertheless, by making use of the 
spectral interpretation of the noncommutative residue as a residual trace, we show that the positivity still pertains when the order is equal to the critical 
dimension $\dim M+1$ (Proposition~\ref{prop:Metric.positivity-cP}). 

When the bracket condition $H+[H,H]=TM$ holds, i.e., $H$ is a 
Carnot-Carath\'eodory distribution, this allows us to get metric estimates in terms of the Carnot-Carath\'eodory metric associated to any given subriemannian 
metric on $H$ (Theorem~\ref{thm:Metric.metric-estimate}). 
This result connects nicely with the work of Fefferman, Stein and their collaborators on metric estimates 
for Green kernels of subelliptic sublaplacians on general Carnot-Carath\'eodory manifolds 
(see, e.g.,~\cite{FS:FSSOSO}, \cite{Ma:EPKLCD}, \cite{NSW:BMDVF1}, \cite{Sa:FSGSSVF}).

Finally, we show that on a Heisenberg manifold $(M,H)$ the Dixmier trace is defined for \psivdos\ of order less than or equal to the 
critical order 
$-(\dim M+1)$ and on such operators agrees with the noncommutative residue (Theorem~\ref{thm:NCG.Dixmier}). 
Therefore, the noncommutative residue allows us to extend the Dixmier trace to the whole algebra of \psivdos\ of integer order. 
In noncommutative geometry the 
Dixmier trace plays the role of the integral on infinitesimal operator of order~$\leq 1$. Therefore, our result allows us to integrate any \psivdo\ even 
though it is not an infinitesimal operator of 
order~$\leq 1$. This is the analogue of a well known result of Connes~\cite{Co:AFNG} for classical \psidos. 

\subsection{Noncommutative residue and contact geometry}
Let $(M^{2n+1},H)$ be a compact orientable contact manifold, so that the 
 hyperplane bundle $H\subset TM$ can be realized as the kernel of a contact form $\theta$ on $M$. The additional datum of a \emph{calibrated} almost 
 complex structure on $H$ defines a Riemannian metric on $M$ whose volume $\Vol_{\theta}M$ depends only on $\theta$. 
 
 Let $\Delta_{b;k}$ be the horizontal sublaplacian associated to the above Riemannian metric acting on horizontal forms of degree $k$, $k\neq n$. This 
operator is hypoelliptic for $k\neq n$ and by making use of the results of~\cite{Po:MAMS1} 
we can explicitly express the noncommutative residue of $\Delta_{b;k}^{-(n+1)}$ as a constant multiple of $\Vol_{\theta}M$ 
(see Proposition~\ref{prop:Contact.residue-Deltab}). 
 
 Next, the contact complex of Rumin~\cite{Ru:FDVC} is a complex of horizontal forms on a contact manifold whose Laplacians are hypoelliptic in every 
 bidegree.  Let $\Delta_{R;k}$ denote the contact Laplacian acting on forms degree $k$, $k=0,\ldots,n$.  
 Unlike the horizontal sublaplacian $\Delta_{R}$ does not act on all horizontal forms, but on  the sections of a subbundle 
 of horizontal forms. Moreover, it is not a sublaplacian and it even has order 4 on forms of degree $n$. Nevertheless, by making use of the results 
 of~\cite{Po:MAMS1} we can show that the noncommutative residues of $\Delta_{R;k}^{-(n+1)}$ for $k\neq n$ and of $\Delta_{R;n}^{-\frac{n+1}{2}}$ are 
 universal constant multiples of the contact volume $\Vol_{\theta}M$ (see Proposition~\ref{prop:Contact.residue-DeltaR}). 

 \subsection{Applications in CR geometry}
 Let $(M^{2n+1},H)$ be a compact orientable $\kappa$-strictly 
 pseudoconvex CR manifold equipped with a pseudohermitian contact form $\theta$, i.e., the hyperplane bundle $H\subset TM$ has an (integrable) complex 
 structure and the Levi form associated to $\theta$ has at every point $n-\kappa$ positive eigenvalues and $\kappa$ negative eigenvalues. If $h$ is a 
 Levi metric on $M$ then the volume with respect to this metric depends only on $\theta$ and is denoted $\Vol_{\theta}M$. 
 
 As in the general contact case we can explicitly relate the pseudohermitian volume $\Vol_{\theta}M$ to the noncommutative residues of the following 
 operators:\smallskip
 
 - $\Box_{b;pq}^{-(n+1)}$, where $\Box_{b;pq}$ denotes the Kohn Laplacian acting on $(p,q)$-forms with $q\neq \kappa$ and 
 $q\neq n-\kappa$ (see Proposition~\ref{prop:CR.residue-Boxb1});\smallskip
 
 - $\Delta_{b;pq}^{-(n+1)}$, where $\Delta_{b;pq}$ denotes the horizontal sublaplacian acting on $(p,q)$-forms with $(p,q)\neq (n-\kappa,\kappa)$ and 
 $(p,q)\neq (\kappa,n-\kappa)$ (see Proposition~\ref{prop:CR.residue-Deltab1}).\smallskip
 
 From now  on we assume $M$ strictly pseudoconvex (i.e.~we have $\kappa=0$) and consider the following operators:\smallskip 
 
-  $\Box_{b;pq}^{-n}$, with $q\neq 0$ and $q\neq n$,;\smallskip

- $\Delta_{b;pq}^{-n}$, with $(p,q)\neq (n,0)$ and $(p,q)\neq (0,n)$.\smallskip 

\noindent Then we can make use of the results of~\cite{BGS:HECRM} to express 
 the noncommutative residues of these operators as universal constant multiple of the integral $\int_{M}R_{n}d\theta^{n}\wedge \theta$, where $R_{n}$ 
 denotes the scalar curvature of the connection of Tanaka~\cite{Ta:DGSSPCM} 
 and Webster~\cite{We:PHSRH}  (see Propositions~\ref{prop:CR.residue-Boxb2} and~\ref{prop:CR.residue-Deltab2}). 
  These last results provide us with a spectral interpretation of the Einstein-Hilbert action in pseudohermitian geometry, which is 
  analogous to that of Connes~(\cite{Co:GCMFNCG}, \cite{KW:GNGWR}, \cite{Ka:DOG}) in the Riemannian case. 
 
%
%
 
 Finally, by using an idea of Connes~\cite{Co:GCMFNCG} we can make use of the noncommutative residue for classical \psidos\ to define the $k$-dimensional 
  volumes Riemannian manifold of dimension $m$ for $k=1,\ldots,m-1$, e.g.~we can give sense to the area in any dimension (see~\cite{Po:LMP07}). 
  Similarly, we can make use of the noncommutative residue for the Heisenberg calculus to define the $k$-dimensional pseudohermitian volume 
  $\Vol^{(k)}_{\theta}M$ for any $k=1,\ldots,2n+2$. The argument involves noncommutative geometry, but we can give a purely differential geometric 
  expression of these lower dimensional volumes (see~Proposition~\ref{prop:CR.lower-dim.-volumes}). Furthermore, in dimension 3 the area (i.e.~the 
  2-dimensional volume) is a constant multiple 
  of the integral of the Tanaka-Webster scalar curvature (Theorem~\ref{thm:spectral.area}). In particular, we find that the area of the sphere 
  $S^{3}\subset \C^{2}$ endowed with its standard pseudohermitian structure has area $\frac{\pi^{2}}{8\sqrt{2}}$.

  \subsection{Potential geometric applications}
  The boundaries of a strictly pseudoconvex domain of $\C^{n+1}$ naturally carry strictly pseudoconvex CR structures, so one can expect 
  the above  results to be useful for studying from the point of view of noncommutative geometry strictly pseudoconvex boundaries, and more generally Stein 
  manifolds with boundaries and the asymptotically complex hyperbolic manifolds of~\cite{EMM:RLSPD}. Similarly, the boundary of a symplectic manifold
  naturally inherits a contact structure, so we could also use the results of this 
  papers to give a noncommutative geometric study of symplectic manifolds with boundary. 
  
  Another interesting potential application concerns a special class of Lorentzian manifolds, the Fefferman's spaces~(\cite{Fe:MAEBKGPCD}, \cite{Le:FMPHI}). 
  A Fefferman's Lorentzian space $\cF$ can be realized as the total space of a circle bundle over a strictly 
  pseudoconvex   CR manifold $M$ and it carries a Lorentzian metric naturally associated to any pseudohermitian contact form on $M$. 
  For instance, the
  curvature tensor of $\cF$ can be explicitly expressed in terms of the curvature and torsion tensors of the Tanaka-Webster connection of $M$ and the 
  Dalembertian of $\cF$ pushes 
  down to the horizontal sublaplacian on $M$. This strongly suggests that one could deduce a noncommutative geometric study of Fefferman spaces 
  from a noncommutative geometric study of strictly pseudoconvex CR manifolds. An item of special interest would be to get a spectral interpretation of 
  the Einstein-Hilbert action in this setting.
  
  Finally, it would be interesting to extend the results of this paper to other subriemannian geometries such as the quaternionic contact manifolds of 
  Biquard~\cite{Bi:MEAS}. 
 
\subsection{Organization of the paper}
The rest of the paper is organized as follows. 

In Section~\ref{sec:Heisenberg-calculus}, 
we recall the main facts about Heisenberg manifold and the Heisenberg 
calculus. 

In Section~\ref{sec:NCR}, we study the logarithmic singularity of the Schwartz kernel of a \psivdo\ and show that it gives rise to a well defined density. We 
then construct the noncommutative residue for the Heisenberg calculus as the residual trace induced on integer order \psivdos\ by the analytic 
extension of the usual trace to non-integer order \psivdos. Moreover, we show that the noncommutative residue of an integer order \psivdo\ 
agrees with the integral of the density defined by the logarithmic singularity of its Schwartz kernel. We end the section by proving that, when the Heisenberg 
manifold is connected, the noncommutative residue is the only trace up to constant multiple. 

In Section~\ref{sec:Analytic-Applications}, we give some analytic applications of the construction of the noncommutative residue. First, we deal with zeta functions of 
hypoelliptic \psivdos\ and relate their singularities to the heat kernel asymptotics of the corresponding operators. Second,  
we prove logarithmic metric estimates for Green kernels of hypoelliptic \psivdos. Finally, we show that the noncommutative residue 
allows us to extend the Dixmier trace to \emph{all} integer order \psivdos.

In Section~\ref{sec:Contact}, we present examples of 
computations of noncommutative residues of some powers of the horizontal sublaplacian and of the contact Laplacian of Rumin on contact manifolds. 

In Section~\ref{sec:CR}, we present some applications in CR geometry. First, we give 
some examples of geometric computations of noncommutative residues of some powers of the horizontal sublaplacian and of the Kohn 
Laplacian.  Second, we make use of the framework of noncommutative geometry and of the noncommutative residue for the Heisenberg calculus 
to define lower dimensional volumes in pseudohermitian geometry. 

Finally, in Appendix for reader's convenience we present a detailed proof of Lemma~\ref{lem:Heisenberg.extension-symbol} about the  
extension of a homogeneous symbol into a homogeneous distribution. This is needed for the analysis of the logarithmic singularity of 
the Schwartz kernel of a \psivdo\ in Section~\ref{sec:NCR}. 

\begin{acknowledgements}
Part of the results of this paper were announced in~\cite{Po:CRAS1} and~\cite{Po:CRAS2} and were presented as part of my PhD thesis at the University of 
Paris-Sud (Orsay, France) in December 2000. I am grateful to my advisor, Alain Connes, and to 
Charlie Epstein, Henri Moscovici and Michel Rumin, for stimulating and helpful discussions related to the subject matter of this 
paper.  In addition, I would like to thank Olivier Biquard, Richard Melrose and Pierre Pansu for their interests in the results of this paper. 
\end{acknowledgements}

\section{Heisenberg calculus}
 \label{sec:Heisenberg-calculus}
 The Heisenberg calculus is the relevant pseudodifferential calculus to study hypoelliptic 
 operators on Heisenberg manifolds. It was independently introduced by Beals-Greiner~\cite{BG:CHM} and Taylor~\cite{Ta:NCMA} (see 
 also~\cite{BdM:HODCRPDO}, \cite{Dy:POHG}, \cite{Dy:APOHSC}, \cite{EM:HAITH}, \cite{FS:EDdbarbCAHG}, \cite{Po:MAMS1}, \cite{RS:HDONG}).
In this section we recall the main facts about the Heisenberg calculus following the point of view of~\cite{BG:CHM} and~\cite{Po:MAMS1}. 

\subsection{Heisenberg manifolds}
In this subsection we gather the main definitions and examples concerning Heisenberg manifolds and their tangent Lie group bundles. 

\begin{definition}
   1) A Heisenberg manifold is a pair $(M,H)$ consisting of a manifold $M$ together with a distinguished hyperplane bundle $H 
\subset TM$. 

2) Given Heisenberg manifolds  $(M,H)$ and $(M',H')$ a diffeomorphism $\phi:M\rightarrow M'$ is said to be a Heisenberg diffeomorphism when 
$\phi_{*}H=H'$.  
\end{definition}

\noindent Following are the main examples of Heisenberg manifolds:\smallskip

\noindent \emph{- Heisenberg group.} The $(2n+1)$-dimensional Heisenberg group
$\bH^{2n+1}$ is the 2-step nilpotent group consisting of $\R^{2n+1}=\R \times \R^{2n}$ equipped with the 
group law,
\begin{equation}
    x.y=(x_{0}+y_{0}+\sum_{1\leq j\leq n}(x_{n+j}y_{j}-x_{j}y_{n+j}),x_{1}+y_{1},\ldots,x_{2n}+y_{2n}).  
\end{equation}
A left-invariant basis for its Lie algebra $\fh^{2n+1}$ is then
provided by the vector fields, 
\begin{equation}
    X_{0}=\frac{\partial}{\partial x_{0}}, \quad X_{j}=\frac{\partial}{\partial x_{j}}+x_{n+j}\frac{\partial}{\partial 
        x_{0}}, \quad X_{n+j}=\frac{\partial}{\partial x_{n+j}}-x_{j}\frac{\partial}{\partial 
        x_{0}}, \quad 1\leq j\leq n.
\end{equation}
 For $j,k=1,\ldots,n$ and $k\neq j$ we have the Heisenberg relations
$[X_{j},X_{n+k}]=-2\delta_{jk}X_{0}$ and $[X_{0},X_{j}]=[X_{j},X_{k}]=[X_{n+j},X_{n+k}]=0$.
In particular, the subbundle spanned by the vector field 
$X_{1},\ldots,X_{2n}$ yields a left-invariant Heisenberg structure on 
$\bH^{2n+1}$.\smallskip

\noindent \emph{- Foliations.} A (smooth) foliation is a manifold $M$ together with a subbundle $\cF \subset TM$ 
integrable in Frobenius' sense, i.e., the space of sections of $H$ is closed under the Lie bracket of vector fields.
Therefore, any codimension 1 foliation is a Heisenberg manifold.\smallskip  

\noindent \emph{- Contact manifolds.} Opposite to foliations are contact manifolds. A contact manifold is a Heisenberg manifold $(M^{2n+1}, H)$ such that $H$ can 
be locally realized as 
the kernel of a contact form, that is, a $1$-form $\theta$ such that $d\theta_{|H}$ is nondegenerate. When $M$ is orientable it is equivalent to 
require $H$ to be globally the kernel of a contact form. Furthermore, by
Darboux's theorem any contact manifold is locally Heisenberg-diffeomorphic to the Heisenberg group $\bH^{2n+1}$ equipped with the standard contact
form $\theta^{0}= dx_{0}+\sum_{j=1}^{n}(x_{j}dx_{n+j}-x_{n+j}dx_{j})$.\smallskip

\noindent \emph{- Confoliations.} According to Elyashberg-Thurston~\cite{ET:C} a \emph{confoliation structure} on an oriented manifold
$M^{2n+1}$ is given by a global non-vanishing $1$-form $\theta$ on $M$ such that
$(d\theta)^{n}\wedge \theta\geq 0$. In particular, if we let $H=\ker \theta$ then $(M,H)$ is a Heisenberg manifold which is a foliation when 
$d\theta\wedge \theta=0$ and a contact manifold when $(d\theta)^{n}\wedge \theta>0$.\smallskip

\noindent \emph{- CR manifolds.} A CR
structure on an orientable manifold $M^{2n+1}$ is given by a rank $n$
complex subbundle $T_{1,0}\subset T_{\C}M$ such that $T_{1,0}$ is integrable in  Frobenius' sense and we have
$T_{1,0}\cap T_{0,1}=\{0\}$, where we have set $T_{0,1}=\overline{T_{1,0}}$. 
Equivalently, the subbundle $H=\Re (T_{1,0}\otimes T_{0,1})$ has the 
structure of a complex bundle of (real) dimension $2n$. In
particular, $(M,H)$ is a Heisenberg manifold. The main example of a CR manifold is that of the (smooth) boundary
$M=\partial D$ of a bounded complex domain $D \subset \C^{n+1}$. In particular,
when $D$ is strongly pseudoconvex with 
defining function $\rho$ the 1-form $\theta=i(\partial
-\bar{\partial})\rho$ is a contact form on $M$.\smallskip 

Next, the terminology Heisenberg manifold stems from the fact that the relevant tangent structure in this setting is that of a bundle $GM$ of graded nilpotent Lie 
groups (see~\cite{BG:CHM}, \cite{Be:TSSRG}, \cite{EMM:RLSPD}, \cite{FS:EDdbarbCAHG},  \cite{Gr:CCSSW}, \cite{Po:Pacific1}, \cite{Ro:INA}, 
\cite{Va:PhD}).  
This tangent Lie group bundle can be described as follows. 

First, there is an intrinsic Levi form  $\cL:H\times H\rightarrow TM/H$ such that, for any point $a 
\in M$ and any sections $X$ and $Y$ of $H$ near $a$, we have 
\begin{equation}
    \cL_{a}(X(a),Y(a))=[X,Y](a) \qquad \bmod H_{a}.
     \label{eq:Heisenberg.Levi-form}
\end{equation}
In other words the class of $[X,Y](a)$ modulo $H_{a}$ depends only on the values $X(a)$ and $Y(a)$, not on the germs of $X$ and $Y$ near $a$ 
(see~\cite{Po:Pacific1}). This allows us to define the tangent Lie algebra bundle $\fg M$ as the vector bundle $(TM/H)\oplus H$ together with the 
grading and field  of Lie brackets such that, for sections $X_{0}$, $Y_{0}$ of $TM/H$ and $X'$, $Y'$ of $H$, we have  
 \begin{gather}
     t.(X_{0}+X')=t^{2}X_{0}+t X', \qquad t\in \R,   
     \label{eq:Heisenberg.Heisenberg-dilations}\\
     [X_{0}+X',Y_{0}+Y']_{\fg M}=\cL(X',Y').
   \end{gather}

Since each fiber $\fg_{a}M$ is 2-step nilpotent, $\fg M$ is the Lie algebra bundle of a Lie group bundle $GM$ which can be realized as  $(TM/H)\oplus H$ 
together with the field of group law such that, 
for sections $X_{0}$, $Y_{0}$ of $TM/H$ and $X'$, $Y'$ of $H$, we have 
 \begin{equation}
     (X_{0}+X').(Y_{0}+Y')=X_{0}+Y_{0}+\frac{1}{2}\cL(X',Y')+X'+Y'.
      \label{eq:Heisenberg.group-law}
 \end{equation}
We call $GM$ the \emph{tangent Lie group bundle} of $M$.

Let $\phi$ be a Heisenberg diffeomorphism from $(M,H)$ onto a Heisenberg manifold $(M',H')$. Since we have $\phi_{*}H=H'$ the linear differential $\phi'$ induces 
linear vector bundle isomorphisms $\phi':H\rightarrow H'$ and $\overline{\phi'}:TM/H\rightarrow TM'/H'$, so that we get a linear vector bundle 
isomorphism $\phi_{H}':(TM/H)\oplus H\rightarrow (TM'/H')\oplus H'$ by letting 
\begin{equation}
    \phi_{H}'(a).(X_{0}+X')= \overline{\phi'}(a)X_{0}+\phi'(a)X',
     \label{eq:Heisenberg.tangent-map}
\end{equation}
for any $a \in M$ and any $X_{0}$ in $(T_{a}M/H_{a})$ and $X'$ in $H_{a}$. This isomorphism commutes with the dilations 
in~(\ref{eq:Heisenberg.Heisenberg-dilations}) 
and it can be further shown that it gives rise to a Lie group isomorphism from $GM$ onto $GM'$ (see~\cite{Po:Pacific1}). 

The above description of $GM$ can be related to the extrinsic approach of~\cite{BG:CHM} as follows. 

\begin{definition}
    A local frame $X_{0},X_{1},\ldots,X_{d}$ of $TM$ such that $X_{1},\ldots,X_{d}$ span $H$ is called a $H$-frame. 
\end{definition}

Let $U \subset \Rd$ be an open of local coordinate equipped with a $H$-frame $X_{0},\ldots,X_{d}$.

   \begin{definition}\label{def:Heisenberg-privileged-coordinates}
For $a\in U$ 
we let $\psi_{a}:\Rd \rightarrow \Rd$ denote the unique affine change of variable such that $\psi_{a}(a)=0$ and 
 $(\psi_{a})_{*}X_{j}(0)=\frac{\partial}{\partial x_{j}}$ for $j=0,\ldots,d$.  The coordinates provided by  the map 
    $\psi_{a}$  are called privileged coordinates centered at $a$.   
   \end{definition}    
    
    In addition, on $\Rd$ we consider the dilations, 
    \begin{equation}
        t.x=(t^{2}x_{0},tx_{1},\ldots,tx_{d}), \qquad t \in \R.
        \label{eq:Heisenberg.Heisenberg-dilations-Rd}
    \end{equation}
    
    In privileged coordinates centered at $a$ we can write $X_{j}=\frac{\partial}{\partial x_{j}}+\sum_{k=0}^{d}a_{jk}(x)\frac{\partial}{\partial x_{j}}$ with 
$a_{jk}(0)=0$. Let $X_{0}^{(a)}=\frac{\partial}{\partial x_{0}}$ and for $j=1,\ldots,d$ let 
$X_{j}^{(a)}=\frac{\partial}{\partial_{x_{j}}}+\sum_{k=1}^{d}b_{jk}x_{k} 
    \frac{\partial}{\partial_{x_{0}}}$, 
where $b_{jk}= \partial_{x_{k}}a_{j0}(0)$.  
With respect to the dilations~(\ref{eq:Heisenberg.Heisenberg-dilations-Rd}) the vector field $X_{j}^{(a)}$ is homogeneous of degree $w_{0}=-2$ for $j=0$ 
and of degree $w_{j}=-1$ for $j=1,\ldots,d$. In fact, using Taylor expansions at $x=0$ we get a formal expansion 
$X_{j} \sim X_{j}^{(a)}+X_{j,w_{j}-1}+\ldots$, 
with $X_{j,l}$ homogeneous vector field of degree $l$. 

The subbundle spanned by the vector fields $X_{j}^{(a)}$ is a 2-step nilpotent Lie algebra under the Lie bracket of vectors fields. Its associated 
Lie group $G^{(a)}$ can be realized as $\Rd$ equipped with the group law, 
\begin{equation}
    x.y=(x_{0}+\sum_{j,k=1}^{d}b_{kj}x_{j}x_{k},x_{1},\ldots,x_{d}).
\end{equation}
    
On the other hand, the vectors $X_{0}(a),\ldots,X_{d}(a)$ provide us with a linear basis of the space $(T_{a}M/H_{a})\oplus H_{a}$. 
This allows us to identify $G_{a}M$ with 
$\Rd$ equipped with the group law,
\begin{equation}
    x.y=(x_{0}+y_{0}+\frac{1}{2}L_{jk}(a)x_{j}y_{k},x_{1}+y_{1},\ldots,x_{d}+y_{d}). 
    \label{eq:Heisenberg.group-law-tangent-group-coordinates}
\end{equation}
Here the functions $L_{jk}$ denote the coefficients of the Levi form~(\ref{eq:Heisenberg.Levi-form}) 
with respect to the $H$-frame $X_{0},\ldots,X_{d}$, i.e., we have $\cL(X_{j},X_{k})=[X_{j},X_{k}]=L_{jk}X_{0} 
\bmod H$. 

The Lie group $G^{(a)}$ is isomorphic to $G_{a}M$ since one can check that $L_{jk}=b_{jk}-b_{kj}$. An explicit isomorphism is given by
\begin{equation}
    \phi_{a}(x_{0},\ldots,x_{d})= (x_{0}-\frac{1}{4}\sum_{j,k=1}^{d}(b_{jk}+b_{kj})x_{j}x_{k},x_{1},\ldots,x_{d}). 
\end{equation}

\begin{definition}\label{def:Heisenberg-Heisenberg-coordinates}
 The local coordinates 
provided by 
the map $\varepsilon_{a}:=\phi_{a}\circ\psi_{a}$ are called Heisenberg coordinates centered at $a$.    
\end{definition}

The Heisenberg coordinates refines the privileged coordinates in such way that the above realizations of $G^{(a)}$ and $G_{a}M$ agree. In particular, 
the vector fields $X_{j}^{(a)}$ and $X_{j}^{a}$ agree in these coordinates. This allows us to see $X_{j}^{a}$ as a first order approximation of $X_{j}$. 
For this reason $X_{j}^{a}$ is called the \emph{model vector field of $X_{j}$} at $a$. 

   \subsection{Left-invariant pseudodifferential operators} Let $(M^{d+1},H)$ be a Heisenberg manifold and let $G$ be the tangent group $G_{a}M$ of $M$ 
  at a given point $a\in M$. We briefly recall the calculus for homogeneous left-invariant \psidos\ on the nilpotent group $G$.
  
    Recall that if $E$ is a finite dimensional vector space the Schwartz class $\cS(E)$ carries a natural Fr\'echet space topology 
    and the Fourier transform of a function $f\in \cS(E)$
    is the element $\hat{f}\in \cS(E^{*})$ such that $\hat{f}(\xi)=\int_{E}e^{i\acou{\xi}{x}}f(x)dx$ for any $\xi \in E^{*}$, where $dx$ denotes the
    Lebesgue measure of $E$.  In the case where $E=(T_{a}M/H_{a})\oplus H_{a}$ the Lebesgue measure actually agrees with the Haar measure of $G$,
    so $\cS(E)$ and $\cS(G)$ agree.  Furthermore, as $E^{*}=(T_{a}M/H_{a})^{*}\otimes H_{a}^{*}$ is just the linear dual $\fg^{*}$ of the Lie algebra of $G$, 
    we also see that $\cS(E^{*})$ agrees with $\cS(\fg^{*})$.
    
      Let $\cS_{0}(G)$ denote the closed subspace of $\cS(G)$ consisting of functions $f \in \cS(G)$ 
  such that for any differential operator $P$ on $\fg^{*}$ we have $(P\hat{f})(0)=0$.
  Notice that the image $\hat{\cS}_{0}(G)$ of $\cS(G)$ under the Fourier transform consists of functions $v\in \cS(\fg^{*})$ such that, given any norm 
  $|.|$ on $G$,  near $\xi=0$ we have $|g(\xi)|=\op{O}(|\xi|^{N})$ for any $N\in \N$.
    
   We endow $\fg^{*}$ with the dilations $\lambda.\xi=(\lambda^{2}\xi_{0},\lambda\xi')$ coming from~(\ref{eq:Heisenberg.Heisenberg-dilations}). For 
   $m\in \C$ we let $S_{m}(\fg^{*}M)$ denote the closed subspace of $C^{\infty}(\fg^{*}\setminus 0)$ consisting in functions 
   $p(\xi)\in C^{\infty}(\fg^{*}\setminus 0)$  such that $p(\lambda.\xi)=\lambda^{m}p(\xi)$ for any $\lambda>0$. 
   
    If $p(\xi)\in S_{m}(\fg^{*})$ then it defines an element of $\hat{\cS}_{0}(\fg^{*})'$ by letting
\begin{equation}
     \acou{p}{g}= \int_{\fg^{*}}p(\xi)g(\xi)d\xi, \qquad g \in \hat{\cS}_{0}(\fg^{*}).
\end{equation}
This allows us to define the inverse Fourier transform of $p$ as the element $\check{p}\in \cS_{0}(G)'$ such that 
$\acou{\check{p}}{f}=\acou{p}{\check{f}}$ for any $f \in \cS_{0}(G)$. It then can be shown (see, e.g., \cite{BG:CHM}, \cite{CGGP:POGD}) 
that the left-convolution with $p$ defines a continuous endomorphism of $\cS_{0}(G)$ via the formula, 
\begin{equation}
    \op{Op}(p)f(x)=\check{p}*f(x)=\acou{\check{p}(y)}{f(xy)}, \qquad f\in \cS_{0}(G).
     \label{Heisenberg.left-invariant-PDO}
\end{equation}
Moreover, we have a bilinear product,
\begin{equation}
    *:S_{m_{1}}(\fg^{*})\times S_{m_{2}}(\fg^{*}) 
\longrightarrow S_{m_{1}+m_{2}}(\fg^{*}),
\label{eq:Heisenberg.product-symbols}
\end{equation}
in such way that, for any $p_{1}\in S_{m_{1}}(\fg^{*})$ and any $p_{2}\in S_{m_{2}}(\fg^{*})$, we have
\begin{equation}
    \op{Op}(p_{1})\circ \op{Op}(p_{2})=\op{Op}(p_{1}*p_{2}).''
\end{equation}

 In addition, if $p \in S_{m}(\fg^{*})$ then  $\op{Op}(p)$ really is a pseudodifferential operator. 
 Indeed, let $X_{0}(a),\ldots,X_{d}(a)$ be a (linear) basis of $\fg$ so that $X_{0}(a)$ is in $T_{a}M/H_{a}$ and 
 $X_{1}(a),\ldots,X_{d}(a)$ span $H_{a}$. For $j=0,\ldots,d$ let $X_{j}^{a}$ be the left-invariant vector field on $G$ such that 
 $X^{a}_{j|_{x=0}}=X_{j}(a)$. The basis $X_{0}(a),\ldots,X_{d}(a)$ yields a linear isomorphism $\fg\simeq \Rd$, hence a global chart of $G$. 
 In the corresponding local coordinates $p(\xi)$ is a homogeneous symbol on $\Rd\setminus 0$ with respect to 
 the dilations~(\ref{eq:Heisenberg.Heisenberg-dilations-Rd}).  
Similarly, each vector field $\frac{1}{i}X_{j}^{a}$, $j=0,\ldots,d$, corresponds to a vector field on $\Rd$ with symbol
 $\sigma_{j}^{a}(x,\xi)$. If we set $\sigma^{a}(x,\xi)=(\sigma_{0}^{a}(x,\xi),\ldots,\sigma_{d}^{a}(x,\xi))$, then
 it can be shown that in these local coordinates we have
 \begin{equation}
     \op{Op}(p)f(x)= (2\pi)^{-(d+1)}\int_{\Rd} e^{i\acou{x}{\xi}}p(\sigma^{a}(x,\xi))\hat{f}(\xi)d\xi, \qquad f \in \cS_{0}(\Rd).
      \label{eq:PsiHDO.PsiDO-convolution}
 \end{equation}
 In other words $\op{Op}(p)$ is the pseudodifferential operator $p(-iX^{a}):=p(\sigma^{a}(x,D))$ acting on $\cS_{0}(\Rd)$.

  \subsection{The \psivdo\ operators}
 The original idea in the Heisenberg calculus, which goes back  to Elias Stein, is to construct a 
 class of operators on a given Heisenberg manifold $(M^{d+1},H)$, called \psivdos, which at any point $a \in M$ are modeled in a suitable sense 
 on left-invariant  pseudodifferential operators on the tangent group $G_{a}M$. 

 Let $U \subset \Rd$ be an open of local coordinates equipped with a $H$-frame $X_{0},\ldots,X_{d}$.

\begin{definition} $S_{m}(\URd)$, $m\in\C$, consists of functions 
    $p(x,\xi)$ in $C^{\infty}(U\times\Rdo)$ which are homogeneous of degree $m$ in the $\xi$-variable with respect to the 
    dilations~(\ref{eq:Heisenberg.Heisenberg-dilations-Rd}), i.e.,  
    we have $p(x,t.\xi)=t^m p(x,\xi)$ for any $t>0$.
\end{definition}

In the sequel we endow $\Rd$ with the pseudo-norm,
\begin{equation}
    \|\xi\|=(\xi_{0}^{2}+\xi_{1}^{4}+\ldots+\xi_{d}^{4})^{1/4}, \qquad \xi\in \Rd.
\end{equation}
In addition, for any multi-order $\beta\in \N^{d+1}_{0}$ we set $\brak\beta=2\beta_{0}+\beta_{1}+\ldots+\beta_{d}$.

\begin{definition} $S^m(\URd)$,  $m\in\C$, consists of functions  $p(x,\xi)$ in $C^{\infty}(\URd)$ with
an asymptotic expansion $ p \sim \sum_{j\geq 0} p_{m-j}$, $p_{k}\in S_{k}(\URd)$, in the sense that, 
for any integer $N$, any compact $K \subset U$ and any multi-orders $\alpha$, $\beta$, there exists $C_{NK\alpha\beta}>0$ such that, 
for any $x\in K$ and any $\xi\in \Rd$ so that $\|\xi \| \geq 1$, we have
\begin{equation}
    | \partial^\alpha_{x}\partial^\beta_{\xi}(p-\sum_{j<N}p_{m-j})(x,\xi)| \leq 
    C_{NK\alpha\beta }\|\xi\|^{\Re m-\brak\beta -N}.
    \label{eq:Heisenberg.asymptotic-expansion-symbols}
\end{equation}
\end{definition}

Next, for $j=0,\ldots,d$ let  $\sigma_{j}(x,\xi)$ denote the symbol (in the 
classical sense) of the vector field $\frac{1}{i}X_{j}$  and set  $\sigma=(\sigma_{0},\ldots,\sigma_{d})$. Then for $p \in S^{m}(\URd)$ we let $p(x,-iX)$ be the 
continuous linear operator from $C^{\infty}_{c}(U)$ to $C^{\infty}(U)$ such that 
    \begin{equation}
          p(x,-iX)f(x)= (2\pi)^{-(d+1)} \int e^{ix.\xi} p(x,\sigma(x,\xi))\hat{f}(\xi)d\xi,
    \qquad f\in C^{\infty}_{c}(U).
    \end{equation}

 In the sequel we let $\psinf(U)$ denote the space of smoothing operators on $U$, that is, the space of continuous operators $P:\cE'(U)\rightarrow \cD'(U)$ with 
a smooth Schwartz kernel. 
 
\begin{definition}
   $\pvdo^{m}(U)$, $m\in \C$, consists of operators $P:C^{\infty}_{c}(U)\rightarrow C^{\infty}(U)$ of the form
\begin{equation}
    P= p(x,-iX)+R, 
\end{equation}
with $p$ in $S^{m}(\URd)$ (called the symbol of $P$) and $R$ smoothing operator.
\end{definition}

 The class of \psivdos\ is invariant under changes of $H$-framed charts (see~\cite[Sect.~16]{BG:CHM}, \cite[Appendix A]{Po:MAMS1}). Therefore, we can 
 extend the definition of \psivdos\ to the Heisenberg manifold $(M^{d+1},H)$ and let them act on sections of a vector bundle $\cE^{r}$ over $M$ as 
 follows. 

\begin{definition}
  $\pvdo^{m}(M,\cE)$, $m\in \C$, consists of continuous operators $P$ from $C^{\infty}_{c}(M,\cE)$ to $C^{\infty}(M,\cE)$ such that:\smallskip
  
  (i) The Schwartz kernel of $P$ is smooth off the diagonal;\smallskip
  
  (ii) For any $H$-framed local chart $\kappa:U\rightarrow V\subset \Rd$ over which there is a trivialization $\tau:\cE_{|U}\rightarrow U\times 
  \C^{r}$ the operator $\kappa_{*}\tau_{*}(P_{|U})$ belongs to $\pvdo^{m}(V,\C^{r}):=\pvdo^{m}(V)\otimes \End \C^{r}$. 
\end{definition}
 
 
\begin{proposition}[\cite{BG:CHM}]  Let $P\in \pvdo^{m}(M,\cE)$, $m \in \C$.\smallskip 
     
     (1)  Let $Q\in \pvdo^{m'}(M,\cE)$, $m'\in \C$, and suppose that $P$ or $Q$ is uniformly properly supported. Then 
     the operator $PQ$ belongs to $\pvdo^{m+m'}(M,\cE)$.\smallskip
 
 (2) The transpose operator $P^{t}$ belongs to $\pvdo^{m}(M,\cE^{*})$.\smallskip 
 
 (3) Suppose that $M$ is endowed with a density~$>0$ and $\cE$ is endowed with a Hermitian metric. Then the adjoint $P^{*}$ of $P$ 
 belongs to $\pvdo^{\overline{m}}(M,\cE)$.      
\end{proposition}

In this setting the principal symbol of a \psivdo\ can be defined intrinsically as follows. 

Let $\fg^{*}M=(TM/H)^{*}\oplus H^{*}$ denote the (linear) dual of the Lie algebra bundle $\fg M$ of $GM$ 
 with canonical projection  $\text{pr}: M\rightarrow \fg^{*}M$. For $m \in \C$ we let $S_{m}(\fg^{*}M,\cE)$ be the space of sections 
$p\in C^{\infty}(\fg^{*}M\setminus 0,\End \text{pr}^{*}\cE)$ such that $p(x,t.\xi)=t^{m}p(x,\xi)$ for any $t>0$. 

\begin{definition}[See {[56]}]\label{def:Heisenberg.principal-symbol} The principal symbol of an operator $P\in \pvdo^{m}(M,\cE)$ is the unique symbol $\sigma_{m}(P)$ in 
    $S_{m}(\fg^{*}M,\cE)$ such that, for any $a\in M$ and for any trivializing $H$-framed local coordinates near $a$, in 
    Heisenberg coordinates centered at $a$ we have
$\sigma_{m}(P)(0,\xi)=p_{m}(0,\xi)$, where $p_{m}(x,\xi)$ is the principal symbol of $P$ in the sense of~(\ref{eq:Heisenberg.asymptotic-expansion-symbols}).
\end{definition}


Given a point $a\in M$ the principal symbol $\sigma_{m}(P)$ allows us to define the model operator of $P$ at $a$ as the left-invariant 
 \psido\ on $\cS_{0}(\fg^{*}M,\cE_{a})$ with symbol $p_{m}^{a}(\xi):=\sigma_{m}(P)(a,\xi)$ so that, in the notation 
 of~(\ref{Heisenberg.left-invariant-PDO}), the operator $P^{a}$ is just $\op{Op}(p_{m}^{a})$.

For $m \in \C$ let $S_{m}(\fg^{*}_{a}M,\cE_{a})$ be the space of functions $p\in C^{\infty}(\fg^{*}_{a}M\setminus 0,\cE_{a})$ which are homogeneous of 
degree $m$. Then the product~(\ref{eq:Heisenberg.product-symbols}) yields a bilinear product, 
\begin{equation}
    *^{a}:S_{m_{1}}(\fg^{*}_{a}M,\cE_{a})\times S_{m_{2}}(\fg^{*}_{a}M,\cE_{a})\rightarrow S_{m_{1}+m_{2}}(\fg^{*}_{a}M,\cE_{a}).
\end{equation}
This product depends smoothly on $a$ as much so to gives rise to the bilinear product, 
\begin{gather}
    *:S_{m_{1}}(\fg^{*}M,\cE)\times S_{m_{2}}(\fg^{*}M,\cE) \longrightarrow S_{m_{1}+m_{2}}(\fg^{*}M,\cE),
    \label{eq:CPCL.product-symbols}\\
    p_{m_{1}}*p_{m_{2}}(a,\xi)=(p_{m_{1}}(a,.)*^{a}p_{m_{2}}(a,.))(\xi), \qquad p_{m_{j}}\in S_{m_{j}}(\fg^{*}M).
\end{gather}
\begin{proposition}[\cite{Po:MAMS1}]\label{prop:Heisenberg.operations-principal-symbols} 
   Let $P\in \pvdo^{m}(M,\cE)$, $m \in \C$. \smallskip
    
        1) Let $Q\in \pvdo^{m'}(M,\cE)$, $m'\in \C$, and suppose that $P$ or $Q$ is uniformly properly supported. Then we have 
        $\sigma_{m+m'}(PQ)=\sigma_{m}(P)*\sigma_{m'}(Q)$,  and for any $a \in M$ the model operator of $PQ$ at $a$ is $P^{a}Q^{a}$.\smallskip
        
        2) We have $\sigma_{m}(P^{t})(x,\xi)=\sigma_{m}(P)(x,-\xi)^{t}$,  and for any $a \in M$ the model operator of $P^{t}$ at $a$  is 
        $(P^{a})^{t}$.\smallskip 
    
    3) Suppose that $M$ is endowed with a density~$>0$ and $\cE$ is endowed with a Hermitian metric.  Then we have 
    $\sigma_{\overline{m}}(P^{*})(x,\xi)=\sigma_{m}(P)(x,\xi)^{*}$,  and for any $a \in M$ the model operator of $P^{*}$ at $a$ is $(P^{a})^{*}$.
\end{proposition}

In addition, there is a complete symbolic calculus 
for \psivdos\ which allows us to carry out the classical parametrix construction  for an operator $P\in \pvdo^{m}(M,\cE)$ whenever its principal 
symbol $\sigma_{m}(P)$ is invertible with respect to the product $*$ (see~\cite{BG:CHM}).   In general, it may be difficult to determine whether $\sigma_{m}(P)$ is invertible with 
respect to that product. 
 Nevertheless, given a point $a\in M$ we have an invertibility criterion for $P^{a}$ in terms of the representation theory of $G_{a}M$; this is 
 the so-called Rockland condition (see, e.g., \cite{Ro:HHGRTC}, \cite{CGGP:POGD}). We then can completely determine the invertibility of the principal 
 symbol of $P$ in terms of the Rockland conditions for its model operators and those of its transpose (see~\cite[Thm.~3.3.19]{Po:MAMS1}).  

Finally, the \psivdos\ enjoy nice Sobolev regularity properties. These properties are best stated in terms of the weighted Sobolev of~\cite{FS:EDdbarbCAHG} 
and~\cite{Po:MAMS1}. These weighted Sobolev spaces can be explicitly related to the usual Sobolev spaces and allows us to show that if $P\in 
\pvdo^{m}(M,\cE)$, $\Re m>0$, has an invertible principal symbol, then $P$ is maximal hypoelliptic, which implies that $P$  
is hypoelliptic with gain of $\frac{m}{2}$-derivatives. We refer to~\cite{BG:CHM} and~\cite{Po:MAMS1} for the precise statements. 
In the sequel we will only need the following. 

\begin{proposition}[\cite{BG:CHM}]\label{prop:Heisenberg.L2-boundedness}
 Assume $M$ compact and let $P\in \pvdo^{m}(M,\cE)$, $\Re m\geq 0$. Then $P$ extends to a bounded operator from $L^{2}(M,\cE)$ to itself and this 
 operator is compact if we further have $\Re m<0$.  
\end{proposition}


\subsection{Holomorphic families of \psivdos}
In this subsection we recall the main definitions and properties of holomorphic families of \psivdos. Throughout the subsection we let 
$(M^{d+1},H)$ be a Heisenberg manifold, we let $\cE^{r}$ be a vector bundle over $M$ and we let $\Omega$ be  an open subset of $\C$.

Let $U\subset \Rd$ be an open of local coordinates equipped with a  $H$-frame $X_{0},\ldots,X_{d}$. We define holomorphic families of symbols 
on $\URd$ as follows. 
\begin{definition}\label{def:Heisenberg.hol-family-symbols}
A family $(p(z))_{z\in\Omega}\subset S^*(\URd)$ is holomorphic when:  \smallskip  
     
     (i) The order $w(z)$ of $p(z)$ depends analytically on $z$; \smallskip 

     (ii) For any $(x,\xi)\in \URd$ the function $z\rightarrow p(z)(x,\xi)$ is holomorphic on $\Omega$; \smallskip 
     
     (iii) The bounds of the asymptotic expansion~(\ref{eq:Heisenberg.asymptotic-expansion-symbols}) 
     for $p(z)$ are locally uniform with respect to $z$, i.e., we have 
     $p(z) \sim \sum_{j\geq 0} p(z)_{ w(z)-j}$, 
         $p(z)_{w(z)-j}\in S_{w(z)-j}(\URd)$, and, for
         any integer $N$, any compacts $K\subset U$ and $L\subset \Omega$ and any multi-orders $\alpha$ and $\beta$, there exists a 
         constant $ C_{NKL\alpha\beta}>0$ such that, for any $(x,z)\in K\times L$ and any $\xi \in \Rd$ so that $\|\xi\|\geq  1$, we have 
\begin{equation}
   | \partial_{x}^\alpha\partial_{\xi}^\beta (p(z)-\sum_{j<N}  
            p(z)_{w(z)-j})(x,\xi)| \leq C_{NKL\alpha\beta} \|\xi\|^{\Re w(z)-N-\brak\beta}.
    \label{eq:Heisenberg.symbols.asymptotic-expansion-hol-families}
\end{equation}
\end{definition}

In the sequel we let $\Hol(\Omega,S^*(\URd))$ denote the class of holomorphic families with  values in $S^*(\URd)$. Notice also that the 
properties (i)--(iii) imply that each homogeneous symbol $p(z)_{w(z)-j}(x,\xi)$ depends analytically on $z$, that is, it gives rise to a 
 holomorphic family with values in $C^{\infty}(\URdo)$ (see~\cite[Rem.~4.2.2]{Po:MAMS1}).   

Since $\psinf(U)=\cL(\cE'(U),C^{\infty}(U))$ is a Fr\'echet space which is isomorphic to $C^{\infty}(U\times U)$ by Schwartz's Kernel 
Theorem, we can define holomorphic families of smoothing operators  as families of operators given by holomorphic families of smooth Schwartz 
kernels. We let $\Hol(\Omega,\psinf(U))$ denote the class of such families. 

\begin{definition}\label{def:Heisenberg.hol-family-PsiHDO's}
 A family $(P(z))_{z\in \Omega}\subset \pvdo^{m}(U)$ is holomorphic when it can be put in the form, 
\begin{equation}
    P(z) = p(z)(x,-iX) + R(z),  \qquad z \in \Omega, 
\end{equation}
 with $(p(z))_{z\in \Omega}\in \Hol(\Omega, S^{*}(\URd))$ and $(R(z))_{z\in \Omega} \in \Hol(\Omega,\psinf(U))$.
\end{definition}

The above notion of holomorphic families of \psivdos\ is invariant under changes of $H$-framed charts (see~\cite{Po:MAMS1}). Therefore, it makes 
sense to define holomorphic families of \psivdos\ on the Heisenberg manifold $(M^{d+1},H)$ acting on sections of the vector bundle $\cE^{r}$  as 
follows. 

\begin{definition}
 A  family $(P(z))_{z\in \Omega}\subset \pvdo^{*}(M,\cE)$ is holomorphic when:\smallskip 
 
     (i) The order $w(z)$ of $P(z)$ is a holomorphic function of $z$;\smallskip 
 
     (ii)  For $\varphi$ and  $\psi$ in $C^\infty_{c}(M)$ with disjoint supports $(\varphi P(z)\psi)_{z\in \Omega}$ is 
     a holomorphic family of smoothing operators;\smallskip 
 
     (iii) For any trivialization $\tau:\cE_{|_{U}}\rightarrow U\times \C^{r}$ over a 
     local $H$-framed chart $\kappa:U \rightarrow V\subset \Rd$ the family $(\kappa_{*}\tau_{*}(P_{z|_{U}}))_{z\in\Omega}$ belongs to 
     $\Hol(\Omega, \pvdo^{*}(V,\C^{r})):=\Hol(\Omega, \pvdo^{*}(V))\otimes \End \C^{r}$.
\end{definition}

We let $\Hol(\Omega,\pvdo^{*}(M,\cE))$ denote the class of holomorphic families of \psivdos\ on $M$ and acting on the sections of $\cE$. 

\begin{proposition}[{\cite[Chap.~4]{Po:MAMS1}}]
    Let $(P(z))_{z\in\Omega}\subset  \pvdo^{*}(M,\cE)$ be a holomorphic family of \psivdos. 
    
    1) Let $(Q(z))_{z\in\Omega}\subset  \pvdo^{*}(M,\cE)$ be a holomorphic family of \psivdos\ and assume that 
    $(P(z))_{z\in\Omega}$ or $(Q(z))_{z\in\Omega}$ is uniformly 
    properly supported with respect to $z$. Then the family $(P(z)Q(z))_{z \in \Omega}$ belongs to $\Hol(\Omega,\pvdo^{*}(M,\cE))$.\smallskip
    
    2) Let $\phi:(M,H)\rightarrow (M',H')$ be a Heisenberg diffeomorphism. Then the family 
    $(\phi_{*}P(z))_{z\in \Omega}$ belongs to $\Hol (\Omega, \Psi_{H'}^{*}(M',\phi_{*}\cE))$. 
\end{proposition}

\subsection{Complex powers of hypoelliptic \psivdos} 
In this subsection we recall the constructions in~\cite{Po:MAMS1} and~\cite{Po:CPDE1} of complex powers of hypoelliptic \psivdos\ as holomorphic 
families of \psivdos. 

Throughout this subsection we let $(M^{d+1},H)$ be a compact Heisenberg manifold equipped with a density~$>0$ and we let $\cE$ be a Hermitian vector 
bundle over $M$. 

Let $P:C^{\infty}(M,\cE)\rightarrow C^{\infty}(M,\cE)$ be a differential operator of Heisenberg order $m$ which is positive, i.e., we have $\acou{Pu}{u}\geq 
0$ for any $u \in C^{\infty}(M,\cE)$, and assume that the principal symbol of $P$ is invertible, that is, $P$ satisfies the Rockland condition at every 
point. 

By standard functional calculus for any $s\in \C$ we can define the power $P^{s}$ as an unbounded operator on $L^{2}(M,\cE)$ whose domain 
contains $C^{\infty}(M,\cE)$. In particular $P^{-1}$ is the partial inverse of $P$ and we have $P^{0}=1-\Pi_{0}(P)$, where $\Pi_{0}(P)$ denotes the 
orthogonal projection onto the kernel of $P$. Furthermore, we have: 

\begin{proposition}[{\cite[Thm.~5.3.4]{Po:MAMS1}}]\label{prop:Heisenberg.complex-powers.positive}
 Assume that $H$ satisfies the bracket condition $H+[H,H]=TM$.   Then the 
 complex powers $(P^{s})_{s \in \C}$ form a holomorphic 1-parameter group of \psivdos\ such that $\ord P^{s}=ms\  \forall s\in \C$.
\end{proposition}

This construction has been generalized to more general hypoelliptic \psivdos\ in~\cite{Po:CPDE1}. 
Let $P:C^{\infty}(M,\cE)\rightarrow C^{\infty}(M,\cE)$ be a \psivdo\ of order $m>0$. In~\cite{Po:CPDE1} there is a notion of \emph{principal cut} for the 
principal symbol $\sigma_{m}(P)$ of $P$ as a ray $L\subset \C\setminus 0$ such that $P-\lambda$ admits a parametrix in a version of the Heisenberg 
calculus with parameter in a conical neighborhood $\Theta \subset \C\setminus 0$ of $L$. 

Let $\Theta(P)$ be the union set of all principal cuts of 
$\sigma_{m}(P)$. Then $\Theta(P)$ is an open conical subset of $\C\setminus 0$ and for any conical subset $\Theta$ of $\Theta(P)$ such that 
$\overline{\Theta}\setminus 0\subset \Theta(P)$ there are at most finitely many eigenvalues of $P$ in $\Theta$ (see~\cite{Po:CPDE1}). 

Let $L_{\theta}=\{\arg \lambda=\theta\}$, $0\leq \theta<2\pi$, be a principal cut for $\sigma_{m}(P)$ such that no eigenvalue of $P$ lies
in $L$. Then $L_{\theta}$ is ray of minimal growth for $P$, so for $\Re s<0$ we define a bounded operator on $L^{2}(M,\cE)$ by letting
 \begin{gather}
   P_{\theta}^{s}= \frac{-1}{2i\pi} \int_{\Gamma_{\theta}} \lambda^{s}_{\theta}(P-\lambda)^{-1}d\lambda,
     \label{eq:Heisenberg.complex-powers-definition}\\
\Gamma_{\theta}=\{ \rho e^{i\theta}; \infty <\rho\leq r\}\cup\{ r e^{it}; 
\theta\geq t\geq \theta-2\pi \}\cup\{ \rho e^{i(\theta-2\pi)};  r\leq \rho\leq \infty\},  
\label{eq:Heisenberg.complex-powers-definition-Gammat}
\end{gather}
where $r>0$ is such that no nonzero eigenvalue of $P$ lies in the disc $|\lambda|<r$. 

\begin{proposition}[\cite{Po:CPDE1}]\label{prop:Heisenberg.powers2}
  The family~(\ref{eq:Heisenberg.complex-powers-definition}) 
 gives rise to a unique holomorphic family  $(P_{\theta}^{s})_{s\in \C}$ of \psivdos\ such that:\smallskip
    
    (i) We have $\ord P_{\theta}^{s}=ms$ for any $s \in \C$;\smallskip
    
    (ii) We have the 1-parameter group property $P_{\theta}^{s_{1}+s_{2}}=P_{\theta}^{s_{1}} P_{\theta}^{s_{2}}$  $\forall s_{j}\in \C$;\smallskip 
     
    (iii) We have $P_{\theta}^{k+s}=P^{k} P_{\theta}^{s}$ for any $k\in \N$ and any $s \in \C$. 
\end{proposition}

Let $E_{0}(P)=\cup_{j \geq 0} \ker P^{j}$ be the characteristic subspace of $P$ associated to $\lambda=0$. This is a finite  dimensional subspace of 
$C^{\infty}(M,\cE)$ and so the projection $\Pi_{0}(P)$ onto $E_{0}(P)$ and along $E_{0}(P^{*})^{\perp}$ 
 is a smoothing operator (see~\cite{Po:CPDE1}). Then we have: 
\begin{equation}
    P_{\theta}^{0}=1-\Pi_{0}(P), \qquad P_{\theta}^{-k}=P^{-k}, \quad k=1,2,\ldots, 
     \label{eq:PsiDO.complex-powers-integers}
\end{equation}
where $P^{-k}$ denotes the partial inverse of $P^{k}$, i.e., the operator that inverts $P^{k}$ on $E_{0}(P^{*})^{\perp}$ and is zero on $E_{0}(P)$.

Assume further that $0$ is not in the spectrum of $P$. Let $Q\in \pvdo^{*}(M,\cE)$ and for $z \in \C$ set $Q(z)=QP_{\theta}^{z/m}$.  Then 
$(Q(z))_{z \in \C}$ is a holomorphic family of \psivdos\ such that $Q_{0}=Q$ and $\ord Q(z)=z+\ord Q$. Following the terminology of~\cite{Gu:GLD} 
a holomorphic family of \psivdos\ with these properties is called a \emph{holomorphic gauging} for $Q$. 

\section{Noncommutative residue trace for the Heisenberg calculus}
\label{sec:NCR}
In this section we construct a noncommutative residue trace for the algebra of integer order \psivdos\ on a Heisenberg manifold. We start by 
describing the logarithmic singularity near the diagonal of the Schwartz kernel of a \psivdo\ of integer order and we show that it gives rise to 
a well-defined density.  We then construct the noncommutative residue for the Heisenberg calculus as the residual trace induced by the analytic 
continuation of the usual trace to \psivdos\ of non-integer orders. Moreover, we show that it agrees with the integral of the density defined by the 
logarithmic singularity of the Schwartz kernel of the corresponding \psivdo. Finally, we prove that when the manifold is connected then every other trace on the algebra 
of integer order \psivdos\ is a constant multiple of our noncommutative residue. This is the analogue of a well-known result of Wodzicki and Guillemin.

\subsection{Logarithmic singularity of the kernel of a \psivdo}
In this subsection we show that the logarithmic singularity of the Schwartz kernel of any integer order \psivdo\  gives rise to a 
density which makes sense intrinsically. This uses the characterization of \psivdos\ in terms of their Schwartz kernels, which we shall now recall. 
 
First, we extend the notion of homogeneity of functions to distributions. For $K$ in $\cS'(\Rd)$ and for $\lambda >0$ we let $K_{\lambda}$ denote 
  the element of $\cS'(\Rd)$ such that
     \begin{equation}
           \acou{K_{\lambda}}{f}=\lambda^{-(d+2)} \acou{K(x)}{f(\lambda^{-1}.x)} \quad \forall f\in\cS(\Rd). 
            \label{eq:PsiHDO.homogeneity-K-m}
      \end{equation}
It will be convenient to also use the notation $K(\lambda.x)$ for denoting $K_{\lambda}(x)$. We say that  
$K$ is homogeneous of degree $m$, $m\in\C$, when $K_{\lambda}=\lambda^m K$ for any $\lambda>0$. 

In the sequel we let $E$ be the anisotropic radial vector field $2x_{0}\partial_{x_{0}}+\partial_{x_{1}}+\ldots+\partial_{x_{d}}$, i.e., $E$ is the 
infinitesimal generator of the flow $\phi_{s}(\xi)=e^{s}.\xi$.

\begin{lemma}[{\cite[Prop.~15.24]{BG:CHM}},~{\cite[Lem.~I.4]{CM:LIFNCG}}]\label{lem:Heisenberg.extension-symbol}
    Let $p(\xi) \in S_{m}(\Rd)$, $m \in \C$.\smallskip
    
    1) If $m$ is not an integer~$\leq -(d+2)$, then $p(\xi)$ can be uniquely extended into a homogeneous distribution $\tau \in \cS'(\Rd)$.\smallskip
    
    2) If $m$ is an integer $\leq -(d+2)$, then at best we can extend $p(\xi)$ into a distribution $\tau \in \cS'(\Rd)$ such that
\begin{equation}
       \tau_{\lambda}=\lambda^{m}\tau +\lambda^{m}\log \lambda\sum_{\brak\alpha=-(m+d+2)} c_{\alpha}(p)\delta^{(\alpha)} \quad \text{for any $\lambda 
       >0$},
     \label{eq:NCR.log-homogeneity}
\end{equation}
  where we have let $c_{\alpha}(p) = \frac{(-1)^{|\alpha|}}{\alpha!}\int_{\|\xi\|=1}\xi^\alpha p(\xi)i_{E}d\xi$. 
  In particular, $p(\xi)$ admits a homogeneous extension if and only if all the coefficients $c_{\alpha}(p)$ vanish. 
\end{lemma}
\begin{remark}
    For reader's convenience a detailed proof of this lemma is given in Appendix. 
\end{remark}

Let $\tau\in \cS'(\Rd)$ and let $\lambda>0$. Then for any $f \in \cS(\Rd)$ we have
\begin{equation}
   \acou{(\check{\tau})_{\lambda}}{f}= \lambda^{-(d+2)}\acou{\tau}{(f_{\lambda^{-1}})^{\vee}}= 
   \acou{\tau}{(\check{f})_{\lambda}}= \lambda^{-(d+2)}\acou{(\tau_{\lambda^{-1}})^{\vee}}{f}.
    \label{eq:NCR.Fourier-transform-scaling}
\end{equation}
Hence $(\check{\tau})_{\lambda} = \lambda^{-(d+2)}(\tau_{\lambda^{-1}})^{\vee}$. Therefore, if we set $\hat{m}=-(m+d+2)$ then we see that:\smallskip

- $\tau$ is homogeneous of degree $m$ if and only if  $\check{\tau}$ is homogeneous of degree $\hat{m}$;\smallskip

- $\tau$ satisfies~(\ref{eq:NCR.log-homogeneity}) if and only if for any $\lambda>0$ we have
\begin{equation}
    \check{\tau}(\lambda.y)= \lambda^{\hat{m}} \check{\tau}(y) - \lambda^{\hat{m}}\log \lambda 
    \sum_{\brak\alpha =\hat{m}} (2\pi)^{-(d+1)}c_{\alpha}(p) (-iy)^{\alpha} .
     \label{eq:NCR.log-homogeneity-kernel}
\end{equation}

Let $U\subset \Rd$ be an open of local coordinates equipped with a  $H$-frame $X_{0},\ldots,X_{d}$.  
In the sequel we set $\N_{0}=\N\cup\{0\}$ and we let $\cS'_{\reg}(\Rd)$  be the space of tempered distributions on $\Rd$ which are smooth outside the 
origin. We endow $\cS'_{\reg}(\Rd)$ with the weakest locally convex topology that makes continuous the embeddings of $\cS'_{\reg}(\Rd)$ into 
$\cS'(\Rd)$ and $C^{\infty}(\Rdo)$. In addition, recall also that if $E$ is a 
topological vector space contained in $\cD'(\Rd)$ then $C^{\infty}(U)\hotimes E$ can be identified as the space
$C^{\infty}(U,E)$ seen as a subspace of $\cD'(\URd)$.

The discussion above about the homogeneity of the (inverse) Fourier transform leads us to consider the classes of distributions below. 
  
\begin{definition}
 $\cK_{m}(\URd)$, $m\in\C$, consists of distributions $K(x,y)$ in $C^\infty(U)\hotimes\cS'_{\reg}(\Rd)$ such that, for any $\lambda>0$, 
 we have:
\begin{equation}
      K(x,\lambda y)= \left\{
      \begin{array}{ll}
          \lambda^m K(x,y) & \text{if $m\not \in \N_{0}$},  \\
          \lambda^m K(x,y) + \lambda^m\log\lambda
                    \sum_{\brak\alpha=m}c_{K,\alpha}(x)y^\alpha & \text{if $m\in \N_{0}$},
      \end{array}\right.
    \label{eq:NCR.log-homogeneity-Km}
\end{equation} 
where the functions $c_{K,\alpha}(x)$, $\brak\alpha=m$, are in $C^{\infty}(U)$ when $m\in \N_{0}$.
\end{definition}

\begin{remark}\label{rem:NCR.regularity-cKm1}
    For $\Re m>0$ we have $K_{m}(\URd)\subset C^{\infty}(U)\hotimes C^{[\frac{\Re m}{2}]'}(\Rd)$, where $[\frac{\Re m}{2}]'$ denotes the greatest 
    integer~$< \Re m$ (see~\cite[Lemma~A.1]{Po:MAMS1}).
\end{remark}

\begin{definition}
$\cK^{m}(\URd)$, $m\in \C$, consists of distributions $K(x,y)$ in $\cD'(\URd)$ with an asymptotic expansion 
     $K\sim \sum_{j\geq0}K_{m+j}$,  $K_{l}\in \cK_{l}(\URd)$,  in the sense that,  for any integer $N$, as soon as $J$ is large enough 
     $K-\sum_{j\leq J}K_{m+j}$ is in  $C^{N}(\URd)$. 
\end{definition}

\begin{remark}\label{rem:NCR.regularity-cKm2}
    The definition implies that any distribution $K\in\cK^{m}(\URd)$ is smooth on $\URdo$. Furthermore, 
    using Remark~\ref{rem:NCR.regularity-cKm1} we see that for
    $\Re m>0$ we have $\cK^{m}(\URd)\subset C^{\infty}(U)\hotimes C^{[\frac{\Re m}{2}]'}(\Rd)$.
\end{remark}

Using Lemma~\ref{lem:Heisenberg.extension-symbol} we can characterize homogeneous symbols on $\URd$ as follows.
\begin{lemma}\label{lem:NCR.extension-symbolU}
Let $m \in \C$ and set $\hat{m}=-(m+d+2)$.\smallskip     
    
 1) If $p(x,\xi)\in S_{m}(\URd)$ then $p(x,\xi)$ can be extended into a distribution $\tau(x,\xi)\in C^{\infty}(U)\hotimes 
 \cS_{\reg}'(\Rd)$ such that $K(x,y):=\check{\tau}_{\xiy}(x,y)$ belongs to $\cK_{\hat{m}}(\URd)$. Furthermore,  if $m$ is an 
 integer~$\leq -(d+2)$ then, using the notation of~(\ref{eq:NCR.log-homogeneity-Km}), 
 we have $c_{K,\alpha}(x)= (2\pi)^{-(d+1)}\int_{\|\xi\|=1}\frac{(i\xi)^{\alpha}}{\alpha!}p(x,\xi)\iota_{E}d\xi$.\smallskip
 
 2) If $K(x,y)\in \cK_{\hat{m}}(\URd)$ then the restriction of $\hat{K}_{\yxi}(x,\xi)$ to $\URdo$ is a symbol in $S_{m}(\URd)$.
\end{lemma}
 
Next, for any $x \in U$ we let $\psi_{x}$ (resp.~$\varepsilon_{x}$) denote the change of variable to the privileged (resp.~Heisenberg) coordinates centered at $x$ 
(cf.~Definitions~\ref{def:Heisenberg-privileged-coordinates} and~\ref{def:Heisenberg-Heisenberg-coordinates}).
 
 Let $p \in S_{m}(\URd)$ and let $k(x,y)\in C^{\infty}(U)\hotimes \cD'(U)$ denote the Schwartz kernel of $p(x,-iX)$, so that 
 $[p(x,-iX)u](x)=\acou{k(x,y)}{u(y)}$ for any $u\in C^{\infty}_{c}(U)$. Then one can check (see, e.g., \cite[p.~54]{Po:MAMS1}) that we have: 
 \begin{equation}
    k(x,y)=|\psi_{x}'| \check{p}_{\xiy}(x,-\psi_{x}(y))=|\varepsilon_{x}'|\check{p}_{\xiy}(x,\phi_{x}(-\varepsilon_{x}(y))).
    \label{eq:Heisenberg.kernel-quantization-symbol-psiy}
\end{equation}
Combining this with Lemma~\ref{lem:NCR.extension-symbolU} leads us to the characterization of \psivdos\ below. 

\begin{proposition}[{\cite[Thms.~15.39, 15.49]{BG:CHM}}, {\cite[Prop.~3.1.16]{Po:MAMS1}}]\label{prop:PsiVDO.characterisation-kernel1}
 Consider a continuous operator $P:C_{c}^\infty(U)\rightarrow C^\infty(U)$ with Schwartz kernel $k_{P}(x,y)$. Let $m \in 
 \C$ and set $\hat{m}=-(m+d+2)$. Then the following are equivalent:\smallskip
 
 (i) $P$ is a \psivdo\ of order $m$.\smallskip
 
  (ii) We can put $k_{P}(x,y)$ in the form,
  \begin{equation}
     k_{P}(x,y)=|\psi_{x}'|K(x,-\psi_{x}(y)) +R(x,y) ,
      \label{eq:PsiHDO.characterization-kernel.privileged}
 \end{equation}
for some  $K\in\cK^{\hat{m}}(\URd)$, $K\sim \sum K_{\hat{m}+j}$, and some $R \in C^{\infty}(U\times U)$.\smallskip 

(iii) We can put $k_{P}(x,y)$ in the form,
  \begin{equation}
     k_{P}(x,y)=|\varepsilon_{x}'|K_{P}(x,-\varepsilon_{x}(y)) +R_{P}(x,y) ,
      \label{eq:PsiHDO.characterization-kernel.Heisenberg}
 \end{equation}
for some $K_{P}\in \cK^{\hat{m}}(\URd)$, $K_{P}\sim \sum K_{P,\hat{m}+j}$, and some $R_{P} \in C^{\infty}(U\times U)$.\smallskip 

Furthermore, if (i)--(iii) hold  then we have $K_{P,l}(x,y)=K_{l}(x,\phi_{x}(y))$ and 
 $P$ has symbol $p\sim \sum_{j\geq 0} p_{m-j}$, where  $p_{m-j}(x,\xi)$ is the restriction to $\URdo$ of $(K_{m+j})^{\wedge}_{\yxi}(x,\xi)$. 
\end{proposition}

%
%
%

Now, let $U\subset \Rd$ be an open of local coordinates equipped with a  $H$-frame $X_{0},X_{1},\ldots,X_{d}$. 
Let $m\in\Z$ and let $K\in \cK^{m}(\URd)$, $K\sim \sum_{j\geq m} K_{j}$. 
Then: 
 
- For $j \leq -1$ the distribution $K_{j}(x,y)$ is homogeneous of degree $j$ with respect to $y$ and is smooth for $y\neq 0$;\smallskip 

- For $j=0$ and $\lambda>0$ we have $K_{0}(x,\lambda.y)=K_{0}(x,y)-c_{K_{0},0}(x)\log \lambda$, which by setting $\lambda=\|y\|^{-1}$ with 
$y\neq 0$ gives 
\begin{equation}
    K_{0}(x,y)=K_{0}(x,\|y\|^{-1}.y)-c_{K_{0},0}\log \|y\|.
     \label{eq:Log.behavior-K0}
\end{equation}

- The remainder term $K-\sum_{j\geq 1}K_{j}$ is in $C^{0}(\URd)$ (cf.~Remarks~\ref{rem:NCR.regularity-cKm1} and~\ref{rem:NCR.regularity-cKm2}).\smallskip

\noindent It follows that $K(x,y)$ has a behavior near $y=0$ of the form, 
\begin{equation}
    K(x,y)=\sum_{m\leq j\leq -1} K_{j}(x,y)-c_{K}(x)\log \|y\| +\op{O}(1), \qquad c_{K}(x)=c_{K_{0},0}(x).
     \label{eq:Log.behavior-K}
\end{equation}
%
%

\begin{lemma} Let $P\in \pvdo^{m}(U)$ have kernel $k_{P}(x,y)$ and set $\hat{m}=-(m+d+2)$.\smallskip
    
  1) Near the diagonal $k_{P}(x,y)$ has a behavior of the form,
    \begin{equation}
        k_{P}(x,y)=\sum_{\hat{m}\leq j \leq -1}a_{j}(x,-\psi_{x}(y)) -c_{P}(x) \log \|\psi_{x}(y)\| +\op{O}(1),
         \label{eq:Log.behavior-kP}
\end{equation}
with $a_{j}(x,y)\in C^{\infty}(\URdo)$ homogeneous of degree $j$ in $y$ and $c_{P}(x)\in C^{\infty}(U)$.\smallskip

2) If  we write $k_{P}(x,y)$ in the forms~(\ref{eq:PsiHDO.characterization-kernel.privileged}) 
and~(\ref{eq:PsiHDO.characterization-kernel.Heisenberg}) with $K(x,y)$ and $K_{P}(x,y)$ in $\cK^{\hat{m}}(\URd)$, then we have
\begin{equation}
    c_{P}(x)=|\psi_{x}'|c_{K}(x)=|\varepsilon_{x}'|c_{K_{P}}(x)=\frac{|\psi_{x}'|}{(2\pi)^{d+1}}\int_{\|\xi\|=1}p_{-(d+2)}(x,\xi)\imath_{E}d\xi,
     \label{eq:NCR.formula-cP}
\end{equation}
where $p_{-(d+2)}$ denotes the symbol of degree $-(d+2)$ of $P$. 
\end{lemma}
\begin{proof}
If we put $k_{P}(x,y)$ in the form~(\ref{eq:PsiHDO.characterization-kernel.privileged}) with $K\in \cK^{\hat{m}}(\URd)$, $K\sim \sum K_{\hat{m}+j}$, 
then it follows from~(\ref{eq:Log.behavior-K}) that $k_{P}(x,y)$ has a behavior 
near the diagonal of the form~(\ref{eq:Log.behavior-kP}) with $c_{P}(x)=|\psi_{x}'|c_{K}(x)=|\psi_{x}'|c_{K_{0},0}(x)$. Furthermore, by 
Proposition~\ref{prop:PsiVDO.characterisation-kernel1} 
the symbol $p_{-(d+2)}(x,\xi)$ of degree $-(d+2)$ of $P$ is the restriction to $\URdo$ of $(K_{0})^{\wedge}_{\yxi}(x,\xi)$, so by 
Lemma~\ref{lem:NCR.extension-symbolU} 
we have $c_{K}(x)=c_{K_{0},0}(x)=(2\pi)^{-(d+1)}\int_{\|\xi\|=1}p_{-(d+2)}(x,\xi)\imath_{E}d\xi$. 

Next, if we put  $k_{P}(x,y)$ in the form~(\ref{eq:PsiHDO.characterization-kernel.Heisenberg}) with 
$K_{P}\in \cK^{\hat{m}}(\URd)$, $K_{P}\sim \sum K_{P,\hat{m+j}}$ then by Proposition~\ref{prop:PsiVDO.characterisation-kernel1} 
we have $K_{P,0}(x,y)=K_{0}(x,\phi_{x}(y))$. Let 
$\lambda>0$. Since $\phi_{x}(\lambda.y)=\lambda.\phi_{x}(y)$, using~(\ref{eq:NCR.log-homogeneity-Km}) we get
\begin{equation}
    K_{P,0}(x,\lambda.y)-K_{P,0}(x,y)=  K_{0}(x,\lambda.\phi_{x}(y))-K_{0}(x,\phi_{x}(y))=c_{K_{0}}(x) \log \lambda . \\
\end{equation}
Hence $c_{K_{P,0}}(x)=c_{K,0}(x)$. As $|\varepsilon_{x}'|=|\phi_{x}'|.|\psi_{x}'|=|\psi_{x}'|$ we see that  
$|\psi_{x}'|c_{K}(x)=|\varepsilon_{x}'|c_{K_{P}}(x)$. The proof is thus achieved.
\end{proof}

 \begin{lemma}\label{lem:log-sing.invariance}
Let $\phi : U\rightarrow \tilde{U}$ be a change of $H$-framed local coordinates. Then for any $\tilde{P}\in 
  \pvdo^{m}(\tilde{U})$ we have $c_{\phi^{*}\tilde{P}}(x)=|\phi'(x)|c_{\tilde{P}}(\phi(x))$. 
\end{lemma}
\begin{proof}
Let $P=\phi^{*}\tilde{P}$. Then $P$ is a \psivdo\ of order $m$ on $U$ (see~\cite{BG:CHM}). Moreover, by~\cite[Prop.~3.1.18]{Po:MAMS1}
if we write the Schwartz kernel $k_{\tilde{P}}(\tilde{x},\tilde{y})$ in the 
form~(\ref{eq:PsiHDO.characterization-kernel.Heisenberg}) with $K_{\tilde{P}}(\tilde{x},\tilde{y})$ in $\cK^{\hat{m}}(\tilde{U}\times \Rd)$, 
then the Schwartz kernel $k_{P}(x,y)$ of $P$ can be put in 
the form~(\ref{eq:PsiHDO.characterization-kernel.Heisenberg}) with $K_{P}(x,y)$ in $\cK^{\hat{m}}(U\times \Rd)$  such that
   \begin{equation}
       K_{P}(x,y) \sim \sum_{\brak\beta\geq \frac{3}{2}\brak\alpha} \frac{1}{\alpha!\beta!} 
       a_{\alpha\beta}(x)y^{\beta}(\partial_{\tilde{y}}^{\alpha}K_{\tilde{P}})(\phi(x),\phi_{H}'(x).y),
        \label{eq:PsiHDO.asymptotic-expansion-KP}
   \end{equation}
   where we have let $a_{\alpha\beta}(x)=\partial^{\beta}_{y}[|\partial_{y}(\varepsilon_{\phi(x)}\circ \phi\circ \tilde{\varepsilon}_{x}^{-1})(y)| 
   (\tilde{\varepsilon}_{\phi(x)}\circ \phi\circ \varepsilon_{x}^{-1}(y)-\phi_{H}'(x)y)^{\alpha}]_{|_{y=0}}$, the map $\phi_{H}'(x)$ is the tangent 
   map~(\ref{eq:Heisenberg.tangent-map}), and $\tilde{\varepsilon}_{\tilde{x}}$ 
   denotes the change to the Heisenberg coordinates at $\tilde{x}\in \tilde{U}$. In particular, we have
\begin{equation}
 K_{P}(x,y)= a_{00}(x) K_{\tilde{P}}(\phi(x),\phi_{H}'(x).y) \quad 
           \bmod y_{j}\cK^{\hat{m}+1}(\URd),
     \label{eq:PsiHDO.asymptotic-expansion-KP2}
\end{equation}
where $a_{00}(x)=|\varepsilon_{\phi(x)}'||\phi'(x)||\varepsilon_{x}'|^{-1}$. 

Notice that $\tilde{K}(x,y):= K_{\tilde{P}}(\phi(x),\phi_{H}'(x).y)$ is an element of 
   $\cK^{\hat{m}}(\URd)$, since we have $\phi'_{H}(x).(\lambda.y)=\lambda.(\phi'_{H}(x).y)$ for any $\lambda>0$.
  Moreover, the distributions in $y_{j}\cK^{*}(\URd)$, $j=0,..,d$, cannot have a logarithmic singularity near $y=0$. To see this it is enough to look 
  at a distribution $H(x,y)\in \cK^{-l}(\URd)$, $l\in \N_{0}$. Then $H(x,y)$ has a behavior near $y=0$ of the form:
\begin{equation}
    H(x,y)=\sum_{-l\leq k \leq -1}b_{k}(x,y)-c_{H}(x)\log \|y\|+\op{O}(1),
\end{equation} 
with $b_{k}(x,y)$ homogeneous of degree $k$ with respect to the $y$-variable. Thus, 
\begin{equation}
    y_{j}H(x,y)=\sum_{-l\leq k \leq -1}y_{j}b_{k}(x,y)-c_{H}(x)y_{j}\log \|y\|+\op{O}(1).
\end{equation}
Observe that each term $y_{j}b_{k}(x,y)$ is homogeneous of degree $k+1$ with respect to $y$ and the term $y_{j}\log \|y\|$ converges to $0$ as $y\rightarrow 0$. 
Therefore, we see that the singularity of $y_{j}H(x,y)$ near $y=0$ cannot contain a logarithmic term. 
  
Combining the above observations with~(\ref{eq:PsiHDO.asymptotic-expansion-KP}) shows that   
  the coefficients of the logarithmic singularities of $K_{P}(x,y)$ and $a_{00}(x)\tilde{K}(x,y)$ must agree, i.e., 
   we have $c_{K_{P}}(x)=c_{a_{00}\tilde{K}}(x)=a_{00}(x)c_{\tilde{K}}(x)=|\varepsilon_{\phi(x)}'||\phi'(x)||\varepsilon_{x}'|^{-1}c_{\tilde{K}}(x)$. 
 Furthermore, the only contribution to the logarithmic singularity of $\tilde{K}(x,y)$ comes from
 \begin{multline}
     c_{K_{\tilde{P}}}(\phi(x))\log\|\phi_{H}'(x)y\|= c_{K_{\tilde{P}}}(\phi(x))\log[\|y\| 
     \|\phi_{H}'(x).(\|y\|^{-1}.y\|)] \\ =c_{K_{\tilde{P}}}(\phi(x))\log\|y\| +\op{O}(1).
 \end{multline}
   Hence $c_{\tilde{K}}(x)=c_{K_{\tilde{P}}}(\phi(x))$. Therefore, we get 
   $c_{K_{P}}(x)=|\varepsilon_{\phi(x)}'||\phi'(x)||\varepsilon_{x}'|^{-1}c_{K_{\tilde{P}}}(\phi(x))$,
   which by combining with~(\ref{eq:NCR.formula-cP}) shows that $c_{P}(x)=|\phi'(x)|c_{\tilde{P}}(\phi(x))$ as desired.  
\end{proof}
   
Let $P\in \pvdo^{m}(M,\cE)$ and let $\kappa:U\rightarrow V$ be a $H$-framed chart over which there is a trivialization $\tau:\cE_{|_{U}}\rightarrow U\times \C^{r}$. 
Then the Schwartz kernel of $P_{\kappa,\tau}:=\kappa_{*}\tau_{*}(P_{|_{U}})$ has a singularity near the diagonal of the form~(\ref{eq:Log.behavior-kP}). 
Moreover, if 
$\tilde{\kappa}:\tilde{U}\rightarrow \tilde{V}$ be a $H$-framed chart over which there is a trivialization $\tau:\cE_{|_{\tilde{U}}}\rightarrow 
\tilde{U}\times \C^{r}$ and if we let $\phi$ denote the Heisenberg diffeomorphism $\tilde{\kappa}\circ \kappa^{-1}:\kappa(U\cap \tilde{U})\rightarrow 
\tilde{\kappa}(U\cap \tilde{U})$, then by Lemma~\ref{lem:log-sing.invariance} 
we have $c_{P_{\kappa,\tau}}(x)=|\phi'(x)|c_{P_{\tilde{\kappa},\tilde{\tau}}}(\phi(x))$ for 
any $x \in U$. Therefore, on $U\cap \tilde{U}$ we have the equality of densities, 
\begin{equation}
    \tau^{*}\kappa^{*}(c_{P_{\kappa,\tau}}(x)dx)= \tilde{\tau}^{*}\tilde{\kappa}^{*}(c_{P_{\tilde{\kappa},\tilde{\tau}}}(x)dx).
\end{equation}

Now, the space  $C^{\infty}(M,|\Lambda|(M)\otimes \End \cE)$ of $\END \cE$-valued densities is a sheaf, so 
there exists a unique density $c_{P}(x) \in C^{\infty}(M,|\Lambda|(M)\otimes \End \cE)$ such that, for any local 
$H$-framed chart $\kappa:U\rightarrow V$ and any trivialization $\tau:\cE_{|_{U}}\rightarrow U\times \C^{r}$, we have 
\begin{equation}
    c_{P}(x)|_{U}=\tau^{*}\kappa^{*}(c_{\kappa_{*}\tau_{*}(P_{|_{U}})}(x)dx).
\end{equation}
Moreover, this density is functorial with respect to Heisenberg diffeomorphisms, i.e., for any Heisenberg diffeomorphism $\phi:(M,H)\rightarrow 
(M',H')$ we have 
\begin{equation}
    c_{\phi_{*}P}(x)=\phi_{*}(c_{P}(x)).
     \label{eq:Log.functoriality-cP}
\end{equation}

Summarizing all this we have proved:

\begin{proposition}\label{thm:NCR.log-singularity}
    Let $P\in \pvdo^{m}(M,\cE)$, $m \in \Z$. Then:\smallskip 
    
    1) On  any trivializing $H$-framed local coordinates the Schwartz kernel $k_{P}(x,y)$ of $P$ 
    has a behavior near the diagonal  of the form, 
    \begin{equation}
        k_{P}(x,y)=\sum_{-(m+d+2)\leq j\leq -1}a_{j}(x,-\psi_{x}(y)) - c_{P}(x)\log \|\psi_{x}(y)\| + 
        \op{O}(1), 
    \end{equation}
    where $c_{P}(x)$ is  given by~(\ref{eq:NCR.formula-cP}) and each function 
    $a_{j}(x,y)$ is smooth for $y\neq 0$ and homogeneous of degree $j$ with respect to $y$.\smallskip 
    
    2) The coefficient $c_{P}(x)$ makes sense globally on $M$ as a smooth $\END\cE$-valued density which is functorial with respect to
    Heisenberg diffeomorphisms.
\end{proposition}

Finally, the following holds.

\begin{proposition}\label{prop.Sing.transpose-adjoint}
      Let $P\in \pvdo^{m}(M,\cE)$, $m \in \Z$.\smallskip 
     
     1) Let $P^{t}\in \pvdo^{m}(M,\cE^{*})$ be the transpose of $P$. Then we have $c_{P^{t}}(x)=c_{P}(x)^{t}$.\smallskip
     
     2) Suppose that $M$ is endowed with a density $\rho>0$ and $\cE$ is endowed with a Hermitian metric. Let $P^{*}\in \pvdo^{m}(M,\cE)$ be the 
    adjoint of $P$. Then we have $c_{P^{*}}(x)=c_{P}(x)^{*}$. 
\end{proposition}
\begin{proof}
   Let us first assume that $\cE$ is the trivial line bundle. Then it is enough to prove the result in $H$-framed local coordinates $U\subset \Rd$, so 
   that the Schwartz kernel $k_{P}(x,y)$ can be put in the form~(\ref{eq:PsiHDO.characterization-kernel.Heisenberg}) 
   with $K_{P}(x,y)$ in $\cK^{\hat{m}}(\URd)$. 
   
   We know that $P^{t}$ is a \psivdo\ 
   of order $m$ (see~\cite[Thm.~17.4]{BG:CHM}). Moreover, by~\cite[Prop.~3.1.21]{Po:MAMS1} 
   we can put its Schwartz kernel $k_{P^{t}}(x,y)$ in the form~(\ref{eq:PsiHDO.characterization-kernel.Heisenberg}) 
   with $K_{P^{t}}(x,y)$ in $\cK^{\hat{m}}(\URd)$ such that
   \begin{equation}
       K_{P^{t}}(x,y) \sim  \sum_{\frac{3}{2}\brak\alpha \leq \brak \beta} \sum_{|\gamma|\leq |\delta| \leq 2|\gamma|} 
       a_{\alpha\beta\gamma\delta}(x) y^{\beta+\delta}  
       (\partial^{\gamma}_{x}\partial_{y}^{\alpha}K_{P})(x,-y), 
   \end{equation}
   where 
   $a_{\alpha\beta\gamma\delta}(x)=\frac{|\varepsilon_{x}^{-1}|}{\alpha!\beta!\gamma!\delta!}
   [\partial_{y}^{\beta}(|\varepsilon_{\varepsilon_{x}^{-1}(-y)}'|(y-\varepsilon_{\varepsilon_{x}^{-1}(y)}(x))^{\alpha})
   \partial_{y}^{\delta}(\varepsilon_{x}^{-1}(-y)-x)^{\gamma}](x,0)$. In particular, we have 
   $K_{P^{t}}(x,y)=K_{P}(x,-y) \bmod y_{j}\cK^{\hat{m}+1}(\URd)$. Therefore, in the same way as in the proof of Lemma~\ref{lem:log-sing.invariance}, we see 
   that the logarithmic singularity near $y=0$ of $K_{P}(x,y)$ agrees with that of $K_{P^{t}}(x,-y)$, hence with that of $K_{P^{t}}(x,y)$. Therefore, we 
   have $c_{K_{P^{t}}}(x)=c_{K_{P}}(x)$. Combining this 
   with~(\ref{eq:NCR.formula-cP}) then shows that $c_{P^{t}}(x)=c_{P}(x)$. 
   
   Next, suppose that $U$ is endowed with a smooth density $\rho(x)>0$. Then the corresponding adjoint $P^{*}$ is a \psivdo\ of order $m$ on $U$ with 
   Schwartz kernel $k_{P^{*}}(x,y)=\rho(x)^{-1}\overline{k_{P^{t}}(x,y)}\rho(y)$. Thus $k_{P^{*}}(x,y)$ can be put in the 
   form~(\ref{eq:PsiHDO.characterization-kernel.Heisenberg})  with 
   $K_{P^{*}}(x,y)$ in $\cK^{\hat{m}}(\URd)$ such that
   \begin{multline}
       K_{P^{*}}(x,y)=[\rho(x)^{-1}\rho(\varepsilon_{x}^{-1}(-y))]\overline{K_{P^{t}}(x,y)}\\
       =\overline{K_{P^{t}}(x,y)} \  \bmod y_{j}\cK^{\hat{m}+1}(\URd).
   \end{multline}
   Therefore, $K_{P^{*}}(x,y)$ and $\overline{K_{P^{t}}(x,y)}$ same logarithmic singularity near $y=0$, so that we have 
   $c_{K_{P^{*}}}(x)=\overline{c_{K_{P^{t}}}(x)}=\overline{c_{K_{P}}(x)}$. Hence $c_{P^{*}}(x)=\overline{c_{P}(x)}$. 
   
   Finally, when $\cE$ is a general vector bundle, we can argue as above to show that we still have $c_{P^{t}}(x)=c_{P}(x)^{t}$, and if  $P^{*}$ is
   the adjoint of $P$ with respect to the density $\rho$ and some Hermitian metric on $\cE$, then we have 
   $c_{P^{*}}(x)=c_{P}(x)^{*}$. 
\end{proof}

\subsection{Noncommutative residue}
Let $(M^{d+1},H)$ be a Heisenberg manifold and let $\cE$ be a vector bundle over $M$. We shall now construct a noncommutative residue trace on 
the algebra $\pvdoz(M,\cE)$ as the residual trace induced by the analytic extension of the operator trace to  \psivdos\ of non-integer order.

Let $\pvdoi(M,\cE) := \cup_{\Re m < -(d+2)}\pvdo^{m}(M,\cE)$ the class of \psivdos\ whose symbols are integrable with respect to the $\xi$-variable 
(this notation is borrowed from~\cite{CM:LIFNCG}). If $P$ belongs to this class, then it follows from Remark~\ref{rem:NCR.regularity-cKm2} 
 that the restriction  to the diagonal of $M\times M$ of its Schwartz kernel  defines a smooth density $k_{P}(x,x)$ with values in $\End \cE$. 
 Therefore, when $M$ is compact then $P$ is a trace-class operator on $L^{2}(M,\cE)$ and we have 
\begin{equation}
    \Tra (P) = \int_{M} \tr_{\cE}k_{P}(x,x).
\end{equation}

We shall now construct an analytic extension of the operator trace to the class $\pvdocz(M,\cE)$ of \psivdos\ of non-integer order. 
As in~ \cite{Gu:GLD} (see also \cite{KV:GDEO}, \cite{CM:LIFNCG}) the approach consists in working directly 
at the level of densities by constructing an analytic extension
of the map $P\rightarrow k_{P}(x,x)$ to $\pvdocz(M,\cE)$. 
Here analyticity is meant with respect to holomorphic families of \psivdos, e.g., 
the map $P \rightarrow k_{P}(x,x)$ is analytic since for any holomorphic family $(P(z))_{z\in \Omega}$ with values in $\pvdoi(M,\cE)$ the output densities 
$k_{P(z)}(x,x)$ depend analytically on $z$ in the Fr\'echet space $C^{\infty}(M,|\Lambda|(M)\otimes \End \cE)$.

Let $U \subset \Rd$ be an open of trivializing local coordinates equipped with equipped with a $H$-frame $X_{0}, \ldots, X_{d}$, and for any $x\in U$ 
let $\psi_{x}$ denote the affine change of variables to the privileged coordinates at $x$. Any $P \in \pvdo^{m}(U)$ can be written as 
$P=p(x,-iX)+R$ with $p\in S^{m}(\URd)$ and $R\in \Psi^{-\infty}(U)$. Therefore, if $\Re m <-(d+2)$ then 
using~(\ref{eq:Heisenberg.kernel-quantization-symbol-psiy}) we get 
\begin{equation}
    k_{P}(x,x)=|\psi_{x}'| (2\pi)^{-(d+2)}\int p(x,\xi)d\xi +k_{R}(x,x). 
    \label{eq:NCR.kP(x,x)}
\end{equation}
This leads us to consider  the functional, 
\begin{equation}
      L(p):=(2\pi)^{-(d+2)}\int p(\xi) d\xi, \qquad p\in\Shi(\Rd). 
\end{equation}

In the sequel, as in Section~\ref{sec:Heisenberg-calculus} for \psivdos, 
we say that a holomorphic family of symbols $(p(z))_{z\in\C}\subset S^{*}(\Rd)$ is a \emph{gauging} for a given symbol $p\in 
S^{*}(\R^{d+1})$ when we have $p(0)=p$ and $\ord p(z)=z+\ord p$ for any $z\in \C$. 

\begin{lemma}[{\cite[Prop. I.4]{CM:LIFNCG}}]\label{lem:Spectral.tildeL}
    1) The functional $L$ has a unique analytic continuation $\tilde{L}$ to $S^{\CZ}(\Rd)$. The value of $\tilde{L}$ on 
    a symbol $p\sim \sum_{j\geq 0} p_{m-j}$ of order $m \in\CZ$ is given by 
\begin{equation}
    \tilde{L}(p)= (p- \sum_{j\leq N}\tau_{m-j})^{\vee}(0),  \qquad N\geq  \Re{m}+d+2,
     \label{eq:NCR.L-tilda}
\end{equation}
    where the value of the integer $N$ is irrelevant and the distribution
    $\tau_{m-j}\in \cS'(\Rd)$ is the unique homogeneous extension of $p_{m-j}(\xi)$ provided by 
    Lemma~\ref{lem:Heisenberg.extension-symbol}.\smallskip 
    
       2) Let $ p\in S^{\Z}(\Rd)$, $p\sim \sum p_{m-j}$, and let $(p(z))_{z \in \C}\subset S^{*}(\Rd)$ 
       be a  holomorphic gauging for $p$. Then $\tilde{L}(p(z))$ has at worst a simple 
        pole singularity at $z=0$ in such way that 
        \begin{equation}
            \res_{z=0}\tilde{L}(p(z)) =  \int_{\|\xi\|=1} p_{-(d+2)}(\xi) \imath_{E}d\xi, 
        \end{equation}
        where  $p_{-(d+2)}(\xi)$ is the symbol of degree $-(d+2)$ of $p(\xi)$ and $E$ is the anisotropic radial vector field  
        $2\xi_{0}\partial_{x_{0}}+\xi_{1}\partial_{\xi_{1}}+\ldots+\xi_{d}\partial_{\xi_{d}}$.
\end{lemma}
\begin{proof}
 First, the extension is necessarily unique since the functional $L$ is holomorphic on $\Shi(\Rd)$ and each symbol 
 $p\in S^{\CZ}(\Rd)$ can be  connected to $\Shi(\Rd)$ by means of a holomorphic family with values in $S^{\CZ}(\Rd)$. 
 
 Let  $p\in S^{\CZ}(\Rd)$, $p\sim \sum_{j\geq 0} p_{m-j}$, and for $j=0,1,\ldots$ let $\tau_{m-j} \in \cS'(\Rd)$ denote the unique 
 homogeneous extension of $p_{m-j}$ provided by Lemma~\ref{lem:Heisenberg.extension-symbol}. 
 For $N\geq  \Re{m}+d+2$
the distribution $p-\sum_{j\leq N}\tau_{m-j}$ agrees with an integrable function near $\infty$, so its Fourier transform is continuous and we may define
\begin{equation}
    \tilde{L}(p)= (p-\sum_{j\leq N}\tau_{m-j})^\wedge (0).
    \label{eq:Spectral.tildeL}
\end{equation}
Notice that if $j>\Re m +d+2$ then $\tau_{m-j}$ is also integrable 
near $\infty$, so  $\hat\tau_{m-j}(0)$ makes well sense. However, its value must be $0$ for homogeneity reasons. This shows that the value of 
$N$  in~(\ref{eq:Spectral.tildeL}) is irrelevant, so this formula defines a linear functional on $S^{\CZ}(\Rd)$. In particular, if $\Re m<-(d+2)$ then 
we can take $N=0$ to get $\tilde{L}(p)=\check{p}(0)=\int p(\xi)d\xi=L(p)$. Hence 
$\tilde{L}$ agrees with $L$ on $\Shi(\Rd)\cap S^{\CZ}(\Rd)$. 

Let $(p(z))_{z \in \Omega}$ be a holomorphic family of symbols such that $w(z)=\ord p(z)$ is never an integer 
and let us study the analyticity of $\tilde{L}(p(z))$. As the functional $L$ is holomorphic on $\Shi(\Rd)$ we may assume that we have
$|\Re w(z)-m|<1$ for some integer $m\geq -(d+2)$.  In this case in~(\ref{eq:Spectral.tildeL}) we can set $N=m+d+2$  and for 
$j=0,\ldots,m+d+1$ we can represent $\tau(z)_{w(z)-j}$ by $p(z)_{w(z)-j}$. Then, 
picking $\varphi \in C_{c}^\infty(\Rd)$ such that $\varphi=1$ near the origin, we see that $ \tilde{L}(p(z))$ is equal to
\begin{multline}
     \int [p(z)(\xi)- (1-\varphi(\xi) )\sum_{j\leq m+d+2} p(z)_{w(z)-j}(\xi)] d\xi 
                             -  \sum_{j\leq m+d+2} \acou{\tau (z)_{w(z)-j}}{\varphi} \\
       = L(\tilde{p}(z)) -\acou{\tau(z)}{\varphi}- \sum_{j\leq m+d+1} \int p(z)_{ w(z)-j}(\xi)\varphi(\xi)d\xi ,  
    \label{eq:Spectral.tildeLbis} 
\end{multline}
where we have let $\tau(z)= \tau(z)_{w(z)-m-(d+2)}$ and 
$\tilde{p}(z)=p(z)-(1-\varphi)\sum_{j\leq m+d+2} p(z)_{w(z)-j}$. 

In the r.h.s.~of~(\ref{eq:Spectral.tildeLbis}) the only term that may 
fail to be analytic is $- \acou{\tau(z)}{\varphi}$. Notice that by the formulas~(\ref{eq:Appendix.almosthomogeneous-extension}) 
and~(\ref{eq:Appendix1.h'}) in Appendix we have
\begin{equation}
    \acou{\tau(z)}{\varphi}= \int p(z)_{w(z)-m-(d+2)}(\varphi(\xi)-\psi_{z}(\xi))d\xi, 
\end{equation}
with $\psi_{z}\in C^{\infty}(\Rd)$ of the form $\psi_{z}(\xi)= \int^{\infty}_{\log\|\xi\|} 
[(\frac 1{w(z)-m}\frac{d}{ds} +1)g](t)dt$, where $g(t)$ can be any function in $C^{\infty}_{c}(\R)$  such that 
$\int g(t)dt=1$. Without any loss of generality we may suppose that $\varphi(\xi)= \int^{\infty}_{\log\|\xi\|} 
g(t)dt$ with $g\in C^{\infty}_{c}(\R)$ as above. Then  we have $\psi_{z}(\xi)= -\frac 1{w(z)-m} g(\log \|\xi\|) + \varphi(\xi)$, which gives
\begin{multline}
    \acou{\tau(z)}{\varphi}  =\frac 1{w(z)-m} \int p(z)_{w(z)-m-(d+2)}(\xi) g(\log \|\xi\|)d\xi \\ 
     =\frac1{w(z)-m} 
    \int \mu^{w(z)-m}g(\log \mu)\frac{d\mu}{\mu}\ \int_{\|\xi\|=1}  p(z)_{w(z)-m-(d+2)}(\xi) \imath_{E}d\xi. 
    \label{eq:NCR.tau-z}
\end{multline}
Together with~(\ref{eq:Spectral.tildeLbis}) this shows that $\tilde{L}(p(z))$ is an analytic function, so the
the first part of the lemma is proved.

Finally, let $p\sim \sum p_{m-j}$ be  a symbol in $S^{\Z}(\Rd)$ and let $(p(z))_{|\Re z -m|<1}$ be a  holomorphic family which is a gauging for  
       $p$. Since $p(z)$ has order $w(z)=m+z$ it follows from (\ref{eq:Spectral.tildeLbis}) and (\ref{eq:NCR.tau-z}) that  
       $\tilde{L}(p(z))$ has at worst a simple pole singularity at $z=0$ such that
\begin{multline}
   \res_{z=0}\tilde{L}(p(z))= \res_{z=0} \frac{-1}{z} \int \mu^{z}g(\log \mu)\frac{d\mu}{\mu}\ \int_{\|\xi\|=1}  
    p(z)_{z-(d+2)}(\xi) \imath_{E}d\xi \\ = - \int_{\|\xi\|=1}  p_{-(d+2)}(\xi) \imath_{E}d\xi .
\end{multline}
This proves the second part of the lemma. 
\end{proof} 

Now,  for $P\in \pvdocz(U)$ we let 
\begin{equation}
    t_{P}(x) = (2\pi)^{-(d+2)}|\psi_{x}'| \tilde{L}(p(x,.)) + k_{R}(x,x),
    \label{eq:NCR.tP-definition}
\end{equation}
where the pair $(p,R) \in S^{\CZ}(\URd)\times \Psi^{\infty}(U)$ is such that $P=p(x,-iX)+R$. This definition does not depend on 
the choice of $(p,R)$. Indeed, if $(p',R')$ is another such pair then $p-p'$ is in 
$S^{-\infty}(\URd)$, so using~(\ref{eq:NCR.kP(x,x)}) we that $k_{R'}(x,x)-k_{R}(x,x)$ is equal to
\begin{multline}
 k_{(p-p')(x,-iX)}(x,x) = (2\pi)^{-(d+2)}|\psi_{x}'| L((p-p')(x,.))\\ = 
         (2\pi)^{-(d+2)}|\psi_{x}'| (\tilde{L}(p(x,.)) -\tilde{L}(p'(x,.))),  
\end{multline}
which shows that the r.h.s.~of~(\ref{eq:NCR.tP-definition}) is the same for both pairs. 

On the other hand, observe that~(\ref{eq:Spectral.tildeLbis}) and (\ref{eq:NCR.tau-z}) show that $\tilde{L}(p(x,.))$ depends 
smoothly on $x$ and that for any  holomorphic family $(p(z))(z)\in \Omega \subset  S^{\CZ}(\URd)$ the map 
$z\rightarrow \tilde{L}(p(x,.))$ is holomorphic from $\Omega$ to $C^{\infty}(U)$. Therefore, the map $P\rightarrow t_{P}(x)$ is an analytic extension 
to $\pvdocz(U)$ of 
the map $P\rightarrow k_{P}(x,x)$. 

Let $P\in \pvdoz(U)$ and let $(P(z))_{z \in 
\Omega}\subset \pvdo^{*}(U)$ be a  holomorphic gauging for $P$. Then 
it follows from~(\ref{eq:Spectral.tildeLbis}) and (\ref{eq:NCR.tau-z}) that with respect to the topology of $C^{\infty}(M,|\Lambda|(M)\otimes \End 
\cE)$ 
the map $z\rightarrow t_{P(z)}(x)$ has at worst a simple pole singularity at $z=0$ with residue
\begin{equation}
    \res_{z=0} t_{P(z)}(x)=-(2\pi)^{-(d+2)}  \int_{\|\xi\|=1}  p_{-(d+2)}(\xi) \imath_{E}d\xi =-c_{P}(x),
     \label{eq:NCR.residue.t-P.U}
\end{equation}
where $p_{-(d+2)}(\xi)$ denotes the symbol of degree $-(d+2)$ of $P$. 

Next, let $\phi:\tilde{U}\rightarrow U$ be a change of $H$-framed local coordinates. 
Let $P\in \pvdocz(U)$ and let $(P(z))_{z \in\C}$ be a  holomorphic family which is a gauging for $P$. As shown in~\cite{Po:CPDE1} 
the \psivdo\ family $(\phi^{*}P(z))_{z\in \C}$ is holomorphic and is a gauging for $\phi^{*}P$.  
Moreover, as for $\Re z$ negatively large enough we have $k_{\phi^{*}P(z)}=|\phi'(x)|k_{P(z)}(\phi(x),\phi(x))$, an analytic continuation gives 
\begin{equation}
   t_{\phi^{*}P}(x)=|\phi'(x)|t_{P}(\phi(x)).    
     \label{eq:NCR.functoriality-tP}
\end{equation}

Now, in the same way as in the construction of the density $c_{P}(x)$ in the proof of 
Proposition~\ref{thm:NCR.log-singularity}, it follows from all this that if $P\in \pvdo^{m}(M,\cE)$ 
then there exists a unique $\End \cE$-valued density $t_{P}(x)$ such that, for any local 
$H$-framed chart $\kappa:U\rightarrow V$ and any trivialization $\tau:\cE_{|_{U}}\rightarrow U\times \C^{r}$, we have 
\begin{equation}
    t_{P}(x)|_{U}=\tau^{*}\kappa^{*}(t_{\kappa_{*}\tau_{*}(P_{|_{U}})}(x)dx).
\end{equation}
On every trivializing $H$-framed chart the map $P\rightarrow t_{P}(x)$ is analytic and satisfies~(\ref{eq:NCR.residue.t-P.U}). Therefore, we obtain: 

\begin{proposition}\label{thm:NCR.TR.local}
 1) The map $P \rightarrow t_{P}(x)$ is  the unique analytic continuation of the map $P \rightarrow k_{P}(x,x)$ to $\pvdocz(M,\cE)$.\smallskip 

 2) Let $P \in \pvdoz(M,\cE)$ and let $(P(z))_{z\in\Omega}\subset \pvdo^{*}(M,\cE)$ be a holomorphic family which is a gauging 
 for $P$. Then, in $C^{\infty}(M,|\Lambda|(M)\otimes \End \cE)$, the map  $z\rightarrow t_{P(z)}(x)$ has at worst a simple pole 
 singularity at  $z=0$ with residue given by
\begin{equation}
    \res_{z=0} t_{P(z)}(x)=- c_{P}(x),
\end{equation}
where $c_{P}(x)$ denotes the $\End \cE$-valued density on $M$ given by Theorem~\ref{thm:NCR.log-singularity}.\smallskip 

3) The map $P\rightarrow t_{P}(x)$ is functorial with respect to Heisenberg diffeomorphisms as in~(\ref{eq:Log.functoriality-cP}). 
\end{proposition}

\begin{remark}
   Taking residues at $z=0$ in~(\ref{eq:NCR.functoriality-tP}) allows us to recover~(\ref{eq:Log.functoriality-cP}). 
\end{remark}

From now one we assume $M$ compact.  We then define the \emph{canonical trace} for the Heisenberg calculus as the functional $\TR$ on $\pvdocz(M,\cE)$ 
given by the formula,
\begin{equation}
     \TR P := \int_{M}\tr_{\cE} t_{P}(x) \qquad \forall P\in  \pvdocz(M,\cE).
\end{equation}

\begin{proposition}\label{thm:NCR.TR.global}
 The canonical trace $\TR$ has the following properties:\smallskip    
    
1)  $\TR$ is the unique analytic continuation 
 to $\pvdocz(M,\cE)$ of the usual trace.\smallskip 
 
2)  We have $ \TR P_{1}P_{2}=\TR P_{2}P_{1}$ whenever $\ord P_{1}+\ord P_{2}\not\in\Z$.\smallskip

3) $\TR$ is invariant by Heisenberg diffeomorphisms, i.e., for any Heisenberg diffeomorphism $\phi:(M,H)\rightarrow (M',H')$ 
we have $ \TR \phi_{*}P=\TR P$  $\forall P\in \pvdocz(M,\cE)$.
\end{proposition}
\begin{proof}
 The first and third properties are immediate consequences of Theorem~\ref{thm:NCR.TR.local}, so we only have to prove the second one. 
 
 For $j=1,2$ let $P_{j}\in \pvdo^{*}(M,\cE)$ and let $(P_{j}(z))_{z\in \C}\subset \pvdo^{*}(M,\cE)$ be a holomorphic gauging 
 for $P_{j}$. We further assume that 
 $\ord P_{1}+\ord P_{2}\not\in\Z$. Then $P_{1}(z)P_{2}(z)$ and $P_{2}(z)P_{1}(z)$ have non-integer order for $z$ in $\C\setminus \Sigma$, where
 $\Sigma:=-(\ord P_{1}+\ord P_{2})+\Z$. 
 For $\Re z$ negatively large enough we have $\Tra 
 P_{1}(z)P_{2}(z)=\Tra P_{2}(z)P_{1}(z)$, so by analytic continuation we see that 
 $\TR P_{1}(z)P_{2}(z)=\TR P_{2}(z)P_{1}(z)$ for any $z\in\C\setminus \Sigma$.  
Setting $z =0$ then shows that we have $ \TR P_{1}P_{2}=\TR P_{2}P_{1}$ as desired. 
 \end{proof}

Next, we define the \emph{noncommutative residue} for the Heisenberg calculus as the linear functional $\Res$ on $\pvdoz(M,\cE)$ given by the 
formula, 
\begin{equation}
     \Res P :=\int_{M}\tr_{\cE} c_{P}(x) \qquad \forall P\in \pvdoz(M,\cE).
\end{equation}

This functional provides us with the analogue for the Heisenberg calculus of the noncommutative residue trace of 
Wodzicki (\cite{Wo:LISA}, \cite{Wo:NCRF}) and 
Guillemin~\cite{Gu:NPWF}, for we have:  

\begin{proposition}\label{thm:NCR.NCR} 
The noncommutative residue $\Res$ has the following properties:\smallskip
    
1) Let $P \in \pvdoz(M,\cE)$ and let $(P(z))_{z\in \Omega}\subset \pvdo^{*}(M,\cE)$ be a  holomorphic gauging for $P$. 
Then at $z=0$ the function $\TR P(z)$ has at worst a simple pole singularity in such way that we have  
\begin{equation}
    \res_{z=0} \TR P(z)= - \Res P. 
     \label{eq:NCR.residueTR}
\end{equation}
  
 2) We have $\Res P_{1}P_{2}=\Res P_{2}P_{1}$ whenever $\ord P_{1}+\ord P_{2}\in \Z$. Hence $\Res$ is a trace on the algebra $\pvdoz(M,\cE)$.\smallskip
 
 3) $\Res$ is invariant by Heisenberg diffeomorphisms.\smallskip
 
 4) We have $\Res P^{t}=\Res P$ and $\Res P^{*}=\overline{\Res P}$ for any $P\in \pvdoz(M,\cE)$.
\end{proposition}
\begin{proof}
 The first property follows from Proposition~\ref{thm:NCR.TR.local}. The third and fourth properties are immediate consequences of
 Propositions~\ref{thm:NCR.log-singularity} and~\ref{prop.Sing.transpose-adjoint}. 
 
 It remains to prove the 2nd property. Let $P_{1}$ and $P_{2}$ be operators in $\pvdo^{*}(M,\cE)$ such that $\ord P_{1}+\ord P_{2}\in \Z$. For $j=1,2$ 
 let $(P_{j}(z))_{z\in \C}\subset \pvdo^{*}(M,\cE)$ be a holomorphic  gauging for $P_{j}$. Then the family 
 $(P_{1}(\frac{z}{2})P_{2}(\frac{z}{2}))_{z\in \C}$ (resp.~$(P_{2}(\frac{z}{2})P_{1}(\frac{z}{2}))_{z\in \C}$) is a holomorphic gauging for $P_{1}P_{2}$ 
 (resp.~$P_{2}P_{1}$). Moreover, by Proposition~\ref{thm:NCR.TR.global} for any $z \in\CZ$ we have 
 $\TR P_{1}(\frac{z}{2})P_{2}(\frac{z}{2})=\TR P_{\frac{z}{2}}(z)P_{1}(\frac{z}{2}$. Therefore, by 
 taking residues at $z=0$ and using~(\ref{eq:NCR.residueTR})  we get $\Res  P_{1}P_{2}=\Res P_{2}P_{1}$ as desired. 
\end{proof}

\subsection{Traces and sum of commutators}
\label{sec:traces}
Let $(M^{d+1},H)$ be a compact Heisenberg manifold and let $\cE$ be a vector bundle over $M$.  In this subsection, we shall prove that when 
$M$ is connected the noncommutative residue spans the space of traces 
on the algebra $\pvdoz(M,\cE)$. As a consequence this will allow us to characterize the sums of commutators in $\pvdoz(M,\cE)$.

Let $H\subset T\Rd$ be a hyperplane bundle such that there exists a global $H$-frame $X_{0},X_{1},\ldots,X_{d}$ of $T\Rd$. We will now give a series 
of criteria for an operator $P\in \pvdoz(\Rd)$ to be a sum of commutators of the form,
\begin{equation}
        P=  [ x_{0}, P_{0}]+\ldots+ [x_{d}, P_{d}], \qquad P_{j}\in \pvdoz(\Rd). 
         \label{eq:Traces.sum-commutators}
\end{equation}

In the sequel for any $x\in \Rd$ we let $\psi_{x}$ denote the affine change of variables to the privileged coordinates at $x$ with respect to the 
$H$-frame $X_{0},\ldots,X_{d}$. 

\begin{lemma}\label{lem:Traces.criterion-logarithmic-kernel}
    Let $P\in \pvdo^{-(d+2)}(\Rd)$ have a kernel of the form, 
    \begin{equation}
        k_{P}(x,y)=|\psi_{x}'|K_{0}(x,-\psi_{x}(y)),
         \label{eq:Traces.logarithmic-kernel}
    \end{equation}
    where $K_{0}(x,y)\in \cK_{0}(\Rd\times \Rd)$ is homogeneous of degree $0$ with respect to $y$. Then $P$ is a sum 
    of commutators of the form~(\ref{eq:Traces.sum-commutators}). 
\end{lemma}
\begin{proof}
  Set $\psi_{x}(y)=A(x).(y-x)$  with $A\in C^{\infty}(\Rd, GL_{d+1}(\Rd))$ and for 
  $j,k=0,\ldots, d$ define
\begin{equation}
    K_{jk}(x,y):= A_{jk}(x) y_{j}^{\beta_{j}}\|y\|^{-4}K_{0}(x,y), \qquad (x,y)\in \Rd\times \Rdo,
\end{equation}
 where $\beta_{0}=1$ and $\beta_{1}=\ldots=\beta_{d}=3$. As $K_{jk}(x,y)$ is smooth for $y\neq 0$ and is homogeneous with respect to $y$ of 
 degree~$-2$ if $j=0$ and of degree $-1$ otherwise, we see that  it belongs to $\cK_{*}(\R \times \R)$. Therefore, the operator $Q_{jk}$ with 
 Schwartz kernel $k_{Q_{jk}}=|\psi_{x}'| K_{jk}(x,-\psi_{x}(y))$ is a \psivdo. 
 
 Next, set  $A^{-1}(x)=(A^{jk}(x))_{1\leq j,k\leq d}$. Since $x_{k}-y_{k}=-\sum_{l=0}^{d}A^{kl}(x)\psi_{x}(y)_{l}$ we deduce that the  
 Schwartz kernel of  $\sum_{j,k=0}^{d}[x_{k},Q_{jk}]$ is $|\psi_{x}'|K(x,-\psi_{x}(y))$, where 
  \begin{multline}
    K(x,y)=   \sum_{0\leq j,k,l \leq d} A^{kl}(x)y_{l}A_{jk}(x)y_{j}^{\beta_{j}}\|y\|^{-4} K_{0}(x,y)\\ 
       =\sum_{0\leq j\leq d} y_{j}^{\beta_{j}+1}\|y\|^{-4}K_{0}(x,y)=K_{0}(x,y).
  \end{multline}
Hence $P=\sum_{j,k=0}^{d}[x_{k},Q_{jk}]$. The lemma is thus proved. 
\end{proof}

\begin{lemma}\label{lem:Traces.smoothing-operators1}
Any $R\in \psinf(\Rd)$ can be written as a sum of commutators of the form~(\ref{eq:Traces.sum-commutators}).
\end{lemma}
\begin{proof}
 Let $k_{R}(x,y)$ denote the Schwartz kernel of $R$. Since $k_{R}(x,y)$ is smooth we can write 
 \begin{equation}
     k_{R}(x,y)=k_{R}(x,x)+(x_{0}-y_{0})k_{R_{0}}(x,y)+\ldots +(x_{d}-y_{d})k_{R_{d}}(x,y),
      \label{eq:Traces.Taylor}
 \end{equation}
for some smooth functions $k_{R_{0}}(x,y),\ldots,k_{R_{d}}(x,y)$. For $j=0,\ldots,d$ let $R_{j}$ be the smoothing operator with Schwartz 
kernel $k_{R_{j}}(x,y)$, and 
let $Q$ be the operator with Schwartz kernel $k_{Q}(x,y)=k_{R}(x,x)$. Then by~(\ref{eq:Traces.Taylor}) 
we have  
\begin{equation}
    R=Q+[x_{0},R_{0}]+\ldots + [x_{d},R_{d}].
    \label{eq:Traces.commutators.smoothing-RQ}
\end{equation} 

Observe that the kernel of $Q$ is of the form~(\ref{eq:Traces.logarithmic-kernel}) with $K_{0}(x,y)=|\psi_{x}'|^{-1}k_{R}(x,x)$. 
Here $K_{0}(x,y)$ belongs to  $\cK_{0}(\Rd\times \Rd)$ 
and is homogeneous of degree $0$ with respect to $y$, so by Lemma~\ref{lem:Traces.criterion-logarithmic-kernel} 
the operator $Q$ is a sum of commutators of the form~(\ref{eq:Traces.sum-commutators}).  
Combining this with~(\ref{eq:Traces.commutators.smoothing-RQ}) then shows that $R$ is of that form too.
\end{proof}

\begin{lemma}\label{lem:Traces.sum-commutators}
    Any $P\in \pvdo^{\Z}(\Rd)$ such that $c_{P}(x)=0$ is a sum of commutators of the form~(\ref{eq:Traces.sum-commutators}).
\end{lemma}
\begin{proof}
For $j=0,\ldots,d$ we let $\sigma_{j}(x,\xi)=\sum_{k=0}^{d}\sigma_{jk}(x)\xi_{k}$ denote the classical symbol of 
$-iX_{j}$. Notice that $\sigma(x):=(\sigma_{jk}(x))$ belongs to $C^{\infty}(\Rd, GL_{d+1}(\C))$.\smallskip 
    
    (i) Let us first assume that $P= (\partial_{\xi_{j}}q)(x,-iX)$ for some $q\in S^{\Z}(\URd)$. Set $q_{\sigma}(x,\xi)=q(x,\sigma(x,\xi))$. Then we 
    have
\begin{equation*}
[q(x,-iX),x_{k}]=   [q_{\sigma}(x,D),x_{k}]= (\partial_{\xi_{k}}q_{\sigma})(x,D)=\sum_{l} \sigma_{lk}(x)(\partial_{\xi_{l}}q)(x,-iX). 
\end{equation*}
 Therefore, if we let $(\sigma^{kl}(x))$ be the inverse matrix of $\sigma(x)$, then we see that 
\begin{equation*}
     \sum_{k}[\sigma^{jk}(x)q(x,-iX),x_{k}]= \sum_{k,l} \sigma^{jk}(x)\sigma_{lk}(x)(\partial_{\xi_{l}}q)(x,-iX)=(\partial_{\xi_{j}}q)(x,-iX)=P.
\end{equation*}
Hence $P$ is a sum of commutators of the form~(\ref{eq:Traces.sum-commutators}).\smallskip  

(ii) Suppose now that $P$ has symbol $p \sim \sum_{j\leq m}p_{j}$ with $p_{-(d+2)}=0$. 
Since $p_{l}(x,\xi)$ is homogeneous of degree $l$ with respect to $\xi$,  the Euler identity, 
\begin{equation}
       2\xi_{0}\partial_{\xi_{0}}p_{l} + \xi_{1}\partial_{\xi_{1}}p_{l}+\ldots+ \xi_{d}\partial_{\xi_{d}}p_{l} 
            = l p_{l},
\end{equation}
        implies that we have 
        \begin{equation}
            2\partial_{\xi_{0}}(\xi_{0}p_{l}) +   \partial_{\xi_{1}}(\xi_{1}p_{l})+\ldots+ 
                  \partial_{\xi_{d}}(\xi_{d}p_{l})= (l+d+2)p_{l}.
        \end{equation}

      For $j=0, \ldots, d$ let $q^{(j)}$ be a symbol so that $q^{(j)} 
  \sim \sum_{l\neq -(d+2)} (l+d+2)^{-1} \xi_{j}p_{l}$. 
     Then for $l\neq -(d+2)$ the  symbol of degree $l$ of $2\partial_{\xi_{0}}q^{(0)}+  
        \partial_{\xi_{1}} q^{(1)}+\ldots+ \partial_{\xi_{j}} q^{(d)}$ is equal to
\begin{equation}
    (l+d+2)^{-1}(2\partial_{\xi_{0}}(\xi_{0}p_{l}) +   \partial_{\xi_{1}}(\xi_{1}p_{l}))+\ldots+  
                    \partial_{\xi_{d}}(\xi_{d}p_{l}))= p_{l}. 
\end{equation}
 Since $p_{-(d+2)}=0$ this shows that $p- 2\partial_{\xi_{0}}q^{(0)} +  
        \partial_{\xi_{1}} q^{(1)}+\ldots+ \partial_{\xi_{j}} q^{(d)}$ is in $S^{-\infty}(\Rd\times \Rd)$.  Thus, there exists  $R$ in $\psinf(\Rd)$ 
        such that
       \begin{equation}
           P=2(\partial_{\xi_{0}}q^{(0)})(x,-iX) +  
        (\partial_{\xi_{1}} q^{(1)})(x,-iX)+\ldots+ (\partial_{\xi_{j}} q^{(d)})(x,-iX)+R,
       \end{equation}

       Thanks to the part (i) and to Lemma~\ref{lem:Traces.smoothing-operators1}
       the operators $(\partial_{\xi_{j}}q^{(j)})(x,-iX)$ and $R$ 
       are sums of commutators of the form~(\ref{eq:Traces.sum-commutators}), so $P$ is of that form as well.\smallskip
        
        (iii) The general case is obtained as follows. Let $p_{-(d+2)}(x,\xi)$ be the symbol of degree $-(d+2)$ of $P$. 
        Then by Lemma~\ref{lem:NCR.extension-symbolU} we can 
        extend $p_{-(d+2)}(x,\xi)$ into a distribution $\tau(x,\xi)$ in $C^{\infty}(\Rd)\hotimes \cS'(\Rd)$ in such way that 
        $K_{0}(x,y):=\check{\tau}_{\xiy}(x,y)$ belongs to $\cK_{0}(\Rd\times\Rd)$. Furthermore, with the notation of~(\ref{eq:NCR.log-homogeneity-Km}) we have 
        $c_{K,0}(x)=(2\pi)^{-(d+2)}\int_{\|\xi\|=1}p_{-(d+2)}(x,\xi)\iota_{E}d\xi$. Therefore, by using~(\ref{eq:NCR.formula-cP}) and the fact $c_{P}(x)$ is zero, we 
        see that $c_{K,0}(x)=|\psi'_{x}|^{-1}c_{P}(x)=0$. In view of~(\ref{eq:NCR.log-homogeneity-Km}) 
        this show that $K_{0}(x,y)$ is homogeneous of degree $0$ with 
        respect to $y$. 
        
       Let $Q\in \pvdo^{-(d+2)}(\Rd)$ be the \psivdo\ with Schwartz kernel $|\psi_{x}'|K_{0}(x,-\psi_{x}(y))$. 
       Then by Lemma~\ref{lem:Traces.criterion-logarithmic-kernel} the operator $Q$ 
       is a sum of commutators of the form~(\ref{eq:Traces.sum-commutators}). Moreover, observe that by Proposition~\ref{prop:PsiVDO.characterisation-kernel1}
       the operator $Q$ has symbol $q\sim 
       q_{-(d+2)}$, where for $\xi \neq 0$ we have $q_{-(d+2)}(x,\xi)=(K_{0})^{\wedge}_{\yxi}(x,\xi)=p_{-(d+2)}(x,\xi)$. Therefore $P-Q$ is a 
       \psivdo\ whose symbol of degree $-(d+2)$ is zero. It then follows from the part (ii) of the proof that $P-Q$ is a sum of 
       commutators of the form~(\ref{eq:Traces.sum-commutators}). 
       All this shows that $P$ is the sum of two operators of the form~(\ref{eq:Traces.sum-commutators}), so $P$ is of that form too. 
\end{proof} 

In the sequel we let $\pvdoc^{*}(\Rd)$ and $\psinfc(\Rd)$ respectively denote the classes of \psivdos\ and smoothing operators on $\Rd$ with 
compactly supported Schwartz kernels. 

\begin{lemma}\label{lem:Traces.sum-commutators.compact}
    There exists $\Gamma \in \pvdoc^{-(d+2)}(\Rd)$ such that, for any $P\in \pvdoc^{\Z}(\Rd)$, we have
    \begin{equation}
        P= (\Res P)\Gamma \quad \bmod [\pvdoc^{\Z}(\Rd), \pvdoc^{\Z}(\Rd)].
         \label{eq:Traces.sum-commutators-compact}
    \end{equation}
\end{lemma}
\begin{proof} Let $P\in \pvdoc^{\Z}(\Rd)$. We will put $P$ into the form~(\ref{eq:Traces.sum-commutators-compact}) in 3 steps.\smallskip 
    
\noindent (i) Assume first that $c_{P}(x)=0$. Then by Lemma~\ref{lem:Traces.sum-commutators} we can write $P$ in the form, 
\begin{equation}
    P=  [ x_{0}, P_{0}]+\ldots+ [x_{d}, P_{d}], \qquad P_{j}\in \pvdoz(\Rd). 
\end{equation}
Let $\chi$ and $\psi$ in $C^{\infty}_{c}(\Rd)$ be such that $\psi(x)\psi(y)=1$ near the support of the kernel of $P$ and $\chi=1$ near $ \supp \psi$. 
Since $\psi P \psi=P$ we obtain
\begin{equation}
    P=\sum_{j=0}^{d} \psi [x_{d}, P_{d}]\psi= \sum_{j=0}^{d} [x_{d}, \psi P_{d}\psi] =\sum_{j=0}^{d}[\chi x_{d}, \psi P_{d}\psi].  
     \label{eq:Traces.commutators.vanishing-cP}
\end{equation}
In particular $P$ is a sum of commutators in $\pvdoc^{\Z}(\Rd)$.

(ii) Let $\Gamma_{0}\in \pvdo^{-(d+2)}$ have kernel $k_{\Gamma_{0}}(x,y)=-\log\|\phi_{x}(y)\|$ and suppose that $P=c \Gamma_{0} \psi$ where $c\in 
C^{\infty}_{c}(\Rd)$ is such that $\int c(x)dx=0$ and $\psi \in C^{\infty}_{c}(\Rd)$ is such that $\psi=1$ near $\supp c$. First, we have: 

\begin{claim}
   If $c\in C^{\infty}_{c}(\Rd)$ is such that $\int c(x)dx=0$, then there exist $c_{0}, \ldots,c_{d}$ in $C_{c}^{\infty}(\Rd)$ such that 
   $c=\partial_{x_{0}}c_{0}+\ldots+\partial_{x_{d}}c_{d}$.
\end{claim}
\begin{proof}[Proof of the Claim] We proceed by induction on the dimension $d+1$. In dimension $1$ the proof follows from the the fact that if $c\in 
    C^{\infty}_{c}(\R)$ is such that $\int_{-\infty}^{\infty}c(x_{0})dx_{0}=0$, 
    then $\tilde{c}(x_{0})=\int_{-\infty}^{x_{0}}c(t)dt$ is an antiderivative of  $c$ with compact support. 
  
 Assume now that the claim is true in dimension $d$ and under this assumption let us prove it in dimension $d+1$. 
 Let $c\in C^{\infty}_{c}(\R^{d+1})$ be such that 
 $\int_{\R^{d+1}} c(x)dx=0$. For any $(x_{0},\ldots,x_{d-1})$ in $\R^{d}$ we let 
 $\tilde{c}(x_{0},\ldots,x_{d-1})=\int_{\R}c(x_{0},\ldots,x_{d-1},x_{d})dx_{d}$. This defines a function in
 $C^{\infty}_{c}(\R^{d})$ such that 
\begin{equation}
     \int_{\R^{d}}\tilde{c}(x_{0},\ldots,x_{d-1})dx_{0}\ldots dx_{d-1}= \int_{\Rd}c(x_{0},\ldots,x_{d})dx_{0}\ldots dx_{d}=0.
\end{equation}
Since the claim is assumed to hold in dimension $d$, it follows that there exist $\tilde{c}_{0}, 
 \ldots,\tilde{c}_{d-1}$ in $C_{c}^{\infty}(\R^{d})$ such that $ \tilde{c}=\sum_{j=0}^{d-1}\partial_{x_{j}}\tilde{c}_{j}$. 
 
Next, let $\varphi\in C^{\infty}_{c}(\R)$ be such that $\varphi(x_{d})dx_{d}=1$.  For any $(x_{0},\ldots,x_{d})$ in $\Rd$ we let 
\begin{equation}
     b(x_{0},\ldots,x_{d})=c(x_{0},\ldots,x_{d})-\varphi(x_{d})\tilde{c}(x_{0},\ldots,x_{d-1}). 
\end{equation}
This defines a function in $C^{\infty}_{c}(\R^{d+1})$ such that
\begin{equation}
    \int_{-\infty}^{\infty}b(x_{0},\ldots,x_{d})dx_{d}=\int_{-\infty}^{\infty}c(x_{0},\ldots,x_{d})dx_{d}-\tilde{c}(x_{0},\ldots,x_{d-1})=0.
\end{equation}
Therefore, we have $b=\partial_{x_{d}}c_{d}$, where $c_{d}(x_{0},\ldots,x_{d}):=\int_{-\infty}^{x_{d}}b(x_{0},\ldots,x_{d-1},t)dt$ is a function 
in $C^{\infty}_{c}(\R^{d+1})$.

In addition, for $j=0,\ldots,d-1$ and for $(x_{0},\ldots,x_{d})$ in $\Rd$ we let $c_{j}(x_{0},\ldots,x_{d})=\varphi(x_{d})\tilde{c}(x_{0},\ldots,x_{d-1})$.
Then $c_{0},,\ldots,c_{d-1}$ belong to $C^{\infty}_{c}(\Rd)$ and we have 
 \begin{multline}
  c(x_{0},\ldots,x_{d})=    b(x_{0},\ldots, x_{d})+\varphi(x_{d})\tilde{c}(x_{0},\ldots,x_{d-1})  \\ = \partial_{x_{d}} c_{d}(x_{0},\ldots,x_{d}) 
     +\varphi(x_{d})\sum_{j=0}^{d-1}\partial_{x_{j}}\tilde{c}_{j}(x_{0},\ldots,x_{d-1}) 
     = \sum_{j=0}^{d}\partial_{x_{j}}{c}_{j}.
 \end{multline}
This shows that the claim is true in dimension $d+1$. The proof is now complete. 
\end{proof}

Let us now go back to the proof of the lemma. Since we have $\int c(x)dx=0$ the above claim tells us that $c$ can be written in the form
$c=\sum_{j=0}^{d}\partial_{j}c_{j}$ with $c_{0},\ldots,c_{d}$ in $C^{\infty}_{c}(\Rd)$. 
Observe also that the Schwartz kernel 
   of $[\partial_{x_{j}},\Gamma_{0}]$ is equal to 
  \begin{multline}
          ( \partial_{x_{j}}- \partial_{y_{j}})[- \log \|\psi_{x}(y)\|] \\ = \sum_{k,l}  ( \partial_{x_{j}}- \partial_{y_{j}}) 
          [\varepsilon_{kl}(x)(x_{l}-y_{l})][\partial_{z_{k}}\log \|z\|]_{z=-\psi_{x}(y)}\\
          = \sum_{k,l} (x_{k}-y_{k}) (\partial_{x_{j}}\varepsilon_{kl})(x) 
          \gamma_{k}(-\psi_{x}(y))\|\psi_{x}(y)\|^{-4},
  \end{multline}
  where we have let $\gamma_{0}(y)=\frac{1}{2} y_{0}$ and $\gamma_{k}(y)=y_{k}^{3}$, $k=1,\ldots,d$. In particular 
  $k_{[\partial_{x_{j}},\Gamma_{0}]}(x,y)$  has no logarithmic singularity near the diagonal, that is, we have
  $c_{[\partial_{x_{j}},\Gamma_{0}]}(x)=0$.

Next, let $\psi \in C^{\infty}_{c}(\Rd)$ be such that $\psi=1$ near $\supp c\cup\supp c_{1}\cup \cdots \cup \supp c_{d}$ and let $\chi \in 
C^{\infty}_{c}(\Rd)$ be such that $\chi=1$ near $\supp \psi$. Then we have
\begin{multline}
  [\chi \partial_{x_{j}}, c_{j}\Gamma_{0}\psi]= [\partial_{x_{j}}, c_{j}\Gamma_{0}\psi] =
    [\partial_{x_{j}},c_{j}]\Gamma_{0} \psi + c_{j} 
          [\partial_{x_{j}}, \Gamma_{0}]\psi+ c_{j}\Gamma_{0}  [\partial_{x_{j}},\psi ]\\
         = \partial_{x_{j}}c_{j} \Gamma_{0} \psi + c_{j} [\partial_{x_{j}}, \Gamma_{0}]\psi + c_{j}\Gamma_{0}  \partial_{x_{j}}\psi.
\end{multline}
Since $ c_{j}\Gamma_{0}  \partial_{x_{j}}\psi$ is smoothing and $c_{c_{j} [\partial_{x_{j}}, \Gamma_{0}]\psi}(x)=c_{j}c_{[\partial_{x_{j}}, 
\Gamma_{0}]}(x)=0$ we deduce from this that $P$ is of the form
$P= \sum_{j=0}^{d}[\chi \partial_{x_{j}}, c_{j}\Gamma_{0}\psi] +Q$
with $Q\in \pvdoc^{\Z}(\Rd)$ such that $c_{Q}(x)=0$. It then follows from the part (i) that $P$ belongs to the commutator space of 
$\pvdoc^{\Z}(\Rd)$. 

(iii) Let $\rho \in C_{c}^{\infty}(\Rd)$ be such that $\int \rho(x)dx=1$, let $\psi \in C^{\infty}_{c}(\Rd)$ be such that $\psi=1$ near $\supp \rho$, 
and set $\Gamma=\rho \Gamma_{0}\psi$. Let $P\in \pvdoc^{\Z}(\Rd)$ and let $\tilde{\psi}\in C^{\infty}_{c}(\Rd)$ be such that 
$\tilde{\psi}=1$ near $\supp c_{P}\cup \supp \psi$. Then we have 
\begin{equation}
    P=(\Res P) \Gamma + (\Res P)\rho \Gamma_{0}(\tilde{\psi}-\psi)+ (c_{P}-(\Res P)\rho)\Gamma_{0}\tilde{\psi}+ P-c_{P}\Gamma_{0}\tilde{\psi}. 
    \label{eq:Traces.decomposition-P-compact-support}
\end{equation}

Notice that $(\Res P)\rho \Gamma_{0}(\tilde{\psi}-\psi)$ belongs to $\psinfc(\Rd)$. Observe also that the logarithmic singularity of 
$P-c_{P}\Gamma_{0}\tilde{\psi}$ is equal to $c_{P}(x)-\tilde{\psi}(x)c_{P}(x)=0$. Therefore, it follows from (i) that these operators 
belong to commutator space of $\pvdoc^{\Z}(\Rd)$. In addition, as $\int (c_{P}(x)-(\Res P)\rho(x))dx=0$ we see that $(c_{P}-(\Res 
P)\rho)\Gamma_{0}\tilde{\psi}$ is as in (ii), so it also belongs to the commutator space of $\pvdoc^{\Z}(\Rd)$. Combining all this 
with~(\ref{eq:Traces.decomposition-P-compact-support}) then shows that $P$ agrees with $(\Res P) \Gamma $ modulo 
a sum of commutators in $\pvdoc^{\Z}(\Rd)$. The lemma is thus proved. 
\end{proof}

Next, we quote the well known lemma below. 
\begin{lemma}[{\cite[Appendix]{Gu:RTCAFIO}}]\label{lem:Traces.smoothing-operators2}
    Any $R \in \psinf(M,\cE)$ such that $\Tr R=0$ is the sum of two commutators in $\psinf(M,\cE)$. 
\end{lemma}

We are now ready to prove the main result of this section. 
\begin{theorem}\label{thm:Traces.traces}
   Assume that $M$ is connected. Then any trace on $\pvdoz(M,\cE)$  is a constant multiple of the noncommutative residue. 
\end{theorem}
\begin{proof}
Let $\tau$ be a trace on $\pvdoz(M,\cE)$.  By Lemma~\ref{lem:Traces.sum-commutators.compact} there exists $\Gamma \in \pvdoc^{-(d+2)}(\Rd)$ such that any 
$P=(P_{ij})$ in $\pvdoc^{\Z}(\Rd,\C^{r})$ can be written as  
\begin{equation}
      P=\Gamma \otimes R \bmod [\pvdoc^{\Z}(\Rd), \pvdoc^{\Z}(\Rd)]\otimes M_{r}(\C),
\end{equation}
 where we have let $R=(\Res P_{ij})\in M_{r}(\C)$.  Notice that $\Tr R= \sum 
\Res P_{ii}=\Res P$. Thus $R-\frac{1}{r}(\Res P)I_{r}$ has a vanishing trace, hence belongs to the commutator space 
of $M_{r}(\C)$. Therefore, we have
\begin{equation}
    P=(\Res P) \Gamma\otimes (\frac{1}{r}I_{r}) \quad \bmod [\pvdoc^{\Z}(\Rd, \C^{r}), \pvdoc^{\Z}(\Rd,\C^{r})].
     \label{eq:Traces.decomposition-P-Rd-Cr}
\end{equation}

Let $\kappa:U\rightarrow \Rd$ be a local $H$-framed chart mapping onto $\Rd$ and such that $\cE$ is trivializable over its domain. For sake of 
terminology's brevity we shall call such a chart a \emph{nice $H$-framed chart}. As $U$ is $H$-framed and is Heisenberg diffeomorphic to $\Rd$ and as$\cE$ 
is trivializable over $U$, it follows from~(\ref{eq:Traces.decomposition-P-Rd-Cr}) that there exists $\Gamma_{U}\in 
\pvdoc^{-(d+2)}(U,\cE_{|_{U}})$ such that, for any $P \in \pvdoc^{\Z}((U,\cE_{|_{U}})$, we have
\begin{equation}
    P=(\Res P)\Gamma_{U} \quad \bmod [\pvdoc^{\Z}(U,\cE_{|_{U}}), \pvdoc^{\Z}(U,\cE_{|_{U}})].
     \label{eq:Traces.local-form-Psivdos}
\end{equation}
If we apply the trace $\tau$, then we see that, for any $P \in \pvdoc^{\Z}(U,\cE_{|_{U}})$, we have 
\begin{equation}
    \tau(P) = \Lambda_{U}\Res P, \qquad \Lambda_{U}:=\tau(\Gamma_{U}).
\end{equation}

Next, let $\cU$ be the set of points $x \in M$ near which there a domain $V$ of a nice $H$-framed chart such that $\Lambda_{V}=\Lambda_{U}$. Clearly 
$\cU$ is a non-empty open subset of $M$. Let us prove that $\cU$ is closed. Let $x \in \overline{\cU}$ and let $V$ be an open neighborhood of $x$ which 
is the domain a nice $H$-framed chart (such a neighborhood always exists). Since $x$ belongs to the closure of $\cU$ the set $\cU\cup V$ is 
non-empty. Let $y \in \cU\cup V$. As $y$ belongs to $\cU$ there exists an open neighborhood $W$ of $y$ which is the domain a nice $H$-frame chart such 
that $\Lambda_{W}=\Lambda_{U}$. Then for any $P$ in $\pvdoc^{\Z}(V\cap W, \cE_{|V\cap W})$ we have $\tau(P)=\Lambda_{V}\Res P =\Lambda_{W}\Res P$. 
Choosing $P$ so that $\Res P\neq 0$ then shows that $\Lambda_{V}=\Lambda_{W}=\Lambda_{U}$. Since $V$ contains $x$ and 
is a domain of a nice $H$-framed chart we deduce that $x$ belongs to $\cU$.  Hence $\cU$ is both closed and open. As $M$ is connected it follows that 
$\cU$ agrees with $M$. Therefore, if we set $\Lambda=\Lambda_{U}$ then, for any domain $V$ of a nice $H$-framed chart, we have  
\begin{equation}
    \tau(P)=\Lambda \Res P \qquad \forall P \in \pvdoc^{\Z}(V,\cE_{|_{V}}).
     \label{eq:Traces.mutiple-local}
\end{equation}

Now, let $(\varphi_{i})$ be a finite partition of the unity subordinated to an open covering $(U_{i})$ of $M$ by domains of nice $H$-framed charts. 
For each index $i$ let $\psi_{i}\in C^{\infty}_{c}(U_{i})$ be such that $\psi_{i}=1$ near $\supp \varphi_{i}$. Then any $P \in 
\pvdoz(M,\cE)$ can be written as  $P=\sum \varphi_{i} P\psi_{i}+R$, where $R$ is a smoothing operator whose kernel vanishes near the diagonal of $M\times M$. In 
particular we have $\Tra R=0$, so by Lemma~\ref{lem:Traces.smoothing-operators2} 
the commutator space of $\pvdoz(M,\cE)$ contains $R$. Since each operator $\varphi_{i} P\psi_{i}$ 
can be seen as an element of $\pvdoc^{\Z}(U_{i},\cE_{|_{U_{i}}})$, using~(\ref{eq:Traces.mutiple-local}) we get 
\begin{equation}
    \tau(P)=\sum \tau(\varphi_{i}P\psi_{i}) = \sum \Lambda \Res \varphi_{i}P\psi_{i}=\Lambda \Res P.
\end{equation}
Hence we have $\tau =\Lambda \Res$. This shows that any trace on $\pvdoz(M,\cE)$ is proportional to the noncommutative residue. 
\end{proof}

Since the dual of $\pvdoz(M,\cE)/[\pvdoz(M,\cE),\pvdoz(M,\cE)]$ is isomorphic to the space of traces on $\pvdoz(M,\cE)$, as a consequence 
of Theorem~\ref{thm:Traces.traces} we get: 
\begin{corollary}
    Assume $M$ connected. Then an operator $P \in \pvdoz(M,\cE)$ is a sum of commutators in $\pvdoz(M,\cE)$ if and only if its noncommutative 
    residue vanishes. 
\end{corollary}

\begin{remark}
In~\cite{EM:HAITH} Epstein and Melrose computed the Hochschild homology of the algebra of symbols $\pvdoz(M,\cE)/\psinf(M,\cE)$ when $(M,H)$ is a 
contact manifold. In fact, as the 
algebra $\psinf(M,\cE)$ is $H$-unital and its Hochschild homology is known, the long exact sequence of~\cite{Wo:LESCHAEA} holds and allows us to relate the 
Hochschild homology of $\pvdoz(M,\cE)$ to that of $\pvdoz(M,\cE)/\psinf(M,\cE)$. In particular, we can recover from this that the space of traces on 
$\pvdoz(M,\cE)$ is one-dimensional when the manifold is connected.
 \end{remark}

\section{Analytic Applications on general Heisenberg manifolds}
\label{sec:Analytic-Applications}
 In this section we derive several  analytic applications of the construction of the noncommutative residue trace for the 
 Heisenberg calculus.  
 First, we deal with zeta functions of hypoelliptic \psivdos\ and relate their singularities to the heat kernel asymptotics of the corresponding 
 operators. 
Second, we give logarithmic metric estimates for Green kernels of hypoelliptic \psivdos\  whose order is equal to the Hausdorff dimension $\dim M +1$. 
 This connects nicely with previous results of Fefferman, Stein and their students and collaborators. 
 Finally, we show that the noncommutative 
 residue for the Heisenberg calculus allows us to extend the Dixmier trace to the whole algebra of integer order \psivdos. This is the analogue for 
 the Heisenberg calculus of a well-known result of Alain Connes.

\subsection{Zeta functions of hypoelliptic \psivdos}
Let $(M^{d+1}, H)$ be a compact Heisenberg manifold equipped with a smooth density~$>0$, let $\cE$ be a Hermitian vector bundle over 
$M$ of rank $r$, and let $P:C^{\infty}(M,\cE)\rightarrow C^{\infty}(M,\cE)$ be a \psivdo\ of integer order $m\geq 1$ with an invertible principal 
symbol. In addition, assume that there is a ray $L_{\theta}=\{\arg \lambda =\theta\}$ which is is not through an eigenvalue of $P$ and is 
a principal cut for the principal symbol $\sigma_{m}(P)$ as in 
Section~\ref{sec:Heisenberg-calculus}. 

Let $(P_{\theta}^{s})_{s\in \C}$ be the associated family of 
complex powers associated to $\theta$ as in Proposition~\ref{prop:Heisenberg.powers2}. 
Since $(P_{\theta}^{s})_{s \in \C}$ is a holomorphic family of \psivdos, Proposition~\ref{thm:NCR.TR.global} allows us to
directly define the zeta function $\zeta_{\theta}(P;s)$ as the meromorphic function,
    \begin{equation}
        \zeta_{\theta}(P;s):=\TR P_{\theta}^{-s}, \qquad s \in \C.
    \end{equation}
\begin{proposition}\label{prop:Zeta.zeta-function}
    Let $\Sigma=\{-\frac{d+2}{m}, -\frac{d+1}{m},\ldots, \frac{-1}{m},\frac{1}{m}, 
    \frac{2}{m}, \ldots\}$. Then the function $\zeta_{\theta}(P;s)$ is analytic outside $\Sigma$, and on $\Sigma$ it has at  worst simple pole singularities such that 
 \begin{equation}
     \res_{s=\sigma}\zeta_{\theta}(P;s)=m\Res P^{-\sigma}_{\theta}, \qquad \sigma\in \Sigma. 
     \label{eq:Zeta.residue-zeta}
 \end{equation}
 In particular, $\zeta_{\theta}(P;s)$ is always regular at $s=0$.
\end{proposition}
\begin{proof}
    Since $\ord P_{\theta}^{-s}=ms$ it follows from Proposition~\ref{thm:NCR.NCR} that $\zeta_{\theta}(P;s)$ is analytic outside $\Sigma':=\Sigma\cup\{0\}$ and 
    on $\Sigma'$ has at  worst simple 
    pole singularities satisfying~(\ref{eq:Zeta.residue-zeta}). At $s=0$ we have $ \res_{s=0}\zeta_{\theta}(P;s)=m\Res P^{0}_{\theta}=m\Res 
    [1-\Pi_{0}(P)]$, but as $\Pi_{0}(P)$ is a smoothing operator we have $\Res [1-\Pi_{0}(P)]=-\Res \Pi_{0}(P)=0$. Thus 
    $\zeta_{\theta}(P;s)$  is regular at $s=0$. 
\end{proof}

Assume now that $P$ is selfadjoint and the union set of its principal cuts is $\Theta(P)=\C\setminus [0,\infty)$. This implies that $P$ is bounded from 
below (see~\cite{Po:CPDE1}), so its spectrum is real and contains at most finitely many negative eigenvalues.  
We will use the subscript $\uparrow$  (resp.~$\downarrow$) to refer to a spectral cutting in the upper halfplane $\Im \lambda>0$ (resp.~lower halfplane $\Im 
\lambda<0$). 

Since $P$ is bounded from below it defines a heat semigroup $e^{-tP}$, $t\geq 0$, and, as the principal symbol of $P$ is 
invertible,  for $t>0$ the operator $e^{-tP}$ is smoothing, hence has a smooth Schwartz kernel $k_{t}(x,y)$ in 
$C^{\infty}(M,\cE)\hotimes C^{\infty}(M,\cE^{*}\otimes |\Lambda|(M))$. Moreover, as 
$t\rightarrow 0^{+}$ we have the heat kernel asymptotics, 
\begin{equation}
    k_{t}(x,x)\sim t^{-\frac{d+2}{m}}\sum_{j\geq 0} t^{\frac{j}{m}}a_{j}(P)(x) + \log t\sum_{k\geq 0}t^{k}b_{k}(P)(x),
     \label{eq:Zeta.heat-kernel-asymptotics}
\end{equation}
where the asymptotics takes place in $C^{\infty}(M,\End \cE \otimes |\Lambda|(M))$, and when $P$ is a differential operator we have 
$a_{2j-1}(P)(x)=b_{j}(P)(x)=0$ for all $j\in \N$ (see~\cite{BGS:HECRM}, \cite{Po:MAMS1} when $P$ is a differential operator and 
see~\cite{Po:CPDE1} for the general case). 

\begin{proposition}\label{prop:Zeta.heat-zeta-local} For $j=0,1,\ldots$ set $\sigma_{j}=\frac{d+2-j}{m}$. Then:\smallskip 
    
 1)  When $\sigma_{j}\not \in \Z_{-}$ we have 
 \begin{equation}
     \res_{s=\sigma_{j}}t_{P_{\updown}^{-s}}(x)=m c_{P^{-\sigma_{j}}}(x) = \Gamma(\sigma_{j})^{-1}a_{j}(P)(x).
      \label{eq:Zeta.tPs-heat1}
 \end{equation}
 2) For $k=1,2,\ldots$ we have
 \begin{gather}
       \res_{s=-k}t_{P_{\updown}^{-s}}(x)=m c_{P^{k}}(x) = (-1)^{k+1}k!b_{k}(P)(x),
             \label{eq:Zeta.tPs-heat2}\\
       \lim_{s\rightarrow -k}[t_{P_{\updown}^{-s}}(x)-m (s+k)^{-1}c_{P^{k}}(x)] = (-1)^{k}k! a_{d+2+mk}(P)(x).
\end{gather}
3) For $k=0$ we have
\begin{equation}
    \lim_{s\rightarrow 0} t_{P_{\updown}^{-s}}(x) =a_{d+2}(P)(x)-t_{\Pi_{0}}(x).   
     \label{eq:Zeta.tPs-heat4}
\end{equation}
\end{proposition}
\begin{remark}
    When $P$ is positive and invertible the result is a standard consequence of the Mellin formula (see, e.g., \cite{Gi:ITHEASIT}). Here it is 
    slightly more complicated because we don't 
    assume that $P$ is positive or invertible. 
\end{remark}
\begin{proof}
For $\Re s>0$ set $Q_{s}= \Gamma(s)^{-1}\int_{0}^{1}t^{s-1}e^{-tP}dt$. Then we have:
 \begin{claim}
    The family $(Q_{s})_{\Re s>0}$ can be uniquely extended to a holomorphic family of \psivdos\ parametrized by $\C$ in such way that:\smallskip
    
    (i) The families $(Q_{s})_{s \in \C}$ and $(P_{\updown}^{-s})_{s \in \C}$ agree up to a holomorphic family of smoothing operators;\smallskip
    
    (ii) We have $Q_{0}=1$ and $Q_{-k}=P^{k}$ for any integer $k\geq 1$.
\end{claim}
\begin{proof}[Proof of the claim]
First, let $\Pi_{+}(P)$ and $\Pi_{-}(P)$ denote the orthogonal projections onto the positive and negative eigenspaces of $P$. Notice 
that $\Pi_{-}(P)$ is a smoothing 
operator because $P$ has at most only finitely many negative eigenvalues. For $\Re s>0$ the Mellin formula allows us to write
\begin{equation}       
    P_{\updown}^{-s}=\Pi_{-}(P) P_{\updown}^{-s}+\Gamma(s)^{-1}\int_{0}^{\infty}t^{s}\Pi_{+}(P)e^{-tP}\frac{dt}{t}=Q_{s}+R_{\updown}(s),
    \label{eq:Zeta.claim-heat.PQs}\\ 
\end{equation}
where  $R_{\updown}(s)$ is equal to
\begin{equation}       
         \Pi_{-}(P) P_{\updown}^{-s}-s^{-1}\Gamma(s)^{-1}\Pi_{0}(P)-\Pi_{-}(P)\int_{0}^{1}t^{s}e^{-tP}\frac{dt}{t}+ 
        \int_{1}^{\infty} t^{s}\Pi_{+}(P)e^{-tP}\frac{dt}{t}.
\end{equation}
Notice that $ (\Pi_{-}(P) P_{\updown}^{-s})_{s \in \C}$ and $(s^{-1}\Gamma(s)^{-1}\Pi_{0}(P))_{s \in \C}$ are holomorphic families of smoothing 
operators because $\Pi_{-}(P)$ and $\Pi_{0}(P)$ are smoothing operators. Moreover, upon writing 
\begin{gather}
    \Pi_{-}(P)\int_{0}^{1}t^{s}e^{-tP}\frac{dt}{t}=\Pi_{-}(P)(\int_{0}^{1} t^{s}e^{-tP}\frac{dt}{t})\Pi_{-}(P),\\ 
    \int_{1}^{\infty} t^{s}\Pi_{+}(P)e^{-tP}\frac{dt}{t}=  e^{-\frac{1}{4}P}(\int_{1/2}^{\infty} t^{s}\Pi_{+}(P)e^{-tP}\frac{dt}{t})e^{-\frac{1}{4}P},
\end{gather}
we see that $( \Pi_{-}(P)\int_{0}^{1}t^{s}e^{-tP}\frac{dt}{t})_{\Re s>0}$ and $( \int_{1}^{\infty} t^{s}\Pi_{+}(P)e^{-tP}\frac{dt}{t})_{\Re s>0}$ are 
holomorphic families of smoothing operators. Therefore $(R_{\updown}(s))_{\Re s>0}$ is a holomorphic family of smoothing operators and    
using~(\ref{eq:Zeta.claim-heat.PQs}) we see that $(Q_{s})_{\Re s>0}$ is a holomorphic family of \psivdos. 

Next, an integration by parts gives 
\begin{equation}
        \Gamma(s+1)P Q_{s+1}= \int_{0}^{1} 
        t^{s}\frac{d}{dt}(e^{-tP})=  e^{-P} + s  \int_{0}^{1} 
        t^{s-1} e^{-tP}dt. 
\end{equation}
    Since $\Gamma(s+1)=s\Gamma(s)$ we get 
\begin{equation}
        Q_{s}= P Q_{s+1} - \Gamma(s+1)^{-1}e^{-P}, \qquad \Re s>0.
    \label{eq:Zeta.Qs-extension1}
\end{equation}
An easy induction then shows that for $k=1,2,\ldots$ we have 
\begin{equation}
     Q_{s}= P^{k} Q_{s+k} - \Gamma(s+k)^{-1}P^{k-1}e^{-P}+\ldots +(-1)^{k}\Gamma(s+1)^{-1}e^{-P}.
     \label{eq:Zeta.Qs-extension2}
\end{equation}
It follows that the family $(Q_{s})_{\Re s>0}$ has a unique analytic continuation to each half-space $\Re s>-k$ for $k=1,2,\ldots$, so it admits a unique analytic 
continuation to $\C$. Furthermore, as for $\Re s>-k$ we have $P^{-s}_{\updown}=P^{k}P^{-(s+k)}_{\updown}$  we get 
\begin{equation}
    Q_{s}-P^{-s}_{\updown}=P^{k} R_{\updown}(s+k) -\Gamma(s+k)^{-1}P^{k-1}e^{-P}+\ldots +(-1)^{k}\Gamma(s+1)^{-1}e^{-P},
\end{equation}
from which we deduce that $(Q_{s}-P^{-s}_{\updown})_{\Re s >-k}$ is a holomorphic family of smoothing operators. Hence 
the families $(Q_{s})_{s \in \C}$ and $(P^{-s}_{\updown})_{s \in \C}$ agree up to a holomorphic family of smoothing operators. 

Finally, we have 
\begin{equation}
    Q_{1}=\Pi_{0}(P)+\int_{0}^{1}(1-\Pi_{0}(P))e^{-tP}dt=\Pi_{0}(P)-P^{-1}(e^{-P}-1).
\end{equation}
Thus setting $s=1$ in~(\ref{eq:Zeta.Qs-extension1}) gives
\begin{multline}
    Q_{0}=P[\Pi_{0}(P)-P^{-1}(e^{-P}-1)]+e^{-P} =  -(1-\Pi_{0}(P))(e^{-P}-1)+e^{-P}\\ = 1-\Pi_{0}(P)+\Pi_{0}e^{-P}=1.
     \label{eq:Zeta.Q0}
\end{multline}
Furthermore, as $\Gamma(s)^{-1}$ vanishes at every non-positive integer, from~(\ref{eq:Zeta.Qs-extension2}) and~(\ref{eq:Zeta.Q0}) we see that we have 
$Q_{-k}=P^{k}Q_{0}=P^{k}$ for any integer $k\geq 1$. 
The proof of the claim is thus achieved. 
\end{proof}

Now, for $j=0,1,\ldots$ we set $\sigma_{j}=\frac{d+2-j}{m}$. As $(R_{\updown}(s))_{s \in \C}:=(P_{\updown}^{-s}-Q_{s})_{s \in \C}$ 
is a holomorphic family of smoothing operators, the map $s \rightarrow 
t_{R_{\updown}(s)}(x)$ is holomorphic from $\C$ to $C^{\infty}(M, |\Lambda|(M)\otimes \End \cE)$. By combining this with 
Proposition~\ref{prop:Zeta.zeta-function} we deduce that for $j=0,1,\ldots$ we have
\begin{equation}
     \res_{s=\sigma_{j}}t_{P_{\updown}^{-s}}(x)=m c_{P^{-\sigma_{j}}}(x) =  \res_{s=\sigma_{j}}t_{Q_{s}}(x),
     \label{eq:Zeta.tPs-TQs'1}
\end{equation}
Moreover, as for $k=1,2,\ldots$ we have $R_{\updown}(-k)=0$ we also see that 
\begin{multline}
  \lim_{s\rightarrow -k}[t_{P_{\updown}^{-s}}(x)-m (s+k)^{-1}c_{P^{k}}(x)] \\ =     \lim_{s\rightarrow -k}[t_{Q_{s}}(x)-(s+k)^{-1}\res_{s=-k}t_{Q_{s}}(x)].
\end{multline}
Similarly, as $P_{\updown}^{0}=1-\Pi_{0}(P)=Q_{0}-\Pi_{0}(P)$ we get
\begin{equation}
       \lim_{s\rightarrow 0} t_{P_{\updown}^{-s}}(x) =\lim_{s\rightarrow 0}t_{Q_{s}}(x)-t_{\Pi_{0}}(x).
           \label{eq:Zeta.tPs-TQs'3}
 \end{equation}

 Next, let $k_{Q_{s}}(x,y)$ denote the kernel of $Q_{s}$. As $Q_{s}$ has order $-m s$, for $\Re s>-\frac{d+2}{m}$ this is a trace-class operator 
 and thanks to~(\ref{eq:Zeta.heat-kernel-asymptotics}) we have 
 \begin{equation}
     \Gamma(s) k_{Q_{s}}(x,x)= \int_{0}^{1}t^{s-1}k_{t}(x,x) dt.
 \end{equation}
 Moreover~(\ref{eq:Zeta.heat-kernel-asymptotics}) 
 implies that, for any integer $N\geq 0$, in $C^{\infty}(M,\End \cE \otimes |\Lambda|(M))$ we have
 \begin{equation}
     k_{t}(x,x)=\sum_{-\sigma_{j}<N}t^{-\sigma_{j}}a_{j}(P)(x)+\sum_{k<N}(t^{k}\log t)b_{k}(P)(x) +\op{O}(t^{N}). 
 \end{equation}
Therefore, for $\Re s >\frac{d+2}{m}$ the density  $\Gamma(s)k_{Q_{s}}(x,x)$ is of the form
\begin{equation}
       \sum_{\sigma_{j}<N}(\int_{0}^{1}t^{s-\sigma_{j}}\frac{dt}{t})a_{j}(P)(x)+ 
      \sum_{k<N}(\int_{0}^{1}t^{k+s}\log t \frac{dt}{t})b_{k}(P)(x) + \Gamma(s) h_{N,s}(x),
\end{equation}
 with $h_{N,s}(x) \in \Hol(\Re s>-N, C^{\infty}(M,\End \cE \otimes |\Lambda|(M))$. Since for $\alpha >0$ we have 
 \begin{equation}
     \int_{0}^{1}t^{\alpha}\log t \frac{dt}{t}=-\frac{1}{\alpha}\int_{0}^{1}t^{\alpha-1}dt = -\frac{1}{\alpha},
 \end{equation}
 we see that $k_{Q_{s}}(x,x)$ is equal to
 \begin{equation}
      \Gamma(s)^{-1} \sum_{\sigma_{j}<N}\frac{1}{s+\sigma_{j}}a_{j}(P)(x)-
     \Gamma(s)^{-1}  \sum_{k<N}\frac{1}{(s+k)^{2}}b_{k}(P)(x) + h_{N,s}(x). 
 \end{equation}
 Since $\Gamma(s)$ is analytic on $\C\setminus (\Z_{-}\cup\{0\})$ and for $k=0,1,\ldots$ near $s=-k$ we have 
 $\Gamma(s)^{-1}\sim (-1)^{k}k!(s+k)^{-1}$ , we deduce that:\smallskip

 - when $\sigma_{j}\not \in \N$ we have $ \res_{s=\sigma_{j}}t_{Q_{s}}(x)= \Gamma(\sigma_{j})^{-1}a_{j}(P)(x)$.\smallskip
 
 \indent - for $k=1,2,\ldots$ we have
 \begin{gather}
 \res_{s=-k}t_{Q_{s}}(x)= (-1)^{k+1}k!b_{k}(P)(x),\\
       \lim_{s\rightarrow -k}[t_{Q_{s}}(x)-(s+k)^{-1}\res_{s=-k}t_{Q_{s}}(x)] = (-1)^{k}k! a_{d+2+mk}(P)(x).
   \end{gather}
\indent - for $k=0$ we have $\lim_{s\rightarrow 0} t_{Q_{s}}(x) =a_{d+2}(P)(x)$.\smallskip

\noindent Combining this with~(\ref{eq:Zeta.tPs-TQs'1})--(\ref{eq:Zeta.tPs-TQs'3}) then proves 
 the equalities~(\ref{eq:Zeta.tPs-heat1})--(\ref{eq:Zeta.tPs-heat4}). 
 \end{proof}

From Proposition~\ref{prop:Zeta.heat-zeta-local} we immediately get:

\begin{proposition}\label{prop:Zeta.heat-zeta-global}
 1) For $j=0,1,\ldots$ let $\sigma_{j}=\frac{d+2-j}{m}$. When $\sigma_{j}\not \in \Z_{-}$ we have: 
 \begin{equation}
    \res_{s=\sigma_{j}}\zeta_{\updown}(P;s) =m \Res P^{-\sigma_{j}}= \Gamma(\sigma_{j})^{-1}\int_{M}\tr_{\cE}a_{j}(P)(x).
      \label{eq:Zeta.heat-zeta-global1}
 \end{equation}
 2) For $k=1,2,\ldots$ we have
 \begin{gather}
      \res_{s=-k}\zeta_{\updown}(P;s) =m \Res P^{k} = (-1)^{k+1}k!\int_{M}\tr_{\cE}b_{k}(P)(x), \\
       \lim_{s\rightarrow -k}[\zeta_{\updown}(P;s)-m (s+k)^{-1} \Res P^{k}] = (-1)^{k}k! \int_{M}\tr_{\cE} a_{d+2+mk}(P)(x).
\end{gather}
3) For $k=0$ we have 
\begin{equation}
    \zeta_{\updown}(P;0)=\int_{M}\tr_{\cE} a_{d+2}(P)(x)-\dim \ker P.
\end{equation}
\end{proposition}

Next, for $k=0,1,\ldots$ let $\lambda_{k}(P)$ denote the $(k+1)$'th eigenvalue of $P$ counted with multiplicity. 
Then by~\cite{Po:MAMS1} and \cite{Po:CPDE1} as  $k\rightarrow \infty$ we have the Weyl asymptotics, 
\begin{equation}
   \lambda_{k}(P)\sim \left(\frac{k}{\nu_{0}(P)}\right)^{\frac{m}{d+2}}, \qquad \nu_{0}(P)=\Gamma(1+\frac{d+2}{m})^{-1} \int_{M}\tr_{\cE}a_{0}(P)(x).
     \label{eq:Zeta.Weyl-asymptotics1}
\end{equation}

Now, by Proposition~\ref{prop:Zeta.heat-zeta-global} we have 
\begin{equation}
   \int_{M}\tr_{\cE}a_{0}(P)(x)=m\Gamma(\frac{d+2}{m})\Res P^{-\frac{d+2}{m}}=\frac{1}{d+2}  \Gamma(1+\frac{d+2}{m})\Res P^{-\frac{d+2}{m}},
\end{equation}
Therefore, we obtain: 

\begin{proposition}\label{prop:Zeta.Weyl-asymptotics}
    As $k\rightarrow \infty$ we have 
    \begin{equation}
           \lambda_{k}(P)\sim \left(\frac{k}{\nu_{0}(P)}\right)^{\frac{m}{d+2}}, \qquad \nu_{0}(P)= (d+2)^{-1}\Res P^{-\frac{d+2}{m}}.
    \end{equation}
\end{proposition}

Finally, we can make use of Proposition~\ref{prop:Zeta.heat-zeta-global} 
to prove a local index formula for hypoelliptic \psivdos\ in the following setting. Assume that $\cE$ admits a $\Z_{2}$-grading
$\cE=\cE^{+}\oplus \cE_{-}$ and  let $D:C^{\infty}(M,\cE)\rightarrow C^{\infty}(M,\cE)$ be a selfadjoint \psivdo\ of integer order $m\geq 1$ 
with an invertible principal 
symbol and of the form, 
\begin{equation}
    D= \left( 
    \begin{array}{cc}
        0& D_{-}  \\ 
        D_{+}& 0
    \end{array}
    \right), \qquad D_{\pm}:C^{\infty}(M,\cE_{\pm}) \rightarrow C^{\infty}(M,\cE_{\mp}). 
\end{equation}
Notice that the selfadjointness of $D$ means that $D_{+}^{*}=D_{-}$. 

Since $D$ has an invertible principal symbol and $M$ is compact we see that 
$D$ is invertible modulo finite rank operators, hence is Fredholm. 
Then we let 
\begin{equation}
    \ind D:= \ind D_{+}=\dim \ker D_{+}-\dim \ker D_{-}.
\end{equation}
\begin{proposition}
    Under the above assumptions we have
    \begin{equation}
        \ind D=\int_{M} \str_{\cE} a_{d+2}(D^{2})(x), 
    \end{equation}
where $\str_{\cE}:=\tr_{\cE^{+}}-\tr_{\cE^{-}}$ denotes the supertrace on the fibers of $\cE$.
\end{proposition}
\begin{proof}
    We have $D^{2} = \left( 
    \begin{array}{cc}
        D_{-}D_{+}& 0  \\
        0 &  D_{+}D_{-}
    \end{array}
    \right)$ and  $D_{\mp} D_{\pm}=D_{\pm}^{*}D_{\pm}$. In particular, $D_{\mp} D_{\pm}$ is a positive operators with
   an invertible principal symbol.  Moreover, for $\Re s >\frac{d+2}{2m}$ the difference $  \zeta(D_{-}D_{+};s) -\zeta(D_{+}D_{-};s) $ is equal to
\begin{equation}
     \sum_{\lambda>0} \lambda^s (\dim\ker (D_{-}D_{+} -\lambda) - 
            \dim\ker (D_{+}D_{-} -\lambda)) = 0,
\end{equation}
    for $D$ induces for any $\lambda>0$ a bijection between $\ker (D_{-}D_{+}-\lambda)$ and $\ker 
    (D_{+}D_{-} -\lambda)$ (see, e.g.,~\cite{BGV:HKDO}). By analytic continuation this yields 
    $ \zeta(D_{-}D_{+};0) -\zeta(D_{+}D_{-};0)=0$. 

On the other hand, by Proposition~\ref{prop:Zeta.heat-zeta-global} we have 
\begin{equation}
    \zeta(D_{\mp} D_{\pm};0) = \int_{M} \tr_{\cE_{\pm}}a_{d+2}(D_{\mp} D_{\pm})(x)-\dim \ker D_{\mp} D_{\pm}.
\end{equation}
Since $\dim \ker D_{\mp} D_{\pm}=\dim \ker D_{\pm}$ we deduce that $\ind D$ is equal to
\begin{equation}
     \int_{M} \tr_{\cE_{+}}a_{d+2}(D_{+}D_{-})(x)- \int_{M} \tr_{\cE_{-}}a_{d+2}(D_{-}D_{+})(x)= \int_{M} \str_{\cE}a_{d+2}(D^{2})(x).
\end{equation}
The proof is thus achieved. 
\end{proof}

\subsection{Metric estimates for Green kernels of hypoelliptic \psivdos}
Consider a compact Heisenberg manifold $(M^{d+1},H)$ endowed with a positive density and let $\cE$ be a Hermitian vector bundle over $M$. 
In this subsection we shall prove that the positivity of a hypoelliptic \psivdo\ pertains in its logarithmic singularity when it has order $-(\dim 
M+1)$. As a consequence this will allow us to derive some metric estimates for Green kernels of hypoelliptic \psidos. 

Let $P:C^{\infty}(M,\cE)\rightarrow C^{\infty}(M,\cE)$ be a \psivdo\ of order $m>0$ whose principal symbol is invertible and is positive in the sense 
of~\cite{Po:MAMS1}, i.e., we can write $\sigma_{m}(P)=q*q^{*}$ with $q\in S_{\frac{m}{2}}(\fg^{*}M,\cE)$. The main technical result of this section 
is the following. 

\begin{proposition}\label{prop:Metric.positivity-cP}
     The density $\tr_{\cE}c_{P^{-\frac{d+2}{m}}}(x)$ is~$>0$. 
\end{proposition}

We will prove Proposition~\ref{prop:Metric.positivity-cP} later on in the section. As a first consequence, by combining 
with Proposition~\ref{prop:Zeta.heat-zeta-local} we get:

\begin{proposition}
    Let $a_{0}(P)(x)$ be the leading coefficient in the small time heat kernel asymptotics~(\ref{eq:Zeta.heat-kernel-asymptotics}) for $P$. 
    Then the density $\tr_{\cE}a_{0}(P)(x)$ 
    is~$>0$.
\end{proposition}

Assume now that the bracket condition $H+[H,H]=TM$ holds, i.e., $H$ is a Carnot-Carath\'eodory distribution in the sense of~\cite{Gr:CCSSW}. Let $g$ be a 
Riemannian metric on $H$ and let $d_{H}(x,y)$ be the associated Carnot-Carath\'eodory metric on $M$. Recall that for two points $x$ 
and $y$ of $M$ the value of $d_{H}(x,y)$ is the infinum of the lengths of all closed paths joining $x$ to $y$ that are tangent to $H$ at each point 
(such a path always exists by Chow Lemma). Moreover, the Hausdorff dimension of $M$ with respect to $d_{H}$ is equal to $\dim M+1$. 

In the setting of general Carath\'eodory distributions there has been lot of interest by Fefferman, Stein and their collaborators for giving 
metric estimates for the singularities of the Green kernels of hypoelliptic sublaplacians 
(see, e.g., \cite{FS:FSSOSO}, \cite{Ma:EPKLCD}, \cite{NSW:BMDVF1}, \cite{Sa:FSGSSVF}). 
This allows us relate the analysis of the hypoelliptic sublaplacian to the metric geometry of the underlying manifold. 

An important result is that it follows from the maximum principle of Bony~\cite{Bo:PMIHUPCOED} that 
the Green of kernel of a selfadjoint hypoelliptic sublaplacian is positive near the diagonal. 
In general the positivity of the principal symbol does not pertain in the Green kernel. However, by making use of Proposition~\ref{prop:Metric.positivity-cP} 
we shall prove: 

\begin{theorem}\label{thm:Metric.metric-estimate}
  Assume that $H+[H,H]=TM$ and let $P:C^{\infty}(M)\rightarrow C^{\infty}(M)$ be a \psivdo\ of order $m>0$ whose principal symbol is 
  invertible and is positive. Let $k_{P^{-\frac{d+2}{m}}}(x,y)$  be the Schwartz kernel 
  of $P^{-\frac{d+2}{m}}$. Then near the diagonal we have 
  \begin{equation}
      k_{P^{-\frac{d+2}{m}}}(x,y)\sim -c_{P^{-\frac{d+2}{m}}}(x)\log d_{H}(x,y).
       \label{eq:Metric.metric-estimate}
  \end{equation}
  In particular $k_{P^{-\frac{d+2}{m}}}(x,y)$ is~$>0$ near the diagonal.  
\end{theorem}
\begin{proof}
  It is enough to proceed in an open of $H$-framed local coordinates $U\subset \Rd$. For $x \in U$ let $\psi_{x}$ be the affine change to the corresponding privileged 
  coordinates at $x$. Since by Proposition~\ref{prop:Metric.positivity-cP} we have $c_{P^{-\frac{d+2}{m}}}(x)>0$, 
  using Proposition~\ref{thm:NCR.log-singularity} we see that near the diagonal we have  
  $ k_{P^{-\frac{d+2}{m}}}(x,y) \sim -c_{P^{-\frac{d+2}{m}}}(x)\log\|\psi_{x}(y)\|$. 
Incidentally, we see that $k_{P^{-\frac{d+2}{m}}}(x,y)$ is positive near the diagonal.  
  
  On the other hand, since $H$ has codimension one our definition of the privileged coordinates agrees with that of~\cite{Be:TSSRG}. Therefore,  
  it follows from~\cite[Thm.~7.34]{Be:TSSRG} that the ratio $\frac{d_{H}(x,y)}{\|\psi_{x}(y)\|}$ remains bounded in $(0, \infty)$ near the diagonal, that is, 
  we have $\log d_{H}(x,y)\sim \log \|\psi_{x}(y)\|$. 
  It then follows that near the diagonal we have $ k_{P^{-\frac{d+2}{m}}}(x,y)\sim -c_{P^{-\frac{d+2}{m}}}(x)\log d_{H}(x,y)$. The theorem is thus 
  proved.
\end{proof}

It remains now to prove Proposition~\ref{prop:Metric.positivity-cP}. To this end recall that for an operator $Q\in \pvdo^{l}(M,\cE)$, $l\in \C$, 
the model operator $Q^{a}$ at a given point $a\in M$ is defined as the left-invariant \psivdo\ on $\cS_{0}(G_{a}M,\cE)$ with symbol 
$q^{a}(\xi)=\sigma_{l}(Q)(a,\xi)$. 
Bearing this in mind we have:


\begin{lemma}\label{lem:Metric.cP-Heisenberg-coordinates}
    Let $Q\in \pvdo^{-(d+2)}(M,\cE)$ and let $Q^{a}$ be its model operator at a point $a \in M$.\smallskip 
    
   1) We have $c_{Q^{a}}(x)=c_{Q^{a}}dx$, where $c_{Q^{a}}$ is a constant and $dx$ denotes the Haar measure of $G_{a}M$.\smallskip
    
   2) In Heisenberg coordinates centered at $a$ we have $c_{Q}(0)=c_{Q^{a}}$.
\end{lemma}
\begin{proof}
    Let $X_{0},X_{1},\ldots,X_{d}$ be a $H$ frame near $a$. Since $G_{a}M$ has underlying set $(T_{a}M/H_{a})\oplus H_{a}$ the 
    vectors $X_{0}(a),\ldots,X_{d}(a)$ define global coordinates for $G_{a}M$, so that we can identify it with $\Rd$ equipped with the group 
    law~(\ref{eq:Heisenberg.group-law-tangent-group-coordinates}). 
    In these coordinates set $q^{a}(\xi):=\sigma_{-(d+2)}(P)(a,\xi)$. Then~(\ref{eq:PsiHDO.PsiDO-convolution}) tells us that $Q^{a}$ corresponds to 
    the operator $q^{a}(-iX^{a})$ acting on $\cS_{0}(\Rd)$, where $X_{0}^{a},\ldots,X_{d}^{a}$ is the left-invariant tangent frame 
   coming from the model vector fields at $a$ of $X_{0},\ldots,X_{d}$. 
   
   Notice that the left-invariance of the frame $X_{0}^{a},\ldots,X_{d}^{a}$ implies that, with respect to this frame, 
 the affine change of variables to the privileged 
coordinates centered at any given point $x\in \Rd$ is just $\psi^{a}_{x}(y)=y.x^{-1}$. In view of~(\ref{eq:Heisenberg.group-law-tangent-group-coordinates}) 
this implies that $|\psi_{x}^{a'}|=1$. Therefore, from~(\ref{eq:NCR.formula-cP}) we get
\begin{equation}
    c_{Q^{a}}(x)=(2\pi)^{-(d+1)}\int_{\|\xi\|=1}q^{a}(\xi)\iota_{E}d\xi.
     \label{eq:Metric.cQa}
\end{equation}
Since the Haar measure of $G_{a}M$ corresponds to the Lebesgue measure of $\Rd$ this proves the 1st part of the lemma. 

Next, by Definition~\ref{def:Heisenberg.principal-symbol} 
in Heisenberg coordinates centered at $a$ the principal symbol $\sigma_{-(d+2)}(Q)(x,\xi)$ agrees at $x=0$ 
with the principal symbol 
$q_{-(d+2)}(x,\xi)$ of $Q$ in the sense of~(\ref{eq:Heisenberg.asymptotic-expansion-symbols}), so we have $q^{a}(\xi)=q_{-(d+2)}(0,\xi)$. 
Furthermore, as we already are in Heisenberg coordinates, 
hence in privileged coordinates, we see that, with respect to the $H$-frame $X_{0},\ldots,X_{d}$, the affine change of variables 
$\psi_{0}$ to the privileged coordinates centered at the origin is just the identity. 
Therefore, by using~(\ref{eq:NCR.formula-cP}) and~(\ref{eq:Metric.cQa}) we see that $c_{Q}(0)$ is equal to
\begin{equation}
    (2\pi)^{-(d+1)}\int_{\|\xi\|=1}q(0,\xi)\iota_{E}d\xi=(2\pi)^{-(d+1)}\int_{\|\xi\|=1}q^{a}(\xi)\iota_{E}d\xi=c_{Q^{a}}.
\end{equation}
The 2nd part of the lemma is thus proved. 
\end{proof}

We are now ready to prove Proposition~\ref{prop:Metric.positivity-cP}. 

\begin{proof}[Proof of Proposition~\ref{prop:Metric.positivity-cP}]
     For sake of simplicity we may assume that $\cE$ is the trivial line bundle, since in the general case the proof follows 
     along similar lines. Moreover, for any $a \in M$ by 
     Lemma~\ref{lem:Metric.cP-Heisenberg-coordinates}  in Heisenberg coordinates centered at $a$ we have 
     $c_{P^{-\frac{d+2}{m}}}(0)=c_{(P^{-\frac{d+2}{m}})^{a}}$. Therefore, it is enough to prove that 
     $c_{(P^{-\frac{d+2}{m}})^{a}}$ is  $>0$ for any $a\in M$. 
    
    Let $a \in M$ and let $X_{0},\ldots,X_{d}$ be a $H$-frame near $a$. By using the coordinates provided by the vectors $X_{0}(a),\ldots,X_{d}(a)$ we 
    can identify $G_{a}M$ with $\Rd$ equipped with the group law~(\ref{eq:Heisenberg.group-law-tangent-group-coordinates}). We then let 
    $H^{a}\subset T\Rd$ be the hyperplane bundle spanned by the model vector fields $X_{1}^{a},\ldots,X_{d}^{a}$ seen as left-invariant vector 
    fields on $\Rd$. In addition, for any $z\in \C$ we let $p(z)(\xi):=\sigma_{z}(P^{\frac{z}{m}})(a,\xi)$ be the principal symbol at $a$ of $P^{\frac{z}{m}}$, 
    seen as a homogeneous symbol on $\Rdo$. Notice that by~\cite[Rem.~4.2.2]{Po:MAMS1} the family $(p(z))_{z \in \C}$ is a holomorphic family with values in 
    $C^{\infty}(\Rdo)$. 
    
    Let $\chi\in C^{\infty}_{c}(\Rd)$ be such that $\chi(\xi)=1$ near $\xi=0$. For any $z \in \C$ and for any pair  $\varphi$ and $\psi$  of functions in 
    $C_{c}^{\infty}(\Rd)$ we set
    \begin{equation}
        \tilde{p}(z)(\xi):=(1-\chi)p(z) \qquad \text{and} \qquad P_{\varphi,\psi}(z):=\varphi \tilde{p}(z)(-iX^{a}) \psi.
    \end{equation}
    Then $(\tilde{p}(z))_{z\in \C}$ and $(P_{\varphi,\psi}(z))_{z\in \C}$ are holomorphic families with values in  $S^{*}(\Rd)$ and $\Psi^{*}_{H^{a}}(\Rd)$ 
    respectively. 
    
    Notice that $P_{\varphi,\psi}(z)$ has order $z$ and the support of its Schwartz kernel is contained in the fixed 
    compact set $\supp \varphi \times \supp \psi$, so by Proposition~\ref{prop:Heisenberg.L2-boundedness} 
    the operator $P_{\varphi,\psi}(z)$ is bounded on $L^{2}(\Rd)$ for $\Re  z\leq 0$. In fact, by arguing as in the proof of~\cite[Prop.~4.6.2]{Po:MAMS1} 
    we can show that $(P_{\varphi,\psi}(z))_{\Re z\leq 0}$ 
    actually is a holomorphic family with values in $\cL(L^{2}(\Rd))$. 
    
     Moreover, by~\cite[Prop.~4.6.2]{Po:MAMS1} the family $(P_{\varphi,\psi}(\overline{z})^{*})_{z\in \C}$ is a holomorphic family with values in 
     $\Psi^{*}_{H^{a}}(\Rd)$ such that $\ord P_{\varphi,\psi}(\overline{z})^{*}=z$ for any $z \in \C$. Therefore
     $(P_{\varphi,\psi}(z)P_{\varphi,\psi}(\overline{z})^{*})_{\Re z <-\frac{d+2}{2}}$ is a holomorphic family 
     with values in $\Psi_{H^{a}}^{\text{int}}(\Rd)$. For any $z \in \C$ let $k(z)(x,y)$ denote the Schwartz kernel 
    of $P_{\varphi,\psi}(z)P_{\varphi,\psi}(\overline{z})^{*}$. Then the support of $k(z)(x,y)$ is contained in the fixed compact set $\op{supp}\varphi \times 
    \op{supp}\varphi$, and by using~\ref{eq:NCR.kP(x,x)} we can check that $(k(z)(x,y))_{\Re z <-\frac{d+2}{2}}$ is a holomorphic family of 
    continuous Schwartz kernels. It then follows that $(P_{\varphi,\psi}(z)P_{\varphi,\psi}(\overline{z})^{*})_{\Re z <-\frac{d+2}{2}}$ 
    is a holomorphic family with values in the Banach ideal $\cL^{1}(L^{2}(\Rd))$ of trace-class operators on $L^{2}(\Rd)$.  
    
    Let us now choose $\psi$ so that $\psi =1$ near $\supp \varphi$. For any $t \in \R$ the operator 
    $P^{\frac{t}{m}}$ is selfadjoint, so by Proposition~\ref{prop:Heisenberg.operations-principal-symbols} its principal symbol is real-valued. 
    Therefore, by Proposition~\ref{prop:Heisenberg.operations-principal-symbols}  the principal symbol of 
    $(P_{\varphi,\psi}(t)P_{\varphi,\psi}(t)^{*})$ is equal to 
\begin{equation}
        [ \varphi p(t)\psi]*^{a}[\overline{\psi}\overline{p(t)}\overline{\varphi}] =|\varphi|^{2}p(t)*p(t)=|\varphi|^{2}p(2t).
     \label{eq:Metric.principal-symbol}
\end{equation}
    In particular, the principal symbols of $P_{\varphi,\psi}(-\frac{d+2}{2})P_{\varphi,\psi}(-\frac{d+2}{2})^{*}$ and 
    $P_{|\varphi|^{2},\psi}(-(d+2))$ agree. By combining this with Lemma~\ref{lem:Metric.cP-Heisenberg-coordinates} we see that
     \begin{equation}
        c_{P_{\varphi,\psi}(-\frac{d+2}{2})P_{\varphi,\psi}(-\frac{d+2}{2})^{*}}(x)= 
            c_{P_{|\varphi|^{2},\psi}(-(d+2))}(x)=|\varphi(x)|^{2}c_{(P^{-\frac{(d+2)}{m}})^{a}}.
         \label{eq:Metric.cPphipsi-cPa}
    \end{equation}
    It then follows from Proposition~\ref{thm:NCR.TR.local} that we have:
    \begin{multline}
        c_{(P^{-\frac{(d+2)}{m}})^{a}}(\int |\varphi(x)|^{2}dx) =\int c_{P_{\varphi,\psi}(-\frac{d+2}{2})P_{\varphi,\psi}(-\frac{d+2}{2})^{*}}(x) dx \\
        =\lim_{t \rightarrow \frac{-(d+2)}{2}}\frac{-1}{t+\frac{d+2}{2}}\int t_{P_{\varphi,\psi}(t)P_{\varphi,\psi}(t)^{*})}(x) dx\\ = 
         \lim_{t \rightarrow [\frac{-(d+2)}{2}]^{-}}\frac{-1}{t+\frac{d+2}{2}} \Tra [P_{\varphi,\psi}(-t)P_{\varphi,\psi}(-t)^{*}] \geq 0.
    \end{multline}
    Thus, by choosing $\varphi$ so that $\int |\varphi|^{2}>0$ we obtain that $c_{(P^{-\frac{(d+2)}{m}})^{a}}$ is $\geq 0$. 
     
     Assume now that $c_{(P^{-\frac{(d+2)}{m}})^{a}}$ vanishes, and let us show that this assumption 
     leads us to a contradiction. Observe that 
     $(P_{\varphi,\psi}(\frac{z-(d+2)}{2})P_{\varphi,\psi}(\frac{\overline{z-(d+2)}}{2})^{*})_{z\in \C}$ is  holomorphic gauging for 
     $P_{\varphi,\psi}(-\frac{d+2}{2})P_{\varphi,\psi}(-\frac{d+2}{2})^{*}$. Moreover, by~(\ref{eq:Metric.cPphipsi-cPa}) we have    
     $c_{P_{\varphi,\psi}(-\frac{d+2}{2})P_{\varphi,\psi}(-\frac{d+2}{2})^{*}}(x)=|\varphi(x)|^{2}c_{(P^{-\frac{(d+2)}{m}})^{a}}=0$.
     Therefore, it follows from Proposition~\ref{prop:Heisenberg.operations-principal-symbols} that 
     $\TR P_{\varphi,\psi}(z)P_{\varphi,\psi}(\overline{z})^{*}$ is analytic near $z=-\frac{d+2}{2}$. In particular, the limit
     $\lim_{t\rightarrow {\frac{-(d+2)}{2}}^{-}} \Tra P_{\varphi,\psi}(t)P_{\varphi,\psi}(t)^{*}$ exists and is finite. 
     
     Let $(\xi_{k})_{k\geq 0}$ be an orthonormal basis 
     of $L^{2}(\Rd)$ and let $N\in \N$. For any $t>\frac{d+2}{2}$ the operator $P_{\varphi,\psi}(t)P_{\varphi,\psi}(t)^{*}$ is trace-class and we have 
     \begin{equation}
         \sum_{0\leq k \leq N}\acou{P_{\varphi,\psi}(t)P_{\varphi,\psi}(t)^{*}\xi_{k}}{\xi_{k}} \leq \Tra [P_{\varphi,\psi}(t)P_{\varphi,\psi}(t)^{*}].
          \label{eq:Metric.partial-trace}
     \end{equation}
      As $t \rightarrow {-\frac{d+2}{2}}^{-}$ the operator $P_{\varphi,\psi}(t)P_{\varphi,\psi}(t)^{*}$ converges to 
      $P_{\varphi,\psi}(-\frac{d+2}{2})P_{\varphi,\psi}(-\frac{d+2}{2})^{*}$ in $\cL(L^{2}(\Rd)$. Therefore, letting $t$ go to ${-\frac{d+2}{2}}^{-}$ 
      in~(\ref{eq:Metric.partial-trace}) shows that, for any $N\in \N$, we have
    \begin{equation}
       \sum_{0\leq k \leq N}\acou{P_{\varphi,\psi}(-\frac{d+2}{2})P_{\varphi,\psi}(-\frac{d+2}{2})^{*}\xi_{k}}{\xi_{k}} 
        \leq\lim_{t\rightarrow [\frac{-(d+2)}{2}]^{-}} \Tra 
        [P_{\varphi,\psi}(t)P_{\varphi,\psi}(t)^{*}] <\infty.
    \end{equation}
    This proves that  $P_{\varphi,\psi}(-\frac{d+2}{2})P_{\varphi,\psi}(-\frac{d+2}{2})^{*}$ is a trace-class operator. Incidentally, we see that  
    $P_{\varphi,\psi}(-\frac{d+2}{2})$ is a Hilbert-Schmidt operator on $L^{2}(\Rd)$. 
    
    Next, let $Q \in \Psi_{H^{a}}^{-\frac{d+2}{2}}(\Rd)$ and let $q(x,\xi)\in S_{-\frac{d+2}{2}}(\Rd\times \Rd)$ be the principal symbol of $Q$. The 
    principal symbol of $\varphi Q\psi$  is $\varphi(x) q(x,\xi)$. Moreover, since for any $z \in \C$ we have $p(z)*p(-z)=p(0)=1$, 
    we see that the principal symbol of $\psi Q\psi P_{\psi,\psi}(\frac{d+2}{2})P_{\varphi,\psi}(-\frac{d+2}{2})$ is equal to
    \begin{equation}
        (\psi q \psi)*(\psi p(\frac{d+2}{2})\psi)*(\varphi p(-\frac{d+2}{2})\psi)= \varphi q*p(\frac{d+2}{2})*p(-\frac{d+2}{2})=\varphi q.
    \end{equation}
    Thus $\varphi Q\psi$ and $\psi Q\psi P_{\psi,\psi}(\frac{d+2}{2})P_{\varphi,\psi}(-\frac{d+2}{2})$ have the same principal symbol. Since they 
    both have a compactly supported Schwartz kernel it follows that we can write 
    \begin{equation}
        \varphi Q\psi = \psi Q\psi P_{\psi,\psi}(\frac{d+2}{2})P_{\varphi,\psi}(-\frac{d+2}{2}) +Q_{1},
        \label{eq:Metric.Hilbert-Schmidt-decomposition}
    \end{equation}
    for some operator $Q_{1} \in \Psi_{H^{a}}^{-\frac{d+2}{2}-1}(\Rd)$ with a compactly supported Schwartz kernel. Observe that:\smallskip 
    
    - the operator $ \psi Q\psi P_{\psi,\psi}(\frac{d+2}{2})$ is a zero'th order \psivdo\ with a compactly supported Schwartz kernel, so this is a bounded operator 
    on $L^{2}(\Rd)$;\smallskip
    
    - as above-mentioned $P_{\varphi,\psi}(-\frac{d+2}{2})$ is a Hilbert-Schmidt operator;\smallskip 
    
    - as $Q_{1}^{*}Q_{1}$ belongs to $\Psi_{H,c}^{\text{int}}(\Rd)$ this is a trace-class operator, and so $Q_{1}$ is a Hilbert-Schmidt operator.\smallskip
    
    \noindent Since the space $\cL^{2}(L^{2}(\Rd))$ of Hilbert-Schmidt operators is a two-sided ideal, it follows 
    from~(\ref{eq:Metric.Hilbert-Schmidt-decomposition}) and the above observations that $\varphi Q\psi $ is a 
    Hilbert-Schmidt operator. In particular, by \cite[p.~109]{GK:ITLNSO} the Schwartz kernel of $\varphi Q\psi $ lies in $L^{2}(\Rd\times \Rd)$. 
    
    We now get a contradiction as follows. Let $Q\in \Psi_{H^{a}}^{-\frac{d+2}{2}}(\Rd)$ have Schwartz kernel, 
    \begin{equation}
        k_{Q}(x,y)=|{\psi_{x}^{a}}'| \|\psi_{x}^{a}(y)\|^{-\frac{d+2}{2}},
    \end{equation}
    where $\psi_{x}^{a}$ is the change to the privileged coordinates at $a$ with respect to the $H^{a}$-frame $X_{0}^{a},\ldots,X_{d}^{a}$ 
    (this makes sense since $\|y\|^{-\frac{d+2}{2}}$ is in $\cK_{-\frac{d+2}{2}}(\Rd\times 
    \Rd)$). As alluded to in the proof of Lemma~\ref{lem:Metric.cP-Heisenberg-coordinates} 
    the left-invariance of the frame $X_{0}^{a},\ldots,X_{d}^{a}$ implies that 
    $\psi_{x}^{a}(y)=y.x^{-1}$. Therefore, the Schwartz kernel of $\varphi Q\psi $ is equal to
    \begin{equation}
        k_{\varphi Q\psi}(x,y)=\varphi(x) \|y.x^{-1}\|^{-\frac{d+2}{2}} \psi(y).
    \end{equation}
     However, this not an $L^{2}$-integrable kernel, since
    $\|y.x^{-1}\|^{-(d+2)}$ is not locally integrable near the diagonal.
    
    We have obtained a contradiction, so $c_{(P^{-\frac{d+2}{m}})^{a}}$ cannot be zero. Since we know 
    that $c_{(P^{-\frac{d+2}{m}})^{a}}$ is $\geq 0$, we see that 
$c_{(P^{-\frac{d+2}{m}})^{a}}$ is $>0$. The proof of Proposition~\ref{prop:Metric.positivity-cP} is thus complete.
\end{proof}

\subsection{The Dixmier trace of \psivdos}
\label{sec:Dixmier}
The quantized calculus of Connes~\cite{Co:NCG} allows us to translate into the language of quantum mechanics the main tools of the classical 
infinitesimal calculus. In particular, an important device is the Dixmier trace~(\cite{Di:ETNN}, \cite[Appendix A]{CM:LIFNCG}), 
which is the noncommutative analogue of the standard integral. 
We shall now show that, as in the case of classical \psidos\ (see~\cite{Co:AFNG}),  the noncommutative residue allows us to extend the 
Dixmier trace to the whole algebra of integer order \psivdos. 

Let us first recall the main facts about Connes' quantized calculus and the Dixmier trace. 
The general setting is that of bounded operators on a separable Hilbert space $\cH$. Extending the well known correspondence in quantum mechanics 
between variables and operators, we get the following dictionary between classical notions of infinitesimal calculus and their operator theoretic 
analogues.  
\begin{center}
    \begin{tabular}{c|c}  
        Classical & Quantum \\ \hline
        
      Real variable &  Selfadjoint operator on $\cH $ \\  
       Complex variable & Operator on $\cH $ \\
 Infinitesimal variable & Compact operator on $\cH $ \\
       Infinitesimal of order  $\alpha>0$ & Compact operator $T$ such that \\ 
                       &  $\mu_{n}(T)=\op{O}(n^{-\alpha})$
    \end{tabular}
\end{center}

The third line can be explained as follows. We cannot say that an operator $T$ is an infinitesimal by requiring that $\|T\| \leq \epsilon$ for any 
$\epsilon >0$, for this would give $T=0$. Nevertheless, we can relax this condition by requiring that for any $\epsilon>0$ we have $\|T\|<\epsilon$ 
outside a finite dimensional space. This means that $T$ is in the closure of finite rank operators, i.e., $T$ belongs to the ideal $\cK$ of compact 
operators on $\cH$. 

In the last line $\mu_{n}(T)$ denotes the  $(n+1)$'th characteristic value of $T$, i.e., the $(n+1)$'th eigenvalue 
of $|T|=(T^{*}T)^{\frac12}$. In particular, by the min-max principle  we have 
\begin{eqnarray}
    \mu_{n}(T) & = & \inf\{ \|T_{E^\perp}\|; \dim E=n\}, 
    \nonumber \\
     & = & \op{dist}(T,\mathcal{R}_{n}) , \qquad \mathcal{R}_{n}= 
   \{\text{operators of rank}\leq n\}, 
   \label{eq:NCG.min-max}
\end{eqnarray}
so the decay of $\mu_{n}(T)$ controls the accuracy of the approximation of $T$ by finite rank operators. Moreover, by using~(\ref{eq:NCG.min-max}) we 
also can check that,  for $S$, $T$ in $\cK$ and $A$, $B$ in $\cL(\cH)$, we have
 \begin{equation}
     \mu_{n}(T+S)\leq \mu_{n}(T)+\mu_{n}(S) \qquad \text{and} \qquad \mu_{n}(ATB)\leq \|A\| \mu_{n}(T) \|B\|,
     \label{eq:NCG.inequalities-2sided-ideals}
 \end{equation}
 This implies that the set of infinitesimal operators of order $\alpha$ is a two-sided ideal of $\cL(\cH)$.  
 
 Next, in this setting the analogue of the integral is provided by the Dixmier trace~(\cite{Di:ETNN}, \cite[Appendix A]{CM:LIFNCG}). The latter arises
 in the study of the logarithmic divergency of the partial traces, 
\begin{equation}
    \Trace_{N}(T) = \sum_{n=0}^{N- 1} \mu_{n}(T), \qquad T \in\cK, \quad T\geq 0.  
\end{equation}
 The domain of the Dixmier trace is the Schatten ideal, 
\begin{equation}
    \lunf=\{T\in \cK;     \|T\|_{1,\infty} :=\sup \frac{\sigma_{N}(T)}{\log N} < \infty\}. 
\end{equation}

We extend the definition of $ \Trace_{N}(T)$ by means of the
interpolation formula, 
\begin{equation}
    \sigma_{\lambda}(T) =\inf \{\|x\|_{1}+ \lambda \|y\| ; x+y=T\}, \qquad \lambda>0,
\end{equation}
where $\|x\|_{1}:=\Tra|x|$ denotes the Banach norm of the ideal $\cL^{1}$ of trace-class operators. For any integer $N$ we have 
$\sigma_{N}(T)=\Trace_{N}(T)$. 
In addition, the Ces\={a}ro mean of $\sigma_{\lambda}(T)$ with respect to the Haar measure $\frac{d\lambda}{\lambda}$ of $\R_{+}^{*}$ is
\begin{equation}
    \tau_{\Lambda}(T) = \frac{1}{\log\Lambda}\int_{e}^\Lambda \frac{\sigma_{\lambda}(T)}
    {\log  \lambda}\frac{d\lambda}{\lambda}, \qquad \Lambda\geq e. 
\end{equation}

Let $\cL(\cH)_{+}=\{T \in \cL(\cH); \ T\geq 0\}$. Then by~\cite[Appendix A]{CM:LIFNCG} 
for $T_{1}$ and $T_{2}$ in $\lunf\cap \cL(\cH)_{+}$ we have
\begin{equation}
    |\tau_{\Lambda}(T_{1}+T_{2}) -\tau_{\Lambda}(T_{1}) -\tau_{\Lambda}(T_{2}) | \leq 
    3(\norminf{T_{1}}+\norminf{T_{2}}) \frac{\log\log\Lambda}{\log\Lambda}. 
\end{equation}
Therefore, the functionals $\tau_{\Lambda}$, $\Lambda \geq e$, give rise to an additive homogeneous map,  
\begin{equation}
 \tau:   \cL^{(1,\infty)}\cap \cL(\cH)_{+} \longrightarrow C_{b}[e,\infty)/C_{0}[e,\infty).
\end{equation}
 It follows from this that for any state $\omega$ on the $C^{*}$-algebra 
$C_{b}[e,\infty)/C_{0}[e,\infty)$, i.e., for any positive linear form such that $\omega(1)=1$, there is a unique linear functional 
$\Trw:\cL^{(1,\infty)}\rightarrow \C$ such that 
\begin{equation}
    \Trw T = \omega(\tau(T)) \qquad \forall T \in  \cL^{(1,\infty)}\cap \cL(\cH)_{+}.  
\end{equation}

We gather the main properties of this functional in the following. 

\begin{proposition}[\cite{Di:ETNN}, \cite{CM:LIFNCG}]\label{prop:NCG.properties-Dixmier-trace}
    For any state $\omega$ on $C_{b}[e,\infty)/C_{0}[e,\infty)$ the Dixmier trace $\Trw$ has the following properties:\smallskip  
 
1) If $T$ is trace-class, then $\Trw T=0$.\smallskip 

2) We have $\Trw (T)\geq 0$ for any $T\in \lunf \cap \cL(\cH)_{+}$.\smallskip 

3)  If $S:\cH'\rightarrow \cH$ is a topological isomorphism, then we have 
$\Tr_{\omega,\cH'}(T)=\Tr_{\omega,\cH}(STS^{-1})$ for any $T\in \lunf(\cH')$. In particular, $\Trw$ does not depend on choice of the inner product on 
$\cH$.\smallskip

4) We have $\Trw AT=\Trw TA$ for any $A \in \cL(\cH)$ and any $T\in \lunf$, that is, $\Trw$ is a trace on the ideal $\lunf$.
\end{proposition}

The functional $\Trw$ is called the \emph{Dixmier trace} associated to $\omega$. We also say  that an operator 
$T \in \cL^{(1,\infty)}$ is \emph{measurable} when the value of $\Trw T$ is independent of the choice of the state $\omega$. We 
    then call \emph{the Dixmier trace} of $T$ the common value, 
    \begin{equation}
        \bint T :=\Trw T. 
    \end{equation}
 In addition, we let $\cM$ denote the space of measurable operators. For instance, if $T\in \cK\cap \cL(\cH)_{+}$ is such that 
 $\lim_{N\rightarrow \infty}\frac{1}{\log N} \sum_{n=0}^{N-1} \mu_{n}(T) = L$, then it can be shown that $T$ is  measurable and we have $\bint T=L$.

An important example of measurable operator is due to Connes~\cite{Co:AFNG}. Let $\cH$ be the Hilbert space $L^{2}(M,\cE)$ of 
$L^{2}$-sections of a Hermitian vector bundle over 
a compact manifold $M$ equipped with a smooth positive density and let $P:L^{2}(M,\cE)\rightarrow 
L^{2}(M,\cE)$ be a classical \psido\ of order $-\dim M$. Then $P$ is measurable for the Dixmier trace and we have 
\begin{equation}
    \bint P =\frac{1}{\dim M} \Res P,
     \label{eq:NCG.Trw-NCR-PsiDOs}
\end{equation}
where $\Res P$ denotes the noncommutative residue trace for classical \psidos\ of Wodzicki~(\cite{Wo:LISA}, \cite{Wo:NCRF}) 
and Guillemin~\cite{Gu:NPWF}.  
This allows us to extends the Dixmier trace to all 
\psidos\ of integer order, hence to integrate any such \psido\ even though it is not an infinitesimal of order~$\leq 1$. 

From now one we let $(M^{d+1},H)$ be a compact Heisenberg manifold equipped with a smooth positive density and we let $\cE$ be a Hermitian vector bundle 
over $M$. In addition, we recall that by Proposition~\ref{prop:Heisenberg.L2-boundedness} any $P\in \pvdo^{m}(M,\cE)$ with $\Re m\geq 0$ extends to a bounded operator 
from $L^{2}(M,\cE)$ to itself and this 
 operator is compact if we further have $\Re m<0$.  

Let $P:C^{\infty}(M,\cE) \rightarrow C^{\infty}(M,\cE)$ be a positive \psivdo\ with an invertible principal symbol of order $m>0$, and for 
$k=0,1,..$ let $\lambda_{k}(P)$ denote the $(k+1)$'~th eigenvalue of $P$ counted with multiplicity. By Proposition~\ref{prop:Zeta.Weyl-asymptotics} 
when $k \rightarrow \infty$ we have
\begin{equation}
    \lambda_{k}(P) \sim (\frac{k}{\nu_{0}(P)})^{\frac{m}{d+2}}, \qquad \nu_{0}(P)= \frac{1}{d+2}\Res P^{-\frac{d+2}{m}}.
\end{equation}
It follows that for any $\sigma \in \C$ with $\Re \sigma <0$ the operator $P^{\sigma}$ is an infinitesimal operator of order $\frac{m |\Re \sigma|}{d+2}$. 
Furthermore, for $\sigma=-\frac{d+2}{m}$ using~(\ref{eq:NCG.inequalities-2sided-ideals}) we see that $P^{-\frac{d+2}{m}}$ is measurable and we have 
\begin{equation}
    \bint P^{-\frac{d+2}{m}}=\nu_{0}(P)=\frac{1}{d+2}\Res P^{-\frac{d+2}{m}}.
     \label{eq:NCG.Dixmier-trace-NCR.hypoelliptic}
\end{equation}

These results are actually true for general \psivdos, for we have: 
\begin{theorem}\label{thm:NCG.Dixmier} 
    Let $P: L^{2}(M,\cE) \rightarrow L^{2}(M,\cE)$ be a \psivdo\ order $m$ with $\Re m<0$.\smallskip  
    
     1) $P$ is an infinitesimal operator of order $(\dim M+1)^{-1}|\Re m|$.  \smallskip 
    
       2) If $\ord P=-(\dim M+1)$, then $P$ is measurable and we have 
\begin{equation}
           \bint P = \frac1{\dim M+1} \Res P. 
     \label{eq:NCG.bint-NCR}
\end{equation}
\end{theorem}
\begin{proof}
    First, let $P_{0}\in \pvdo^{1}(M,\cE)$ be a positive and invertible \psivdo\ with an invertible principal symbol 
    (e.g.~$P_{0}=(1+\Delta^{*}\Delta)^{\frac{1}{4}}$, where $\Delta$ is a hypoelliptic sublaplacian). Then $PP_{0}^{m}$ is a zeroth 
    order \psivdo. By Proposition~\ref{}  any zeroth order \psivdo\ is bounded 
    on $L^{2}(M,\cE)$ and as above-mentioned $P^{-m}_{0}$ is an infinitesimal of order $\alpha:=(\dim M+1)^{-1}|\Re m|$. 
    Since we have $P=PP_{0}^{m}.P_{0}^{-m}$ we see 
    that $P$ is the product of a bounded operator and of an infinitesimal operator of order $\alpha$. As (\ref{eq:NCG.inequalities-2sided-ideals})  
    shows that the space of infinitesimal operators of order $\alpha$ is a two-sided ideal,  it follows that $P$ is an infinitesimal of 
    order $\alpha$. In particular, if $\ord P=-(d+2)$ then $P$ is an infinitesimal of order $1$, hence is contained in $\cL^{(1,\infty)}$. 
    
  Next, let $\Trw$ be the Dixmier trace associated to a state $\omega$ on $C_{b}[e,\infty)/C_{0}[e,\infty)$, and let us prove that for any $P \in 
    \pvdo^{-(d+2)}(M,\cE)$ we have $\Trw P=\frac{1}{d+2}\Res P$. 
    
    Let $\kappa:U\rightarrow \Rd$ be a $H$-framed chart mapping onto $\Rd$ such that there is a trivialization $\tau:\cE_{|U} \rightarrow U\times 
    \C^{r}$ of $\cE$ over $U$  (as in the proof of  Theorem~\ref{thm:Traces.traces} 
    we shall call such a chart a \emph{nice $H$-framed chart}). As in Subsection~\ref{sec:traces} we shall use the subscript $c$ to denote \psivdos\ with 
    a compactly supported Schwartz kernel (e.g.~$\pvdoc^{\Z}(\Rd)$ denote the class of integer order \psivdos\ on $\Rd$ whose Schwartz kernels have 
    compact supports). Notice that if $P\in\pvdoc^{\Z}(\Rd,\C^{r})$ then the operator $\tau^{*}\kappa^{*}P$ belongs to $\pvdo^{\Z}(M,\cE)$ and the 
    support of its Schwartz  kernel is a compact subset of $U\times U$.  

  Since $P_{0}$ is a positive \psivdo\ with an invertible principal symbol, Proposition~\ref{prop:Metric.positivity-cP} 
tells us that the density $\tr_{\cE}c_{P_{0}^{-(d+2)}}(x)$ is~$>0$, so we can write $\kappa_{*}[\tr_{\cE}c_{P_{0}^{-(d+2)}}(x)_{|U}]=c_{0}(x)dx$ 
for some positive function $c_{0}\in C^{\infty}(\Rd)$. 
Then for any  $c \in C_{c}^{\infty}(\Rd)$ and any $\psi \in C_{c}^{\infty}(\Rd)$ such that $\psi=1$ near $\supp c$ we let
\begin{equation}
    P_{c,\psi}:=(\frac{c\circ \kappa}{c_{0}\circ \kappa})P_{0}^{-(d+2)} (\psi\circ \kappa).
     \label{eq:Dixmier.Pcpsi}
\end{equation}
Notice that $P_{c,\psi}$ belongs to $\pvdo^{-(d+2)}(M,\cE)$ and it 
depends on the choice $\psi$ only modulo operators in $\psinf(M,\cE)$. Since the latter are trace-class operators and the Dixmier trace 
$\Tr_{\omega}$ vanishes on such operators (cf.~Proposition~\ref{prop:NCG.properties-Dixmier-trace}), 
we see that the value of  $\Trw P_{c,\psi}$ does not depend on the choice of 
$\psi$. Therefore, we define a linear functional $L:C_{c}^{\infty}(\Rd) \rightarrow \C$ by assigning to any $c \in C_{c}^{\infty}(\Rd)$ the value
\begin{equation}
    L(c):=\Trw P_{c,\psi}, 
\end{equation}
where $\psi \in C_{c}^{\infty}(\Rd)$ is such that $\psi=1$ near $\supp c$. 

On the other hand, let $P\in \pvdoc^{-(d+2)}(U,\cE_{|U})$. Then $\tau_{*}P$ belongs to $\pvdoc^{-(d+2)}(U,\C^{r}):=\pvdoc^{-(d+2)}(U)\otimes M_{r}(\C)$. 
Set $\tau_{*}P=(P_{ij})$ and define $\tr P:=\sum P_{ii}$. In addition, for $i,j=1,\ldots,r$ let $E_{ij}\in M_{r}(\C)$ be 
the elementary matrix whose all entries are zero except that on the $i$th row and $j$th column which is equal to $1$. Then we have 
\begin{equation}
    \tau_{*}P=\frac{1}{r}(\tr P)\otimes I_{r} + \sum_{i} P_{ii}\otimes (E_{ii}-\frac{1}{r}I_{r}) +\sum_{i\neq j} P_{ij}\otimes E_{ij}.
\end{equation}
Any matrix $A \in M_{r}(\C)$ with vanishing trace is contained in the commutator space $[M_{r}(\C),M_{r}(\C)]$. Notice also that the space 
$\pvdoc^{-(d+2)}(U)\otimes [M_{r}(\C),M_{r}(\C)]$ is contained in $[\pvdoc^{0}(U,\C^{r}),\pvdoc^{-(d+2)}(U,\C^{r})]$. 
Therefore, we see that 
\begin{equation}
    P= \frac{1}{r}(\tr P)\otimes \op{id}_{\cE} \mod [\pvdoc^{0}(U,\cE_{|U}),\pvdoc^{-(d+2)}(U,\cE_{|U})].
     \label{eq:Dixmier.decomposition-P}
\end{equation}

Let us write $\kappa_{*}[\tr_{\cE}c_{P}(x)]=a_{P}(x)dx$ with $a_{P} \in C^{\infty}_{c}(\Rd)$, and let 
$\psi \in C_{c}^{\infty}(\Rd)$ be such that $\psi=1$ near $\supp a_{P}$. Then we have  
\begin{equation*}
    \kappa_{*}[c_{\tr P_{a_{P},\psi}}(x)=(\frac{a_{P}(x)}{c_{0}(x)})\psi(x)\kappa_{*}[\tr_{\cE}c_{P_{0}^{-(d+2)}}(x)]=a_{P}(x)dx=\kappa_{*}[\tr_{\cE} 
    c_{P}(x)]=\kappa_{*}[c_{\tr P}(x)].
\end{equation*}
 In other words $Q:=\tr P-\tr P_{a_{P},\psi}$ is 
an element of $\pvdoc^{-(d+2)}(U)$ such that $c_{Q}(x)=0$. By the step (i) of the proof of Lemma~\ref{lem:Traces.sum-commutators.compact} we then can  
write $\kappa_{*}Q$ in the form $\kappa_{*}Q=[\chi_{0},Q_{0}]+\ldots+[\chi_{d},Q_{d}]$ for some functions 
$\chi_{0},\ldots,\chi_{d}$ in $C^{\infty}_{c}(\Rd)$ and some operators $Q_{0},\ldots,Q_{d}$ in $\pvdoc^{\Z}(\Rd)$.  In fact, it follows from the proof 
of Lemmas~\ref{lem:Traces.sum-commutators} and~\ref{lem:Traces.sum-commutators.compact} that $Q_{0},..,Q_{d}$ can be chosen to have order~$\leq -(d+2)$. 
This insures us that $\kappa_{*}Q$ is contained 
in $[\pvdoc^{0}(\Rd),\pvdoc^{-(d+2)}(\Rd)]$. Thus, 
\begin{equation}
    \tr P= \tr P_{a_{P},\psi} \mod [\pvdoc^{0}(U),\pvdoc^{-(d+2)}(U)].
\end{equation}
By combining this with~(\ref{eq:Dixmier.decomposition-P}) we obtain 
\begin{equation*}
    P=  \frac{1}{r}(\tr P)\otimes \op{id}_{\cE} =  \frac{1}{r}(\tr P_{a_{P},\psi})\otimes \op{id}_{\cE} =P_{a_{P},\psi} \mod 
    [\pvdoc^{0}(U,\cE_{|U}),\pvdoc^{-(d+2)}(U,\cE_{|U})].
\end{equation*}
Notice that $[\pvdoc^{0}(U,\cE_{|U}),\pvdoc^{-(d+2)}(U,\cE_{|U})]$ is contained in $[\pvdo^{0}(M,\cE),\pvdo^{-(d+2)}(M,\cE)]$, which is itself contained in 
the commutator space $[\cL(L^{2}(M)),\cL^{(1,\infty)}(M)]$ of $\lunf$.  As the Dixmier trace $\Trw$ vanishes on the latter space 
(cf.~Proposition~\ref{prop:NCG.properties-Dixmier-trace}) we deduce that 
\begin{equation}
    \Trw P=\Trw P_{a_{P},\psi}=L(a_{P}). 
     \label{eq:NCG.Trw-tau-cP}
\end{equation}

Now, let $c \in C^{\infty}_{c}(\Rd)$ and set $c_{1}= 
\frac{c}{\sqrt{c_{0}(x)}}$. In addition, let $\psi \in C^{\infty}_{c}(\Rd)$ be such that $\psi\geq 0$ and $\psi=1$ near $\supp c$, and set 
$\tilde{c}_{1}=c\circ \kappa$ and $\tilde{\psi}=\psi\circ\kappa$. Notice that with the notation of~(\ref{eq:Dixmier.Pcpsi}) we have 
$\overline{\tilde{c}_{1}}\tilde{c}_{1}P_{0}^{-(d+2)}\tilde{\psi}=P_{|c|^{2},\psi}$. Observe also that we have
\begin{equation*}
    (\tilde{c}_{1}P_{0}^{-\frac{d+2}{2}}\tilde{\psi})(\tilde{c}_{1}P_{0}^{-\frac{d+2}{2}}\tilde{\psi})^{*} 
    =\tilde{c}_{1}P_{0}^{-\frac{d+2}{2}}\tilde{\psi}^{2}P_{0}^{-\frac{d+2}{2}}\overline{\tilde{c}_{1}}
    =  \tilde{c}_{1}P_{0}^{-(d+2)}\psi \overline{\tilde{c}_{1}} \mod \Psi^{-\infty}(M,\cE).
\end{equation*}
As alluded to earlier the trace $\Trw$ vanishes on smoothing operators, so we get
\begin{multline}
  \Trw  [  (\tilde{c}_{1}P_{0}^{-\frac{d+2}{2}}\tilde{\psi})(\tilde{c}_{1}P_{0}^{-\frac{d+2}{2}}\tilde{\psi})^{*} ]= 
  \Trw [\tilde{c}_{1}P_{0}^{-(d+2)}\tilde{\psi} \overline{\tilde{c}_{1}}]\\ =\Trw 
    [\overline{\tilde{c}_{1}}\tilde{c}_{1}P_{0}^{-(d+2)}\tilde{\psi} ] = \Trw P_{|c|^{2},\psi}=L(|c|).
\end{multline}
Since $\Trw$ is a positive trace (cf.~Proposition~\ref{prop:NCG.properties-Dixmier-trace}) 
it follows that we have $L(|c|^{2})\geq 0$ for any $c \in C^{\infty}_{c}(\Rd)$, i.e., $L$ 
is a positive linear functional on $C_{c}^{\infty}(\Rd)$. Since any such functional uniquely extends to a Radon measure on $C_{0}^{\infty}(\Rd)$, 
this shows that $L$ defines a positive Radon measure.

Next, let $a \in \Rd$ and let $\phi(x)=x+a$ be the translation by $a$ on $\Rd$. Since $\phi'(x)=1$ we see that 
$\phi$ is a Heisenberg diffeomorphism, so for any $P \in \pvdoc^{*}(\Rd)$ the operator $\phi_{*}P$ is in $\pvdoc^{*}(\Rd)$ too. Set 
$\phi_{\kappa}=\kappa^{-1} \circ \phi \circ \kappa$. Then by~(\ref{eq:Log.functoriality-cP}) we have 
\begin{equation}
   \kappa_{*}[ \tr_{\cE} c_{\phi_{\kappa*}P_{c,\psi}}(x)]=\kappa_{*}\phi_{\kappa*}[\tr_{\cE}c_{P_{c,\psi}}(x)]= 
   \phi_{*}[c(x)dx]=c(\phi^{-1}(x))dx.
\end{equation}
Since shows that $a_{\phi_{\kappa*}P_{c,\psi}}(x)= c(\phi^{-1}(x))$, so from~(\ref{eq:NCG.Trw-tau-cP}) we  get 
\begin{equation}
    \Trw \phi_{\kappa*}P_{c,\psi}=L[c\circ \phi^{-1}].
    \label{eq:Dixmier.Lcphi}
\end{equation}

Let $K$ be a compact subset of $\Rd$. Then $\phi_{\kappa}$ gives rise to a continuous linear isomorpshism 
$\phi_{\kappa*}:L^{2}_{\kappa^{-1}(K)}(M,\cE)\rightarrow L^{2}_{\kappa^{-1}(K+a)}(M,\cE)$. By combining it with a continuous linear isomorphism
$L^{2}_{\kappa^{-1}(K)}(M,\cE)^{\perp}\rightarrow L^{2}_{\kappa^{-1}(K+a)}(M,\cE)^{\perp}$ we obtain a continuous linear isomorphism  
$S:L^{2}(M,\cE)\rightarrow L^{2}(M,\cE)$ which agrees with $\phi_{\kappa*}$ on $L^{2}_{\kappa^{-1}(K)}(M,\cE)$. In particular, we have 
$\phi_{\kappa*}P_{c,\psi} = S P_{c,\psi}S^{-1}$. Therefore, by using Proposition~\ref{prop:NCG.properties-Dixmier-trace} 
and~\ref{eq:Dixmier.Lcphi} we see that, for any $c \in C_{K}^{\infty}(\Rd)$, 
we have 
\begin{equation}
    L[c]=\Trw P_{c,\psi}=\Trw S P_{c,\psi}S^{-1}= \Trw \phi_{\kappa*}P_{c,\psi}=L[c\circ \phi^{-1}].
\end{equation}
This proves that $L$ is translation-invariant. Since any translation invariant Radon measure on $C^{\infty}_{c}(\Rd)$ is a constant multiple of the 
Lebesgue measure, it follows that there exists a constant $\Lambda_{U}\in \C$ such that,  for any $c \in C^{\infty}_{c}(\Rd)$, we have 
\begin{equation}
    L(c)=\Lambda_{U} \int c(x)dx. 
     \label{eq:Dixmier.L-Lebesgue}
\end{equation}

Now, combining~(\ref{eq:NCG.Trw-tau-cP}) and~(\ref{eq:Dixmier.L-Lebesgue}) shows that, for any $P\in \pvdoc^{-(d+2)}(U,\cE_{|U})$, we have
\begin{multline}
   \Trw P=  \Lambda_{U}\int_{\Rd} a_{P}(x)dx= \Lambda_{U}\int_{\Rd} \kappa_{*}[\tr_{\cE}c_{P}(x)] \\
    =\Lambda_{U}\int_{M}\tr_{\cE}c_{P}(x)=(2\pi)^{d+1}\Lambda_{U}\Res P.
     \label{eq:Dixmier-Trw-Res}
\end{multline}
This shows that, for any domain $U$ of a nice $H$-framed chart,  on $\pvdoc^{-(d+2)}(U,\cE_{|U})$ the Dixmier trace 
$\Trw $ is a constant multiple of the noncommutative residue. Therefore, if we let $M_{1}, \ldots, M_{N}$ be the connected components of $M$, 
then by arguing as in the proof of Theorem~\ref{thm:Traces.traces} we can prove that on each connected component $M_{j}$ there exists a constant 
$\Lambda_{j}\geq 0$ such that 
\begin{equation}
    \Trw P =\Lambda_{j} \Res P \qquad \forall P\in \pvdo^{-(d+2)}(M_{j},\cE_{|M_{j}}).
\end{equation}
In fact, if we take $P=P_{0|_{M_{j}}}^{-(d+2)}$ then from~(\ref{eq:NCG.Dixmier-trace-NCR.hypoelliptic}) we get $\Lambda_{j}=(d+2)^{-1}$. Thus, 
\begin{equation}
    \Trw P =\frac{1}{d+2}\Res P \qquad \forall P \in \pvdo^{-(d+2)}(M,\cE).
\end{equation}
This proves that any operator $P \in \pvdo^{-(d+2)}(M,\cE)$ is measurable and its Dixmier trace then is 
equal to $(d+2)^{-1}\Res P$. The theorem is thus proved. 
\end{proof}

As a consequence of Theorem~\ref{thm:NCG.Dixmier} we can extend the Dixmier trace to the whole algebra  $\pvdoz(M,\cE)$ by letting
\begin{equation}
        \bint P :=\frac{1}{d+2}\Res P \quad \text{for any $P\in \pvdoz(M,\cE)$}.
\end{equation}

In the language of the quantized calculus this means that we can integrate any \psivdo\ of integer order, even though it is not an infinitesimal operator 
of order~$\geq 1$. This property will be used in Section~\ref{sec:CR} to define lower dimensional volumes in 
pseudohermitian geometry.

\section{Noncommutative residue and contact geometry}
\label{sec:Contact}
In this section we make use of the results of~\cite{Po:MAMS1} to compute the noncommutative residues of
some geometric operators on contact manifolds. 

Throughout this section we let $(M^{2n+1},H)$ be a compact orientable contact manifold, i.e., $(M^{2n+1},H)$ is a Heisenberg manifold and there 
exists a contact 1-form $\theta$ on $M$ such that $H=\ker \theta$ (cf.~Section~\ref{sec:Heisenberg-calculus}). 

Since $M$ is orientable the hyperplane $H$ admits an almost complex structure $J\in C^{\infty}(M,\End H)$, $J^{2}=-1$, 
which is calibrated with respect to $\theta$, i.e., $d\theta(.,J.)$ is positive definite on $H$. 
We then can endow $M$ with the Riemannian metric,
\begin{equation}
    g_{\theta,J}=\theta^{2}+d\theta(.,J.).
     \label{eq:Contact.Riem-metric}
\end{equation}
The volume of $M$ with respect to $g_{\theta,J}$ depends only on $\theta$ and is equal to 
\begin{equation}
    \Vol_{\theta}M:=\frac{1}{n!}\int_{M}d\theta^{n}\wedge \theta.
\end{equation}

In addition, we let $X_{0}$ be the \emph{Reeb field} associated to $\theta$, that is, the unique vector field on $M$ such 
that $\iota_{X_{0}}\theta=1$ and $\iota_{X_{0}}d\theta=0$.

\subsection{Noncommutative residue and the horizontal sublaplacian (contact case)}
In the sequel we shall identify $H^{*}$ with the subbundle of $T^{*}M$ annihilating 
the orthogonal complement $H^{\perp}\subset TM$. This yields the orthogonal splitting, 
\begin{equation}
    \Lambda_{\C}T^{*}M=(\bigoplus_{0\leq k \leq 2n} \Lambda^{k}_{\C}H^{*})\oplus (\theta\wedge  \Lambda^{*}T_{\C}^{*}M).
\end{equation}
The horizontal differential $d_{b;k}:C^{\infty}(M,\Lambda^{k}_{\C}H^{*})\rightarrow C^{\infty}(M,\Lambda^{k+1}_{\C}H^{*})$ is 
\begin{equation}
    d_{b}=\pi_{b;k+1}\circ d,
\end{equation}
where $\pi_{b;k}\in C^{\infty}(M,\End  \Lambda_{\C}T^{*}M)$ denotes the orthogonal projection onto $\Lambda^{k}_{\C}H^{*}$. This is not the 
differential of a chain complex, for we have 
\begin{equation}
    d_{b}^{2}=-\cL_{X_{0}}\varepsilon(d\theta)=-\varepsilon (d\theta)\cL_{X_{0}},
     \label{eq:CR.square-db}
\end{equation}
where $\varepsilon (d\theta)$ denotes the exterior multiplication by $d\theta$. 

The horizontal sublaplacian $\Delta_{b;k}:C^{\infty}(M,\Lambda^{k}_{\C}H^{*})\rightarrow C^{\infty}(M,\Lambda^{k+1}_{\C}H^{*})$ is 
\begin{equation}
    \Delta_{b;k}=d_{b;k}^{*}d_{b;k}+d_{b;k-1}d_{b;k-1}^{*}.
    \label{eq:CR.horizontal-sublaplacian}
\end{equation}
Notice that the definition of $\Delta_{b}$ makes sense on any Heisenberg manifold equipped with a Riemannian metric. This operator was first introduced by 
Tanaka~\cite{Ta:DGSSPCM}, but versions of this operator acting on functions were independently defined by 
Greenleaf~\cite{Gr:FESPM} and Lee~\cite{Le:FMPHI}. 
Since the fact that $(M,H)$ is a contact manifold implies that the Levi form~(\ref{eq:Heisenberg.Levi-form}) 
is nondegenerate, from~\cite[Prop.~3.5.4]{Po:MAMS1} we get:

\begin{proposition}
    The principal symbol of $\Delta_{b;k}$ is invertible if and only if we have $k\neq n$. 
\end{proposition}
Next, for $\mu\in (-n,n)$ we let 
\begin{equation}
    \rho(\mu)= \frac{\pi^{-(n+1)}}{2^{n}n!} \int_{-\infty}^{\infty}e^{-\mu\xi_{0}}(\frac{\xi_{0}}{\sinh \xi_{0}})^{n}d\xi_{0}. 
\end{equation}
Notice that with the notation of~\cite[Eq.~(6.2.29)]{Po:MAMS1} we have $\rho(\mu)=(2n+2)\nu(\mu)$. 
For $q \neq n$ let  $\nu_{0}(\Delta_{b;k})$ be the coefficient $\nu_{0}(P)$ in the Weyl asymptotics~(\ref{eq:Zeta.Weyl-asymptotics1}) for 
$\Delta_{b;k}$, i.e., we have
 $\Res \Delta_{b;k}^{-(n+1)}=(2n+2)\nu_{0}(\Delta_{b;k})$. By~\cite[Prop.~6.3.3]{Po:MAMS1} we have 
$\nu_{0}(\Delta_{b;k})=\tilde{\gamma}_{nk}\Vol_{\theta}M$, where
$\tilde{\gamma}_{nk}:=\sum_{p+q=k}2^{n}\binom{n}{p} \binom{n}{q}\nu(p-q)$. Therefore, we get:

\begin{proposition}\label{prop:Contact.residue-Deltab}
    For $k \neq n$ we have 
    \begin{equation}
         \Res \Delta_{b;k}^{-(n+1)}= \gamma_{nk} \Vol_{\theta}M ,\quad 
 \gamma_{nk}=\sum_{p+q=k}2^{n}\binom{n}{p} \binom{n}{q}\rho(p-q). 
           \label{eq:Contact.residue-Deltab}
    \end{equation}
    In particular $\gamma_{nk}$ is a universal constant depending only on $n$ and $k$.
\end{proposition}

\subsection{Noncommutative residue and the contact Laplacian}
The contact complex of Rumin~\cite{Ru:FDVC} can be seen as an attempt to get a complex of horizontal forms by forcing the equalities 
$d_{b}^{2}=0$ and $(d_{b}^{*})^{2}=0$. Because of~(\ref{eq:CR.square-db}) there are two natural ways to modify $d_{b}$ to get a chain complex. 
The first one is to force the equality $d_{b}^{2}=0$ by restricting 
$d_{b}$ to the subbundle $\Lambda^{*}_{2}:=\ker \varepsilon(d\theta) \cap \Lambda^{*}_{\C}H^{*}$, since the latter
is closed under $d_{b}$ and is annihilated by $d_{b}^{2}$. 
Similarly, we get the equality $(d_{b}^{*})^{2}=0$ by restricting $d^{*}_{b}$ to the subbundle 
$\Lambda^{*}_{1}:=\ker \iota(d\theta)\cap \Lambda^{*}_{\C}H^{*}=(\im \varepsilon(d\theta))^{\perp}\cap \Lambda^{*}_{\C}H^{*}$, where 
$\iota(d\theta)$ denotes the interior product 
with $d\theta$. This amounts to replace $d_{b}$ by $\pi_{1}\circ d_{b}$, where $\pi_{1}$ is the orthogonal projection onto $\Lambda^{*}_{1}$.

In fact, since $d\theta$ is nondegenerate on $H$ the operator $\varepsilon(d\theta):\Lambda^{k}_{\C}H^{*}\rightarrow \Lambda^{k+2}_{\C}H^{*}$  is 
injective for $k\leq n-1$ and surjective for $k\geq n+1$. This implies that $\Lambda_{2}^{k}=0$ for $k\leq n$ and $\Lambda_{1}^{k}=0$ for $k\geq n+1$. 
Therefore, we only have two halves of complexes. 
As observed by Rumin~\cite{Ru:FDVC} we get a full complex by connecting the two halves by means of the differential operator, 
\begin{equation}
B_{R}:C^{\infty}(M,\Lambda_{\C}^{n}H^{*})\rightarrow C^{\infty}(M,\Lambda_{\C}^{n}H^{*}),   \qquad 
B_{R}=\cL_{X_{0}}+d_{b,n-1}\varepsilon(d\theta)^{-1}d_{b,n},
\end{equation}
where $\varepsilon(d\theta)^{-1}$ is the inverse of $\varepsilon(d\theta):\Lambda^{n-1}_{\C}H^{*}\rightarrow \Lambda^{n+1}_{\C}H^{*}$. Notice that 
$B_{R}$ is second order differential operator.  Thus, if we let $\Lambda^{k}=\Lambda_{1}^{k}$ for $k=0,\ldots,n-1$ and we let
$\Lambda^{k}=\Lambda_{1}^{k}$ for $k=n+1,\ldots,2n$, then we get the chain complex,
\begin{multline}
    C^{\infty}(M)\stackrel{d_{R;0}}{\rightarrow}C^{\infty}(M,\Lambda^{1})\stackrel{d_{R;1}}{\rightarrow}
    \ldots 
    C^{\infty}(M,\Lambda^{n-1})\stackrel{d_{R;n-1}}{\rightarrow}C^{\infty}(M,\Lambda^{n}_{1})\stackrel{B_{R}}{\rightarrow}\\  
    C^{\infty}(M,\Lambda^{n}_{2})  \stackrel{d_{R;n}}{\rightarrow}C^{\infty}(M,\Lambda^{n+1})
    \ldots \stackrel{d_{R;2n-1}}{\longrightarrow} C^{\infty}(M,\Lambda^{2n}),
     \label{eq:contact-complex}
\end{multline}
where $d_{R;k}:=\pi_{1}\circ d_{b;k}$ for $k=0,\ldots,n-1$ and $d_{R;k}:=d_{b;k}$ for $k=n,\ldots,2n-1$. 
This complex is called the \emph{contact complex}. 

The contact Laplacian is defined as follows. In degree $k\neq n$ it consists of the differential operator 
$\Delta_{R;k}:C^{\infty}(M,\Lambda^{k})\rightarrow C^{\infty}(M,\Lambda^{k})$ given by
\begin{equation}
    \Delta_{R;k}=\left\{
    \begin{array}{ll}
        (n-k)d_{R;k-1}d^{*}_{R;k}+(n-k+1) d^{*}_{R;k+1}d_{R;k}& \text{$k=0,\ldots,n-1$},\\
         (k-n-1)d_{R;k-1}d^{*}_{R;k}+(k-n) d^{*}_{R;k+1}d_{R;k}& \text{$k=n+1,\ldots,2n$}.
         \label{eq:contact-Laplacian1}
    \end{array}\right.
\end{equation}
In degree $k=n$ it consists of the differential operators $\Delta_{R;nj}:C^{\infty}(M,\Lambda_{j}^{n})\rightarrow C^{\infty}(M,\Lambda^{n}_{j})$, $j=1,2$, 
defined by the formulas,  
\begin{equation}
    \Delta_{R;n1}= (d_{R;n-1}d^{*}_{R;n})^{2}+B_{R}^{*}B_{R}, \quad   \Delta_{R;n2}=B_{R}B_{R}^{*}+  (d^{*}_{R;n+1}d_{R;n}).
         \label{eq:contact-Laplacian2}
\end{equation}

Observe that $\Delta_{R;k}$, $k\neq n$, is a differential operator of order $2$, whereas $\Delta_{R;n1}$ and $\Delta_{R;n2}$ are differential operators of 
order $4$. Moreover, Rumin~\cite{Ru:FDVC} proved that in every degree the contact Laplacian is maximal hypoelliptic in the sense of~\cite{HN:HMOPCV}. 
In fact, in every degree the contact Laplacian has an invertible principal symbol, hence admits a parametrix in the Heisenberg calculus 
(see~\cite{JK:OKTGSU},~\cite[Sect.~3.5]{Po:MAMS1}).

For $k\neq n$ (resp.~$j=1,2$) we let $\nu_{0}(\Delta_{R;k})$ (resp.~$\nu_{0}(\Delta_{R;nj})$) be the coefficient 
$\nu_{0}(P)$ in the Weyl asymptotics~(\ref{eq:Zeta.Weyl-asymptotics1}) 
for $\Delta_{R;k}$ (resp.~$\Delta_{R;nj}$). By Proposition~\ref{prop:Zeta.Weyl-asymptotics} we have 
$\Res \Delta_{R;k}^{-(n+1)}=(2n+2)\nu_{0}(\Delta_{R;k})$ and $ \Res \Delta_{R;nj}^{-\frac{n+1}{2}}=(2n+2)\nu_{0}(\Delta_{R;nj})$. 
Moreover, by~\cite[Thm.~6.3.4]{Po:MAMS1} there exist  universal positive constants $\nu_{nk}$ and $\nu_{n,j}$ depending 
 only on $n$, $k$ and $j$ such that $\nu_{0}(\Delta_{R;k})=\nu_{nk}\Vol_{\theta}M$ and $\nu_{0}(\Delta_{R;nj})= \nu_{n,j}\Vol_{\theta}M$. 
 Therefore, we obtain: 

\begin{proposition}\label{prop:Contact.residue-DeltaR}
    1)  For $k \neq n$ there exists a  universal constant $\rho_{nk}>0$ depending only on $n$ and $k$ such that
     \begin{equation}
           \Res \Delta_{R;k}^{-(n+1)}=\rho_{nk} \Vol_{\theta}M. 
      \end{equation}
    
    2) For $j=1,2$ there exists a  universal constant $\rho_{n,j}>0$ depending only on $n$ and $j$  such that
    \begin{equation}
      \Res \Delta_{R;nj}^{-\frac{n+1}{2}}= \rho_{n,j}  \Vol_{\theta}M.  
    \end{equation}
\end{proposition}
\begin{remark}
    We have $\rho_{nk}=(2n+2)\nu_{nk}$ and $\rho_{n,j}=(2n+2)\nu_{n,j}$, so it follows from the proof of~\cite[Thm.~6.3.4]{Po:MAMS1} that we can 
    explicitly relate the universal constants $\rho_{nk}$ and $\rho_{n,j}$ to the fundamental solutions of the heat operators 
    $\Delta_{R;k}+\partial_{t}$ and $\Delta_{R;nj}+\partial_{t}$ associated to the contact Laplacian on the Heisenberg group $\bH^{2n+1}$ 
    (cf.~\cite[Eq.~(6.3.18)]{Po:MAMS1}). For 
    instance, if $K_{0;k}(x,t)$ denotes the fundamental solution of $\Delta_{R;0}+\partial_{t}$ on $\bH^{2n+1}$ then we have $\rho_{n,0}= 
    \frac{2^{n}}{n!}K_{0;0}(0,1)$. 
\end{remark}

\section{Applications in CR geometry} 
\label{sec:CR}
In this section we present some applications in CR geometry of the noncommutative residue for the Heisenberg calculus. After recalling the geometric 
set-up, we shall compute the noncommutative residues of some powers of the horizontal sublaplacian and of the Kohn 
Laplacian on CR manifolds endowed with a pseudohermitian structure. 
After this we will make use of the framework of noncommutative geometry to define lower dimensional volumes in pseudohermitian geometry. For 
instance, we will give sense to the area of any 3-dimensional pseudohermitian manifold as a constant multiple the integral of the Tanaka-Webster scalar curvature. 
As a by-product this will allow us to get a spectral interpretation of the Einstein-Hilbert action in pseudohermitian geometry.

\subsection{The geometric set-up}
Let $(M^{2n+1},H)$ be a compact orientable CR manifold. Thus $(M^{2n+1},H)$ is a Heisenberg manifold such that $H$ 
admits a complex structure $J\in C^{\infty}(M,\End H)$, $J^{2}=-1$, in such way that $T_{1,0}:=\ker (J+i)\subset T_{\C}M$ is a complex rank $n$ subbundle 
which is integrable in Fr\"obenius' sense (cf.~Section~\ref{sec:Heisenberg-calculus}). In addition, we set $T_{0,1}=\overline{T_{1,0}}=\ker(J-i)$. 

Since $M$ is orientable and $H$ is orientable by means of its complex structure, there exists a global non-vanishing real 1-form $\theta$ such 
that $H=\ker \theta$. Associated to $\theta$ is its Levi form, i.e., the Hermitian form on $T_{1,0}$ such that 
\begin{equation}
    L_{\theta}(Z,W)=-id\theta(Z,\overline{W}) \qquad \forall Z,W \in T_{1,0}.
\end{equation}

\begin{definition}
 We say that  $M$ is strictly pseudoconvex (resp.~$\kappa$-strictly pseudoconvex) when we can choose $\theta$ so that $L_{\theta}$ is 
positive definite (resp.~has signature $(n-\kappa,\kappa,0)$) at every point. 
\end{definition}

If $(M, H)$ is $\kappa$-strictly pseudoconvex then $\theta$ is a contact form on $M$. Then in the terminology of \cite{We:PHSRH} the datum of the contact form 
$\theta$ annihilating $H$ defines a \emph{pseudohermitian structure} on $M$. 

From now we assume that $M$ is $\kappa$-strictly pseudoconvex, and we let $\theta$ be a pseudohermitian contact form 
such that $L_{\theta}$ has signature $(n-\kappa,\kappa,0)$ everywhere. We let $X_{0}$ be the Reeb vector field associated to $\theta$, so that 
$\iota_{X_{0}}\theta=1$ and $\iota_{X_{0}}d\theta=0$ (cf.~Section~\ref{sec:Contact}), and we let $\cN\subset T_{\C}M$ be the complex line bundle spanned by 
$X_{0}$. 

We endow $M$ with a \emph{Levi metric} as follows. First, we always can construct a splitting $T_{1,0}=T_{1,0}^{+}\oplus T_{1,0}^{+}$ with 
subbundles $T_{1,0}^{+}$ and $T_{1,0}^{-}$  which are orthogonal with respect to $L_{\theta}$ and 
such that $L_{\theta}$ is positive definite on $T_{1,0}^{+}$ and negative definite on $T_{1,0}^{-}$ (see, e.g., \cite{FS:EDdbarbCAHG}, \cite{Po:MAMS1}). 
Set $T_{0,1}^{\pm}=\overline{T_{1,0}^{\pm}}$. Then we have the splittings, 
\begin{equation}
    T_{\C}M=\cN\oplus T_{1,0}\oplus T_{0,1}=\cN \oplus T_{1,0}^{+}\oplus T_{1,0}^{-}\oplus T_{0,1}^{+}\oplus T_{0,1}^{-}.
    \label{eq:CR.splitting-TcM}
\end{equation}
Associated to these splittings is the unique Hermitian metric $h$ on $T_{\C}M$ such that:\smallskip

- The splittings~(\ref{eq:CR.splitting-TcM}) are orthogonal with respect to $h$;\smallskip

 - $h$ commutes with complex conjugation;\smallskip

- We have $h(X_{0},X_{0})=1$ and $h$ agrees with $\pm L_{\theta}$ on $T_{1,0}^{\pm}$.\smallskip

\noindent In particular, the matrix of $L_{\theta}$ with respect to $h$ is $\op{diag}(1,\ldots,1,-1,\ldots,-1)$, where $1$ has multiplicity 
$n-\kappa$ and $-1$ multiplicity $-1$. 

Notice that when $M$ is strictly pseudoconvex $h$ is uniquely determined by $\theta$, since in this case $T_{1,0}^{+}=T_{1,0}$ and one can check that 
we have $h=\theta^{2}+d\theta(.,J.)$, that is, $h$ agrees on $TM$ with the Riemannian metric $g_{\theta,J}$ in~(\ref{eq:Contact.Riem-metric}). 
In general, we can check that the volume form of $M$ with respect to $h$ depends only on $\theta$ and is equal to 
\begin{equation}
    v_{\theta}(x):=\frac{(-1)^{\kappa}}{n!}d\theta^{n}\wedge\theta. 
\end{equation}
In particular, the volume of $M$ with respect to $h$ is 
\begin{equation}
    \Vol_{\theta}M:=\frac{(-1)^{\kappa}}{n!}\int_{M} d\theta^{n}\wedge\theta.
\end{equation}

Finally, as proved by Tanaka~\cite{Ta:DGSSPCM} and Webster~\cite{We:PHSRH} the datum of the pseudohermitian contact 
form $\theta$ defines a natural connection, the \emph{Tanaka-Webster connection}, which preserves the pseudohermitian structure of $M$, i.e., it preserves both 
$\theta$ and $J$. It can be defined as follows. 

Let $\{Z_{j}\}$ be a local frame of $T_{1,0}$. Then 
$\{X_{0},Z_{j},Z_{\overline{j}}\}$ forms a frame of $T_{\C}M$ with dual coframe $\{\theta,\theta^{j}, \theta^{\overline{j}}\}$, 
with respect to which we can write $d\theta =ih_{j\overline{k}}\theta^{j}\wedge \theta^{\overline{k}}$. 
Using the matrix $(h_{j\overline{k}})$ and its inverse $(h^{j\overline{k}})$ to lower and raise indices, the connection 1-form 
$\omega=(\omega_{j}^{~k})$ and the 
 torsion form $\tau_{k}=A_{jk}\theta^{j}$ of the Tanaka-Webster connection  are uniquely determined by the relations, 
\begin{equation}
     d\theta^{k}=\theta^{j}\wedge \omega_{j}^{~k}+\theta \wedge \tau^{k}, \qquad 
     \omega_{j\bar{k}} + \omega_{\bar{k}j} 
         =dh_{j\bar{k}}, \qquad  A_{jk}=A_{k j}. 
\end{equation}

The curvature tensor $\Pi_{j}^{~k}:=d\omega_{j}^{~k}-\omega_{j}^{~l}\wedge  \omega_{l}^{~k}$ satisfies the structure equations,
\begin{equation}
  \Pi_{j}^{~k}=R_{j\bar{k} l\bar{m}} \theta^{l}\wedge \theta^{\bar{m}} + 
    W_{j\bar{k}l}\theta^{l}\wedge \theta -  W_{\bar{k}j\bar{l}}\theta^{\bar{l}}\wedge \theta 
    +i\theta_{j}\wedge \tau_{\bar{k}}-i\tau_{j}\wedge \theta_{\bar{k}}.
     \label{eq:CR.TW-curvature}
\end{equation}
The \emph{Ricci tensor} of the Tanaka-Webster connection is $ \rho_{j \bar{k}}:=R_{l~j \bar{k}}^{~l}$,
and its \emph{scalar curvature} is $R_{n}: =\rho_{j}^{~j}$.

\subsection{Noncommutative residue and the Kohn Laplacian}
\label{subsec:CR.NCR-pseudohermitian}
The $\dbarb$-complex  of 
Kohn-Rossi~(\cite{KR:EHFBCM},~\cite{Ko:BCM}) is defined as follows. 

Let $\Lambda^{1,0}$ (resp.~$\Lambda^{0,1}$) be the annihilator of $T_{0,1}\oplus \cN$ (resp.~$T_{0,1}\oplus \cN$) in $T^{*}_{\C}M$. 
For $p,q=0,\ldots,n$ let $\Lambda^{p,q}:=(\Lambda^{1,0})^{p}\wedge (\Lambda^{0,1})^{q}$ be the bundle of $(p,q)$-covectors on $M$, so that 
we have the orthogonal decomposition,
\begin{equation}
    \Lambda^{*}T_{\C}^{*}M=(\bigoplus_{p,q=0}^{n}\Lambda^{p,q})\oplus (\theta\wedge  \Lambda^{*}T_{\C}^{*}M).
     \label{eq:CR-Lambda-pq-decomposition}
\end{equation}
Moreover, thanks to the integrability of $T_{1,0}$, given any local section $\eta$ of $\Lambda^{p,q}$, its differential $d\eta$ can be uniquely 
decomposed as
\begin{equation}
    d\eta =\dbarbpq \eta + \partial_{b;p,q}\eta + \theta \wedge \cL_{X_{0}}\eta,
     \label{eq:CR.dbarb}
\end{equation}
where $\dbarbpq \eta $ (resp.~$\partial_{b;p,q}\eta$) is a section of $\Lambda^{p,q+1}$ (resp.~$\Lambda^{p+1,q}$). 

The integrability of $T_{1,0}$ further implies that $\overline{\partial}_{b}^{2}=0$ on $(0,q)$-forms, so that we get the cochain 
complex $\overline{\partial}_{b;0,*}:C^{\infty}(M,\Lambda^{0,*})\rightarrow C^{\infty}(M,\Lambda^{0,*+1})$. On $(p,q)$-forms with $p\geq 1$ the operator 
$\dbarb^{2}$ is a tensor which vanishes when the complex structure $J$ is invariant under the Reeb flow (i.e., 
when we have $[X_{0},JX]=J[X_{0},X]$ for any local section $X$ of $H$). 

 Let $\dbarbpq^{*}$ be 
 the formal adjoint of $\dbarbpq$ with respect to the Levi metric of $M$.   Then the \emph{Kohn Laplacian}  
$\Boxbpq :C^{\infty}(M,\Lambda^{p,q})\rightarrow C^{\infty}(M,\Lambda^{p,q})$  is defined to be
\begin{equation}
   \Boxbpq=\dbarbpq^{*}\dbarbpq +  \overline{\partial}_{b;p,q-1}\overline{\partial}_{b;p,q-1}^{*}.
\end{equation}
This a differential operator which has order 2 in the Heisenberg calculus sense. Furthermore, we have:

\begin{proposition}[\cite{BG:CHM}]
    The principal symbol of $\Boxbpq$ is invertible if and only if we have $q \neq \kappa$ and $q\neq n-\kappa$.  
\end{proposition}

Next, for $q \not\in\{\kappa,n-\kappa\}$ let $\nu_{0}(\Boxbpq)$ be the coefficient $\nu_{0}(P)$ in the Weyl asymptotics~(\ref{eq:Zeta.Weyl-asymptotics1})
 for $\Boxbpq$. By~\cite[Thm.~6.2.4]{Po:MAMS1} we have
 $  \nu_{0}(\Boxbpq)=\tilde{\alpha}_{n\kappa pq}\Vol_{\theta}M$, where $\tilde{\alpha}_{n\kappa pq}$ is equal to  
\begin{equation}
   \sum_{\max(0,q-\kappa)\leq  k\leq \min(q,n-\kappa)} \frac{1}{2} \binom{n}{p} \binom{n-\kappa}{k}\binom{\kappa}{q-k} 
        \nu(n-2(\kappa-q+2k)).
        \label{eq:CR.talphankpq}
 \end{equation}
Therefore, by arguing as in the proof of Proposition~\ref{prop:Contact.residue-Deltab} we get: 

\begin{proposition}\label{prop:CR.residue-Boxb1}
 For $q \neq \kappa$ and $q\neq n-\kappa$  we have
 \begin{equation}
     \Res \Boxbpq^{-(n+1)}=\alpha_{n\kappa pq}\Vol_{\theta}M, 
      \label{eq:CR.residue-Boxb1}
 \end{equation}
where $\alpha_{n\kappa pq}$ is equal to  
\begin{equation}
   \sum_{\max(0,q-\kappa)\leq  k\leq \min(q,n-\kappa)} \frac{1}{2} \binom{n}{p} \binom{n-\kappa}{k}\binom{\kappa}{q-k} 
        \rho(n-2(\kappa-q+2k)).
          \label{eq:CR.alphapq}  
\end{equation}
In particular $\alpha_{n\kappa pq}$ is a universal constant depending only on $n$, $\kappa$, $p$ and $q$.
\end{proposition}

\begin{remark}\label{rem:CR.residue-Boxb1-local}
    Let $a_{0}(\Boxbpq)(x)$ be the leading coefficient in the heat kernel asymptotics~(\ref{eq:Zeta.heat-kernel-asymptotics}) for $\Boxbpq$. 
    By~(\ref{eq:Zeta.tPs-heat1}) we have $ \nu_{0}(\Boxbpq)= \frac{1}{ (n+1)!} \int_{M}\tr_{\Lambda^{p,q}}a_{0}(\Boxbpq)(x)$. 
    Moreover, a careful look at the proof 
    of~\cite[Thm.~6.2.4]{Po:MAMS1} shows that we have 
    \begin{equation}
        \tr_{\Lambda^{p,q}}a_{0}(\Boxbpq)(x)=(n+1)!\tilde{\alpha}_{n\kappa pq}v_{\theta}(x).
    \end{equation}
     Since by~(\ref{eq:Zeta.tPs-heat1}) 
    we have $2c_{\Boxbpq^{-(n+1)}}(x)=(n!)^{-1}a_{0}(\Boxbpq)(x)$, it follows that the equality~(\ref{eq:CR.residue-Boxb1}) 
    ultimately holds at the level of densities, that is, we have 
    \begin{equation}
        c_{\Boxbpq^{-(n+1)}}(x)=\alpha_{n\kappa pq}v_{\theta}(x).
    \end{equation}
 \end{remark}

Finally, when $M$ is strictly  pseudoconvex, i.e., when $\kappa=0$, we have:

\begin{proposition}\label{prop:CR.residue-Boxb2}
 Assume $M$ strictly pseudoconvex. Then for $q =1,\ldots, n-1$ there exists a universal constant $\alpha_{npq}'$ 
 depending only on $n$, $p$ and $q$ such that 
     \begin{equation}
        \Res \Boxbpq^{-n}=\alpha_{npq}'\int_{M}R_{n}d\theta^{n}\wedge \theta,
    \end{equation}
    where $R_{n}$ denotes the Tanaka-Webster scalar curvature of $M$. 
\end{proposition}
\begin{proof}
For $q=1,\ldots,n-1$ let $a_{2}(\Boxbpq)(x)$ be the coefficient of 
$t^{-n}$ in the heat kernel asymptotics~(\ref{eq:Zeta.heat-kernel-asymptotics}) for $\Boxbpq$. 
By~(\ref{eq:Zeta.tPs-heat1}) we have $2c_{\Boxbpq^{-n}}(x)=\Gamma(n)^{-1}a_{2}(\Boxbpq)(x)$. Moreover, 
by~\cite[Thm.~8.31]{BGS:HECRM} there exists a universal constant $\alpha_{npq}'$ depending only on $n$, $p$ and $q$ such that 
$\tr_{\Lambda^{p,q}}a_{2}(\Boxbpq)(x)=\alpha_{npq}'R_{n}d\theta^{n}\wedge \theta$. Thus, 
\begin{equation}
    \Res \Boxbpq^{-n}= \int_{M}\tr_{\Lambda^{p,q}}c_{\Boxbpq^{-n}}(x)= \alpha_{npq}'\int_{M}R_{n}d\theta^{n}\wedge \theta,
\end{equation}
where $\alpha_{npq}'$ is a universal constant depending only on $n$, $p$ and $q$.
\end{proof}

\subsection{Noncommutative residue and the horizontal sublaplacian (CR case)}
Let us identify $H^{*}$ with the subbundle of $T^{*}M$ annihilating 
the orthogonal supplement $H^{\perp}$, and let $\Delta_{b}:C^{\infty}(M,\Lambda^{*}_{\C}H^{*})\rightarrow C^{\infty}(M,\Lambda^{*}_{\C}H^{*})$ be the 
horizontal sublaplacian on $M$ as defined in~(\ref{eq:CR.horizontal-sublaplacian}). 

Notice that with the notation of~(\ref{eq:CR.dbarb}) we have  $d_{b}=\dbarb +\partial_{b}$. Moreover, we can check that 
$\dbarb \partial_{b}^{*}+\partial_{b}^{*} \dbarb= \dbarb^{*} \partial_{b}+\partial_{b} 
 \dbarb^{*}=0$. Therefore, we have
  \begin{equation}
      \Delta_{b}=\Box_{b}+ \overline{\Box}_{b}, \qquad \overline{\Box}_{b}:=\partial_{b}^{*}\partial_{b}+\partial_{b}\partial_{b}^{*}.
  \end{equation}
In particular, this shows that the horizontal sublaplacian $\Delta_{b}$ preserves the  bidegree, so 
 it induces a differential operator  $\Delta_{b;p,q}:C^{\infty}(M,\Lambda^{p,q})\rightarrow C^{\infty}(M,\Lambda^{p,q})$. Then the following holds. 

\begin{proposition}[{\cite[Prop.~3.5.6]{Po:MAMS1}}] 
    The principal symbol of $\Delta_{b;p,q}$ is invertible if and only if we have $(p,q)\neq (\kappa,n-\kappa)$ 
  and $(p,q)\neq (n-\kappa,\kappa)$.
\end{proposition}

Bearing this in mind we have:
\begin{proposition}\label{prop:CR.residue-Deltab1}
For $(p,q)\neq (\kappa,n-\kappa)$ and $(p,q)\neq (n-\kappa,\kappa)$ we have
\begin{equation}
             \Res \Delta_{b;p,q}^{-(n+1)}= \beta_{n\kappa pq}\Vol_{\theta}M,
     \label{eq:CR.residue-Deltab1a}
\end{equation}
 where $\beta_{n\kappa pq}$ is equal to
\begin{equation}
        \!  \!  \!  \!   \sum_{\substack{\max(0,q-\kappa)\leq  k\leq \min(q,n-\kappa)\\ \max(0,p-\kappa)\leq l\leq \min(p,n-\kappa)}}  \!  \!  \!  \!
        2^{n}\binom{n-\kappa}{l}\binom{\kappa}{p-l} \binom{n-\kappa}{k}\binom{\kappa}{q-k} 
        \rho(2(q-p)+4(l-k)).
         \label{eq:CR.residue-Deltab1b}
\end{equation}
In particular $\beta_{n\kappa pq}$  is a universal constant depending only on $n$, $\kappa$, $p$ and $q$.
\end{proposition}
\begin{proof}
 Let $\nu_{0}(\Delta_{b;p,q})$ be the coefficient $\nu_{0}(P)$ in the Weyl asymptotics~(\ref{eq:Zeta.Weyl-asymptotics1}) 
for $\Delta_{b;p,q}$. By~\cite[Thm.~6.2.5]{Po:MAMS1} we have $\nu_{0}(\Delta_{b;p,q})=\frac{1}{2n+2}\beta_{n\kappa pq}\Vol_{\theta}M$, 
where $\beta_{n\kappa pq}$ is given by~(\ref{eq:CR.residue-Deltab1b}). We then can show that $ \Res \Delta_{b;p,q}^{-(n+1)}= \beta_{n\kappa pq}\Vol_{\theta}M$
by arguing as in the proof of Proposition~\ref{prop:Contact.residue-Deltab}. 
\end{proof}

\begin{remark}\label{rem:CR.residue-Deltab1-local}
In the same way as~(\ref{eq:CR.residue-Boxb1}) (cf.~Remark~\ref{rem:CR.residue-Boxb1-local}) the equality~(\ref{eq:CR.residue-Deltab1a}) holds at the level 
of densities, i.e., we have $c_{\Delta_{b;p,q}^{-(n+1)}}(x)=\beta_{n\kappa pq}v_{\theta}(x)$. 
\end{remark}

\begin{proposition}\label{prop:CR.residue-Deltab2}
 Assume that $M$ is strictly pseudoconvex. For $(p,q)\neq (0,n)$ and $(p,q)\neq (n,0)$ there exists a universal constant $\beta_{npq}'$ 
 depending only $n$, $p$ and $q$ such that 
\begin{equation}
        \Res \Delta_{b;p,q}^{-n}=\beta_{npq}' \int_{M}R_{n}d\theta^{n}\wedge \theta. 
         \label{eq:CR.residue-Deltab2}
\end{equation}
\end{proposition}
\begin{proof}
   The same analysis as that of~\cite[Sect.~8]{BGS:HECRM} for the coefficients in the heat kernel asymptotics~(\ref{eq:Zeta.heat-kernel-asymptotics}) 
   for the Kohn Laplacian can be carried out for the coefficients of the heat kernel asymptotics for $\Delta_{b;p,q}$  (see~\cite{St:SICRM}). 
    In particular, if we let $a_{2}(\Delta_{b;p,q})(x)$ be the coefficient of $t^{-n}$ in the heat kernel 
    asymptotics for $\Delta_{b;p,q}$, then there exists a universal constant $\tilde{\beta}_{npq}$ 
 depending only on $n$, $p$ and $q$ such that $\tr_{\Lambda^{p,q}}a_{2}(\Delta_{b;p,q})(x)=\tilde{\beta}_{npq}R_{n}d\theta^{n}\wedge \theta$. 
 Arguing as in the proof of Proposition~\ref{prop:CR.residue-Boxb2} then shows that  
 $\Res \Delta_{b;p,q}^{-n}=\beta_{npq}' \int_{M}R_{n}d\theta^{n}\wedge \theta$, 
 where$\beta_{npq}'$ is a universal constant depending only $n$, $p$ and $q$.
\end{proof}

\subsection{Lower dimensional volumes in pseudohermitian geometry}
\label{subsec.CR.area}
Following an idea of Connes~\cite{Co:GCMFNCG}
we can make use of the noncommutative residue for classical \psidos\ to define lower dimensional dimensional volumes in 
Riemannian geometry, e.g., we can give sense to the area and the length of a Riemannian manifold even when the dimension is not 1 or 2
(see~\cite{Po:LMP07}). We shall now make use of the noncommutative residue for the Heisenberg calculus to define lower dimensional
volumes in pseudohermitian geometry. 

In this subsection we assume that $M$ is strictly pseudoconvex. In particular, the Levi metric $h$ is uniquely determined 
by $\theta$. In addition, we let $\Delta_{b;0}$ be the horizontal sublaplacian acting on functions.Then, as explained in Remark~\ref{rem:CR.residue-Deltab1-local}, 
we have
$c_{\Delta_{b;0}^{-(n+1)}}(x)=\beta_{n}v_{\theta}(x)$, where $\beta_{n}=\beta_{n000}=2^{n}\rho(0)$. In particular, for any $f \in C^{\infty}(M)$ 
we get  
$c_{f\Delta_{b;0}^{-(n+1)}}(x)=\beta_{n}f(x)v_{\theta}(x)$. Combining this with Theorem~\ref{thm:NCG.Dixmier}
then gives
\begin{equation}
    \bint f\Delta_{b;0}^{-(n+1)}=\frac{1}{2n+2}\int_{M}c_{f\Delta_{b;0}^{-(n+1)}}(x)=\frac{\beta_{n}}{2n+2}\int_{M}f(x)v_{\theta}(x).
\end{equation}
Thus the operator $\frac{2n+2}{\beta_{n}}\Delta_{b;0}^{-(n+1)}$ allows us to recapture the volume form 
$v_{\theta}(x)$. 

Since $-(2n+2)$ is the critical order for a \psivdo\ to be trace-class and $M$ has Hausdorff dimension $2n+2$ with respect to the 
Carnot-Carath\'eodory metric defined by the Levi metric on $H$, it stands for reason to define the \emph{length element} of 
$(M,\theta)$ as the positive selfadjoint operator  
$ds$ such that $(ds)^{2n+2}=\frac{2n+2}{\beta_{n}}\Delta_{b;0}^{-(n+1)}$, that is, 
\begin{equation}
    ds:= c_{n}\Delta_{b;0}^{-1/2}, \qquad c _{n}=(\frac{2n+2}{\beta_{n}})^{\frac{1}{2n+2}}.
\end{equation}

\begin{definition}
For $k=1,2,\ldots,2n+2$ the $k$-dimensional volume of $(M,\theta)$~is
\begin{equation}
    \op{Vol}_{\theta}^{(k)}M:=\bint ds^{k}.
\end{equation}
In particular, for $k=2$ the area of $(M,\theta)$ is $ \op{Area}_{\theta}M:=\bint ds^{2}$. 
\end{definition}

We have $\bint ds^{k}=\frac{(c_{n})^{k}}{2n+2}\int_{M}c_{\Delta_{b;0}^{-\frac{k}{2}}}(x)$ and thanks to~(\ref{eq:Zeta.tPs-heat1}) 
we know that $2 
c_{\Delta_{b;0}^{-\frac{k}{2}}}(x)$ agrees with $\Gamma(\frac{k}{2})^{-1}a_{2n+2-k}(\Delta_{b;0})(x)$, where $a_{j}(\Delta_{b;0})(x)$ 
denotes the coefficient of 
$t^{\frac{2n+2-j}{2}}$ in the heat kernel asymptotics~(\ref{eq:Zeta.heat-kernel-asymptotics}) for $\Delta_{b;0}$. Thus, 
\begin{equation}
     \op{Vol}_{\theta}^{(k)}M=\frac{(c_{n})^{k}}{4(n+1)}\Gamma(\frac{k}{2})^{-1}\int_{M}a_{2n+2-k}(\Delta_{b})(x).
\end{equation}
Since $\Delta_{b;0}$ is a differential operator we have $a_{2j-1}(\Delta_{b;0})(x)=0$ for any $j\in \N$, so $ \op{Vol}_{\theta}^{(k)}M$ 
vanishes when $k$ is odd. Furthermore, as alluded to in the proof of Proposition~\ref{prop:CR.residue-Deltab2} 
the analysis in~\cite[Sect.~8]{BGS:HECRM} of the coefficients of the heat kernel asymptotics for the Kohn Laplacian 
applies \emph{verbatim} to the heat kernel asymptotics for the horizontal sublaplacian. Thus, we can write
\begin{equation}
    a_{2j}(\Delta_{b;0})(x)=\gamma_{nj}(x)d\theta^{n}\wedge \theta(x),
    \label{eq:CR-volumes-gamma-nj}
\end{equation}
where $\gamma_{nj}(x)$ is a universal linear combination, depending only on $n$ and $j$, in complete contractions of  
covariant derivatives of the curvature and torsion tensors of the Tanaka-Webster connection (i.e.~$\gamma_{nj}(x)$ is a local pseudohermitian 
invariant). 
In particular, we have $\gamma_{n0}(x)=\gamma_{n0}$ and $\gamma_{n1}=\gamma_{n1}'R_{n}(x)$, 
where $\gamma_{n0}$ and 
$\gamma_{n1}$ are universal constants and $R_{n}(x)$ is the Tanaka-Webster scalar curvature (in fact the constants $\gamma_{n0}$ and 
$\gamma_{n1}'$ can be explicitly related to the constants $\beta_{n000}$ and $\beta_{n00}'$). Therefore, we obtain: 

\begin{proposition}\label{prop:CR.lower-dim.-volumes}
    1) $\Vol^{(k)}_{\theta}M$ vanishes when $k$ is odd.\smallskip

2) When $k$ is even we have
\begin{equation}
     \op{Vol}_{\theta}^{(k)}M=\frac{(c_{n})^{k}}{4(n+1)}\Gamma(\frac{k}{2})^{-1}\int_{M}\tilde{\gamma}_{nk}(x)d\theta^{n}\wedge \theta(x).
     \label{eq:CR.volumes-even}
\end{equation}
where $\tilde{\gamma}_{nk}(x):=\gamma_{nn+1-\frac{k}{2}}(x)$ is a universal linear combination, depending only on $n$ and $k$, of complete contractions of 
weight $n+1-\frac{k}{2}$ of covariant 
derivatives of the curvature and torsion tensors of the Tanaka-Webster connection. 
\end{proposition}

In particular, thanks to~(\ref{eq:CR.volumes-even}) we have a purely differential-geometric formulation of the $k$-dimensional volume 
$ \op{Vol}_{\theta}^{(k)}M$. Moreover, for $k=2n+2$ we get: 
\begin{equation}
    \op{Vol}_{\theta}^{(2n+2)}M=\frac{(c_{n})^{2n+2}}{4(n+1)}\frac{\gamma_{n0}}{n!}\int_{M}d\theta^{n}\wedge \theta. 
\end{equation}
Since $\op{Vol}_{\theta}^{(2n+2)}M=\op{Vol}_{\theta}M=\frac{1}{n!}\int_{M}d\theta^{n}\wedge \theta$ we see that 
$(c_{n})^{2n+2}=\frac{4(n+1)}{\gamma_{n0}}$, where $\gamma_{n0}$ is above.

On the other hand,  when $n=1$ (i.e.~$\dim M=3$) and $k=2$ we get
\begin{equation}
     \op{Area}_{\theta}M=\gamma_{1}''\int_{M}R_{1}d\theta\wedge \theta, \qquad \gamma_{1}'':=\frac{(c_{1})^{2}}{8}\gamma_{11}' 
     =\frac{\gamma_{11}'}{\sqrt{8\gamma_{10}}}, 
     \label{eq:CR.area-universal}
\end{equation}
where $\gamma_{11}'$ is above. 
To compute $\gamma_{1}''$ it is enough to compute $\gamma_{10}$ and $\gamma_{11}'$ in the special case of the unit sphere $S^{3}\subset \C^{2}$ 
equipped with its standard pseudohermitian structure, i.e., for $S^{3}$ equipped with the CR structure induced by the complex structure of $\C^{2}$ 
and with the 
pseudohermitian contact form $\theta:= \frac{i}2 (z_{1} d\bar{z}_{1}+ z_{2} d\bar{z}_{2})$.  

First, the volume $\Vol_{\theta}S^{3}$ is equal to 
    \begin{equation}
     \int_{S^{3}}d\theta\wedge \theta = \frac{-1}{4}\int_{S^{3}}(z_{2}dz_{1}\wedge d\bar{z_{1}} \wedge d\bar{z_{2}} + 
     z_{1} dz_{1}\wedge dz_{2}\wedge d\bar{z_{2}})=\pi^{2}.
        \label{eq:CR.volumeS3}
    \end{equation}
Moreover, by~\cite{We:PHSRH} the Tanaka-Webster scalar here is $R_{1}=4$, so we get  
\begin{equation}
    \int_{S^{3}}R_{1}d\theta\wedge \theta =4\Vol_{\theta}S^{3}=4\pi^{2}. 
\end{equation}    

Next, for $j=0,1$ set $A_{2j}(\Delta_{b;0})=\int_{S^{3}}a_{2j}(\Delta_{b;0})(x)$. In view of the definition of the constants $\gamma_{10}$ and 
 $\gamma_{11}'$ we have 
 \begin{equation}
     A_{0}(\Delta_{b;0})=\gamma_{10}\int_{S^{3}}d\theta\wedge \theta =\pi^{2}\gamma_{10}, \quad 
      A_{2}(\Delta_{b;0})=\gamma_{11}'\int_{S^{3}}R_{1}d\theta\wedge \theta =4\pi^{2}\gamma_{11}'.
      \label{eq:CR.gamma-A2j}
 \end{equation}
Notice that $A_{0}(\Delta_{b;0})$ and $A_{2}(\Delta_{b;0})$ are the coefficients of $t^{-2}$ and $t^{-1}$ in the asymptotics of $ \Tr 
 e^{-t\Delta_{b;0}}$ as $t\rightarrow 0^{+}$.  Moreover, we have $\Delta_{b;0} 
    =\boxdot_{\theta}-\frac{1}{4}R_{1}=\boxdot_{\theta}-1$, where $\boxdot_{\theta}$ denotes the CR invariant sublaplacian of Jerison-Lee~\cite{JL:YPCRM}, 
    and by~\cite[Thm. 4.34]{St:SICRM} 
     we have $\Tr e^{-t\boxdot_{\theta}} =\frac{\pi^{2}}{16t^{2}}+ 
    \op{O}(t^\infty)$ as $t \rightarrow 0^{+}$. Therefore, as $t\rightarrow 0^{+}$ we have
    \begin{equation}
        \Tr e^{-t\Delta_{b;0}}=e^{t}\Tr e^{-t\boxdot_{\theta}}\sim \frac{\pi^{2}}{16t^{2}}(1+t+\frac{t^{2}}{2}+\ldots).
    \end{equation}
Hence $A_{0}(\Delta_{b;0})=A_{2}(\Delta_{b;0})=\frac{\pi^{2}}{16}$. Combining this with~(\ref{eq:CR.gamma-A2j}) then shows that 
$\gamma_{10}=\frac{1}{16}$ and 
$\gamma_{11}'=\frac{1}{64}$, from which we get $\gamma_{1}''=\frac{1/64}{\sqrt{8.\frac{1}{16}}}=\frac1{32\sqrt2}$. Therefore, we get: 

\begin{theorem}\label{thm:spectral.area}
    If $\dim M=3$, then we have
    \begin{equation}
        \op{Area}_{\theta}M = \frac1{32\sqrt2}\int_{M}R_{1}d\theta \wedge\theta. 
        \label{eq:CR.area-dimension3} 
    \end{equation}
\end{theorem}

For instance, for $S^{3}$ equipped with its standard pseudohermitian structure we obtain 
$\op{Area}_{\theta}S^{3}=\frac{\pi^{2}}{8\sqrt{2}}$.

\appendix 
\section*{Appendix. Proof of Lemma~\ref{lem:Heisenberg.extension-symbol}} 
\setcounter{section}{1}
\setcounter{equation}{0}
In this appendix, for reader's convenience we give a detailed proof of Lemma~\ref{lem:Heisenberg.extension-symbol} about the extension of a homogeneous 
symbol on $\Rdo$ into a homogeneous distribution on $\Rd$. 

Let $p\in C^\infty(\Rdo)$ be homogeneous  of degree $m$, $m\in \C$, so that $p(\lambda.\xi)=\lambda^{m}p(\xi)$ for any $\lambda>0$.  
If $\Re m>-(d+2)$, then $p$ is integrable near the origin, so it defines a tempered distribution which is its unique homogeneous  extension. 

If $\Re m \leq -(d+2)$, then we can extend $p$ into the distribution $\tau\in\cS'(\Rd)$ defined by the formula,  
\begin{equation}
    \acou{\tau}{u}= \int [u(\xi)-\psi(\|\xi\|)\sum_{\brak\alpha\leq k} 
    \frac{\xi^\alpha}{\alpha!} u^{(\alpha)}(0)] p(\xi)d\xi \qquad \forall u\in\cS(\Rd), 
   \label{eq:Appendix.almosthomogeneous-extension}
\end{equation}
where $k$ is an integer $\geq -(\Re m +d+2)$ and $\psi$ is a function in $C_{c}^\infty(\R_{+})$ such that $\psi=1$ near $0$. 
Then in view of~(\ref{eq:PsiHDO.homogeneity-K-m}) for any $\lambda>0$ we have 
\begin{equation*}
   \begin{split}
 \acou{\tau_{\lambda}}{u}-\lambda^{m}  \acou{\tau}{u} &=     \lambda^{-(d+2)}\int [u(\lambda^{-1}.\xi)-\psi(\|\xi\|)\sum_{\brak\alpha\leq k} 
     \frac{\xi^\alpha\lambda^{-\brak\alpha}}{\alpha!} u^{(\alpha)}(0)]p(\xi)d\xi \\  
     &  -\lambda^{m}\int [u(\xi)-\psi(\|\xi\|)\sum_{\brak\alpha\leq k} 
    \frac{\xi^\alpha}{\alpha!} u^{(\alpha)}(0)] p(\xi)d\xi,\\
    &= \lambda^{m} \sum_{\brak\alpha\leq k} \frac{u^{(\alpha)}(0)}{\alpha!} \int [\psi(\|\xi\|)-\psi(\lambda\|\xi\|)] \xi^{\alpha}p(\xi)d\xi,\\
     &=  \lambda^{m} \sum_{\brak\alpha\leq k} \rho_{\alpha}(\lambda) c_{\alpha}(p) \acou{\delta^{(\alpha)}}{u}, 
   \end{split}
\end{equation*}
where we have let 
\begin{equation*}
c_{\alpha}(p) = \frac{(-1)^{|\alpha|}}{\alpha!}\int_{\|\xi\|=1}\xi^\alpha p(\xi)i_{E}d\xi, \qquad     
\rho_{\alpha}(\lambda)=\int_{0}^\infty \mu^{\brak\alpha+m+d+2} 
(\psi(\mu)-\psi(\lambda\mu)) \frac{d\mu}{\mu},
\end{equation*}
and, as in the statement of Lemma~\ref{lem:Heisenberg.extension-symbol}, $E$ is the vector field 
$2\xi_{0}\partial_{\xi_{0}}+\xi_{1}\partial_{\xi_{1}}+\ldots+\xi_{d}\partial_{\xi_{d}}$.

Set $\lambda=e^s$ and assume that $\psi$ is of the form $ \psi(\mu)=h(\log 
\mu) $ with $h\in C^\infty(\R)$ such that $h=1$ near $-\infty$ and $h=0$ 
near $+\infty$. Then, setting $a_{\alpha}=\brak\alpha+m+d+2$,  we have
\begin{equation}
    \frac{d}{ds}\rho_{\alpha}(e^s)=  \frac{d}{ds} \int_{-\infty}^{\infty}(h(t)-h(s+t))e^{a_{\alpha}t}dt=- e^{-as}\int_{-\infty}^{\infty} 
    e^{a_{\alpha}t} h'(t)dt.
     \label{eq:Appendix.differentiation-rhoalpha}
\end{equation}
 As $\rho_{\alpha}(1)=0$ it follows that $\tau$  is homogeneous of degree $m$  provided that 
\begin{equation}
 \int_{-\infty}^{\infty} e^{at} h'(t) ds = 0 \qquad \text{for $a=m+d+2, \ldots, m+d+2+k$}.    
    \label{eq:Appendix.homogeneous-extension-condition}
\end{equation}

    Next, if  $g\in C_{c}^\infty(\Rd)$ is such that $\int g(t)dt=1$, then for any $a \in \C\setminus 0$ we have 
    \begin{equation}
        \int _{-\infty}^{\infty} e^{at}(\frac{1}{a}\frac{d}{dt}+1)g(t)dt=0.
         \label{eq:Appendix1.integral-a-g}
    \end{equation}
    Therefore, if $m \not\in\Z$ then we can check that the conditions~(\ref{eq:Appendix.homogeneous-extension-condition}) are satisfied by 
\begin{equation}
    h'(t)=\prod_{a=m+d+2}^{m+d+2+k} (\frac{1}{a} \frac{d}{dt} +1)g(t).
     \label{eq:Appendix1.h'}
\end{equation}
As $\int_{-\infty}^{\infty}h'(t)dt=1$ we then see that the distribution $\tau$  defined 
by~(\ref{eq:Appendix.almosthomogeneous-extension}) with $\psi(\mu)=\int_{\log \mu}^\infty h'(t)dt$  is a homogeneous extension of $p(\xi)$.  

On the other hand, if $\tilde{\tau}\in\cS'(\Rd)$ is another homogeneous  extension of $p(\xi)$ 
then $\tau-\tau_{1}$ is supported at the origin, so we have $\tau=\tilde{\tau} + \sum b_{\alpha} 
\delta^{(\alpha)}$ for some constants $b_{\alpha}\in\C$. Then, for any $\lambda>0$, we have
\begin{equation}
    \tau_{\lambda}-\lambda^{m}\tau=\tilde{\tau}_{\lambda}-\lambda^{m}\tilde{\tau} + 
    \sum (\lambda^{-(d+2-\brak\alpha)} -\lambda^m) b_{\alpha}\delta^{(\alpha)}.
   \label{eq:Appendix.tau-l-tau-tilde}
\end{equation}
As both $\tau$ and $\tilde{\tau}$ are homogeneous of degree $m$, we deduce that 
$\sum (\lambda^{-(d+2-\brak\alpha)} -\lambda^m) b_{\alpha}\delta^{(\alpha)}=0$. The linear independence of the family $\{\delta^{(\alpha)}\}$ then 
implies that all the constants $b_{\alpha}$ vanish, that is, we have $\tilde{\tau}=\tau$. Thus $\tau$ is the unique homogeneous extension of $p(\xi)$ on $\Rd$. 
   
  Now, assume that $m$ is an integer $\leq -(d+2)$. Then in the formula~(\ref{eq:Appendix.almosthomogeneous-extension}) 
  for $\tau$ we can take $k=-(m+d+2)$ and let $\psi$ be of the form,
\begin{equation}
    \psi(\mu)=\int_{\log \mu}^\infty h'(t) dt, \qquad h'(t)=\prod_{a=m+d+2}^{m+d+2+k} (\frac{1}{a} \frac{d}{dt} +1) g(t),
\end{equation}  
   with $g\in C_{c}^\infty(\Rd)$ such that $\int g(t)dt=1$.  Then thanks to~(\ref{eq:Appendix.differentiation-rhoalpha}) and~(\ref{eq:Appendix1.integral-a-g}) we have 
$\rho_{\alpha}(\lambda)=0$ for $\brak \alpha<-(m+d+2)$, while for $\brak\alpha=-(m+d+2)$ we get
\begin{equation}
    \frac{d}{ds}\rho_{\alpha}(e^s)=\int h'(t)dt= \int g(t)dt=  1. 
\end{equation}
Since $\rho_{\alpha}(1)=0$ it follows that $\rho_{\alpha}(e^{s})=s$, that is, we have $\rho_{\alpha}(\lambda)=\log \lambda$. Thus, 
\begin{equation}
       \tau_{\lambda}=\lambda^{m}\tau +\lambda^{m}\log \lambda\sum_{\brak\alpha=-(m+d+2)} c_{\alpha}(p)\delta^{(\alpha)} \qquad \forall \lambda>0.
     \label{eq:Appendix.taulambda-tau}
\end{equation}
In particular, we see that  if all the coefficients $c_{\alpha}(p)$ with $\brak\alpha=-(m+d+2)$ vanish then $\tau$ is homogeneous of degree $m$. 

Conversely, suppose that $p(\xi)$ admits a homogeneous extension $\tilde{\tau}\in\cS'(\Rd)$. As $\tau-\tilde{\tau}$ is supported at $0$, we can write 
$\tau=\tilde{\tau} + \sum b_{\alpha} 
\delta^{(\alpha)}$ with $b_{\alpha}\in\C$. For any $\lambda>0$ we have $ \tilde{\tau}_{\lambda}=\lambda^{m}\tilde{\tau}$, so by combining this 
with~(\ref{eq:Appendix.tau-l-tau-tilde}) we get 
\begin{equation}
 \tau_{\lambda}-\lambda^{m}\tau  = \sum_{\brak\alpha\neq -(m+d+2)} 
    b_{\alpha}(\lambda^{-(\brak\alpha+d+2)}-\lambda^{m})\delta^{(\alpha)}.
\end{equation}
By comparing this with~(\ref{eq:Appendix.taulambda-tau}) and by using  linear independence of the family $\{\delta^{(\alpha)}\}$ we then deduce that 
we have  $c_{\alpha}(p)=0$ for $\brak\alpha= -(m+d+2)$. Therefore $p(\xi)$ 
admits a homogeneous extension if and only if all the coefficients $c_{\alpha}(p)$ with $\brak\alpha=-(m+d+2)$ vanish. The proof of 
Lemma~\ref{lem:Heisenberg.extension-symbol}  is thus  achieved.


\begin{thebibliography}{GJMS}  

\bibitem[BG]{BG:CHM} Beals, R.; Greiner, P.C.:  \emph{Calculus on Heisenberg manifolds}. 
 Annals of Mathematics Studies, vol. 119. Princeton University Press, Princeton, NJ, 1988.
 
\bibitem[BGS] {BGS:HECRM} Beals, R.; Greiner, P.C. Stanton, N.K.:  \emph{The heat equation on a CR manifold}. 
J. Differential Geom. \textbf{20} (1984), no. 2, 343--387.

\bibitem[Be]{Be:TSSRG} Bella{\"\i}che, A.:  \emph{The tangent space in sub-Riemannian geometry}. 
\emph{Sub-Riemannian geometry}, 1--78, Progr. Math., 144, Birkh\"auser, Basel, 1996.
 
\bibitem[BGV]{BGV:HKDO} Berline, N.; Getzler, E.; Vergne, M.: 
  \emph{Heat kernels and Dirac operators}. Grundlehren der Mathematischen Wissenschaften, vol. 298. Springer-Verlag, Berlin, 1992.
 
\bibitem[Bi]{Bi:MEAS} Biquard, O.: \emph{M\'etriques d'Einstein asymptotiquement
sym\'etriques.} Ast\'erisque No. 265 (2000). 

\bibitem[Bo]{Bo:PMIHUPCOED} Bony, J.M: \emph{Principe du maximum, in\'egalit\'e de Harnack et unicit\'e du probl\`eme de Cauchy pour les op\'erateurs 
elliptiques d\'eg\'en\'er\'es}. Ann. Inst. Fourier \textbf{19} (1969) 277--304.
  
  \bibitem[BdM]{BdM:HODCRPDO} Boutet de Monvel, L.: 
\emph{Hypoelliptic operators with double characteristics and related pseudo-differential operators.}
Comm. Pure Appl. Math. \textbf{27} (1974), 585--639.

\bibitem[BGu]{BG:STTO} Boutet de Monvel, L.; Guillemin, V. \emph{The spectral theory of Toeplitz operators}. Annals of Mathematics Studies, 99. 
Princeton University Press, Princeton, NJ, 1981.

\bibitem[CGGP]{CGGP:POGD}  Christ, M.; Geller, D.; G\l owacki, P.; Polin, L.: \emph{Pseudodifferential operators on groups with dilations.} 
Duke Math. J. \textbf{68} (1992) 31--65.

\bibitem[Co1]{Co:AFNG} Connes, A.: \emph{The action functional in noncommutative geometry}. 
Comm. Math. Phys. \textbf{117} (1988), no. 4, 673--683.

\bibitem[Co2]{Co:NCG} Connes, A.: \emph{Noncommutative geometry}. 
Academic Press, Inc., San Diego, CA, 1994.

\bibitem[Co3]{Co:GCMFNCG} Connes, A.: \emph{Gravity coupled with matter and the foundation of non-commutative geometry}. 
 Comm. Math. Phys. \textbf{182} (1996), no. 1, 155--176.

 \bibitem[CM]{CM:LIFNCG} Connes, A.; Moscovici, H.: \emph{The local index formula in noncommutative geometry}. 
 Geom. Funct. Anal. \textbf{5} (1995), no. 2, 174--243.

\bibitem[Di]{Di:ETNN} Dixmier, J.: \emph{Existence de traces non normales.} 
C. R. Acad. Sci. Paris S\'er. A-B \textbf{262} (1966) A1107--A1108.

\bibitem[Dy1]{Dy:POHG} Dynin, A.:  \emph{Pseudodifferential operators on the
Heisenberg group. } Dokl. Akad. Nauk SSSR \textbf{225} (1975) 1245--1248. 

\bibitem[Dy2]{Dy:APOHSC} Dynin, A.:  \emph{An algebra of pseudodifferential operators on
the Heisenberg groups. Symbolic calculus.} Dokl. Akad. Nauk
SSSR \textbf{227} (1976), 792--795.

\bibitem[ET]{ET:C} Eliashberg, Y.; Thurston, W.: \emph{Confoliations}. University Lecture Series, 13, AMS, Providence, RI, 1998.

\bibitem[EM]{EM:HAITH} Epstein, C.L.;  Melrose, R.B.:  
\emph{The Heisenberg algebra, index theory and homology}. Preprint, 1998.  Available at \texttt{http://www-math.mit.edu/$\tilde{}$rbm}. 

\bibitem[EMM]{EMM:RLSPD} Epstein, C.L.;  Melrose, R.B.; Mendoza, G.:   
\emph{Resolvent of the Laplacian on strictly pseudoconvex domains}. Acta Math. \textbf{167} (1991), no. 1-2, 1--106. 

\bibitem[FGLS] {FGLS:NRMB} Fedosov, B.V.;  Golse, F.; Leichtnam, E.; Schrohe, E.:
\emph{The noncommutative residue for manifolds with boundary}. 
J. Funct. Anal. \textbf{142} (1996), no. 1, 1--31. 

 \bibitem[Fe1]{Fe:MAEBKGPCD} Fefferman, C.: \emph{Monge-Amp\`ere equations, the Bergman kernel, and geometry of pseudoconvex domains}. 
 Ann. of Math. (2) \textbf{103} (1976), no. 2, 395--416. 
 
%

\bibitem[FS]{FS:FSSOSO} Fefferman, C.L.; S\'anchez-Calle, A.: \emph{Fundamental solutions for second order subelliptic operators}. 
Ann. of Math. (2) \textbf{124} (1986), no. 2, 247--272.

\bibitem[FSt]{FS:EDdbarbCAHG} Folland, G.; Stein, E.M.: \emph{Estimates for the $\bar \partial\sb{b}$ complex and analysis on the Heisenberg group.} 
Comm. Pure Appl. Math. \textbf{27} (1974) 429--522.

 \bibitem[Gi]{Gi:ITHEASIT} Gilkey, P.B.:  \emph{Invariance theory, the heat equation, and the Atiyah-Singer index theorem}. 
 Mathematics Lecture Series, 11. Publish or Perish, Inc., Wilmington, Del., 1984.
 
\bibitem[GK]{GK:ITLNSO} Gohberg, I.C.; Kre\u\i n, M.G.: 
\emph{Introduction to the theory of linear nonselfadjoint operators}. 
 Trans. of Math. Monographs 18, AMS, Providence, 1969. 

 \bibitem[GG]{GG:CRIPSL} Gover, A.R.; Graham, C.R.: 
\emph{CR invariant powers of the sub-Laplacian.} J. Reine Angew. Math. \textbf{583} (2005), 1--27.
 
 \bibitem[Gr]{Gr:FESPM} Greenleaf, A.: \emph{The first eigenvalue of a sub-Laplacian on a pseudo-Hermitian manifold}. 
 Comm. Partial Differential Equations \textbf{10} (1985), no. 2, 191--217.
 
\bibitem[Gro]{Gr:CCSSW} Gromov, M.: \emph{Carnot-Carath\'eodory spaces seen from within}. 
\emph{Sub-Riemannian geometry}, 79--323, Progr. Math., 144, Birkh\"auser, Basel, 1996.

\bibitem[Gu1]{Gu:NPWF} Guillemin, V.: 
\emph{A new proof of Weyl's formula on the asymptotic distribution of eigenvalues}. 
Adv. in Math. \textbf{55} (1985), no. 2, 131--160. 

\bibitem[Gu2]{Gu:GLD} Guillemin, V.W.:  \emph{Gauged Lagrangian distributions.} 
Adv. Math. \textbf{102} (1993), no. 2, 184--201.

\bibitem[Gu3]{Gu:RTCAFIO} Guillemin, V.: 
 \emph{Residue traces for certain algebras of Fourier integral operators}. 
 J. Funct. Anal. \textbf{115} (1993), no. 2, 391--417. 
   

\bibitem[HN]{HN:HMOPCV} Helffer, B.; Nourrigat, J.: \emph{Hypoellipticit\'e maximale pour des op\'erateurs polyn\^omes de champs de vecteurs.} Prog. 
Math., No. 58, Birkh\"auser, Boston, 1986.
 
 \bibitem[JL]{JL:YPCRM} Jerison, D.; Lee, J.M.: \emph{The Yamabe problem on CR manifolds}. 
J. Differential Geom. \textbf{25} (1987), no. 2, 167--197. 
  
 \bibitem[JK]{JK:OKTGSU} Julg, P.; Kasparov, G.: \emph{Operator $K$-theory for the group ${\rm SU}(n,1)$}.  
 J. Reine Angew. Math. \textbf{463} (1995), 99--152.
 
 \bibitem[KW]{KW:GNGWR} Kalau, W.; Walze, M.: \emph{Gravity, non-commutative geometry and the Wodzicki residue}. 
 J. Geom. Phys. \textbf{16} (1995), no. 4, 327--344. 
 
\bibitem[Ka]{Ka:RNC} Kassel, C.: \emph{Le r\'esidu non commutatif (d'apr\`es M. Wodzicki)}. 
 S\'eminaire Bourbaki, Vol. 1988/89. Ast\'erisque No. 177-178, (1989), Exp. No. 708, 199--229.
 
\bibitem[Kas]{Ka:DOG} Kastler, D.: \emph{The Dirac operator and gravitation}. 
 Comm. Math. Phys. \textbf{166} (1995), no. 3, 633--643. 
 
\bibitem[Ko]{Ko:BCM} Kohn, J.J.: \emph{Boundaries of complex manifolds}. 1965 Proc. Conf. Complex Analysis (Minneapolis, 1964) pp. 81--94
Springer, Berlin 

\bibitem[KR]{KR:EHFBCM} Kohn, J.J.; Rossi, H.: \emph{On the extension of holomorphic functions from the boundary of a complex manifold}.
Ann. of Math. \textbf{81} (1965) 451--472.

\bibitem[KV]{KV:GDEO} Kontsevich, M.; Vishik, S.: 
\emph{Geometry of determinants of elliptic operators.} 
 \emph{Functional analysis on the eve of the 21st century}, 
 Vol. 1 (New Brunswick, NJ, 1993), 173--197, Progr. Math., 131, 
 Birkh\"auser Boston, Boston, MA, 1995. 

 \bibitem[Le]{Le:FMPHI} Lee, J.M.:  \emph{The Fefferman metric and pseudo-Hermitian invariants}. 
Trans. Amer. Math. Soc. \textbf{296} (1986), no. 1, 411--429. 

 \bibitem[Les]{Le:NCRPDOLPS} Lesch, M.: \emph{On the noncommutative residue for pseudodifferential operators with log-polyhomogeneous symbols}. 
 Ann. Global Anal. Geom. \textbf{17} (1999), no. 2, 151--187.
 
\bibitem[Ma]{Ma:EPKLCD} Machedon, M.: \emph{Estimates for the parametrix of the Kohn Laplacian on certain domains}. 
 Invent. Math. \textbf{91} (1988), no. 2, 339--364.  


 \bibitem[MMS]{MMS:FIT} Mathai, V.; Melrose, R.; Singer, I.: \emph{Fractional index theory}. J. Differential 
 Geom. \textbf{74} (2006) 265--292.
 
 \bibitem[MN]{MN:HPDO1} Melrose, R.; Nistor, V.: \emph{Homology of pseudodifferential operators I. Manifolds with boundary}. 
 Preprint, arXiv, June `96.
 
\bibitem[NSW]{NSW:BMDVF1} Nagel, A.; Stein, E.M.; Wainger, S.: \emph{Balls and metrics defined by vector fields. I. Basic properties}. 
Acta Math. \textbf{155} (1985), no. 1-2, 103--147.

%

\bibitem[PR]{PR:CDBFCF} Paycha, S.; Rosenberg, S.: \emph{Curvature on determinant bundles and first Chern forms}. 
J. Geom. Phys. \textbf{45} (2003), no. 3-4, 393--429.

\bibitem[Po1]{Po:CRAS1} Ponge, R.: \emph{Calcul fonctionnel sous-elliptique et r\'esidu non commutatif sur les vari\'et\'es de 
Heisenberg}. C.~R. Acad. Sci. Paris, S\'erie I, \textbf{332} (2001)  611--614. 

\bibitem[Po2]{Po:CRAS2} Ponge, R.:  \emph{G\'eom\'etrie spectrale et formules d'indices locales pour les vari\'et\'es CR et contact}. 
     C. R. Acad. Sci. Paris, S\'erie I, \textbf{332} (2001) 735--738. 
  
 \bibitem[Po3]{Po:IJM1} Ponge, R.: \emph{Spectral asymmetry, zeta functions and the noncommutative residue}. 
 Int. J. Math. \textbf{17} (2006), 1065-1090.

\bibitem[Po4]{Po:Pacific1} Ponge, R.: \emph{The tangent groupoid of a Heisenberg manifold.} Pacific 
Math. J. \textbf{227} (2006) 151--175.

\bibitem[Po5]{Po:MAMS1}  Ponge, R.: \emph{Heisenberg calculus and spectral theory of hypoelliptic operators on Heisenberg manifolds.} 
E-print, arXiv, Sep.~05, 140 pages. To appear in Mem. Amer. Math. Soc..
 
\bibitem[Po6]{Po:Crelle2} Ponge, R.: \emph{Noncommutative residue invariants for CR and contact manifolds.} E-print, arXiv, Oct.~05, 30 pages. To 
appear in  J. Reine Angew. Math..
%
 
 
 \bibitem[Po7]{Po:LMP07} Ponge, R.: \emph{Noncommutative geometry and lower dimensional volumes in Riemannian geometry}. E-print, 
 arXiv, July~07.
 
 \bibitem[Po8]{Po:CPDE1}  Ponge, R.: \emph{Hypoelliptic functional calculus on Heisenberg manifolds. A resolvent approach.} 
 E-print, arXiv, Sep.~07.
 
 \bibitem[Ro1]{Ro:HHGRTC} Rockland, C.: \emph{Hypoellipticity on the Heisenberg group-representation-theoretic criteria.}  
 Trans. Amer. Math. Soc.  \textbf{240}  (1978) 1--52. 
  
\bibitem[Ro2]{Ro:INA}ÊRockland, C.: \emph{Intrinsic nilpotent approximation.} Acta Appl. Math. \textbf{8} (1987), no. 3, 213--270. 

\bibitem[RS]{RS:HDONG} Rothschild, L.; Stein, E.: \emph{Hypoelliptic differential operators and nilpotent groups.}
Acta Math. \textbf{137} (1976) 247--320.

 \bibitem[Ru]{Ru:FDVC} Rumin, M.: \emph{Formes diff\'erentielles sur les vari\'et\'es de contact}. 
 J. Differential Geom. \textbf{39} (1994), no.2, 281--330.
  
\bibitem[Sa]{Sa:FSGSSVF} S\'anchez-Calle, A.: \emph{Fundamental solutions and geometry of the sum of squares of vector fields}. 
Invent. Math. \textbf{78} (1984), no. 1, 143--160.

 \bibitem[Sc]{Sc:NCRMCS} Schrohe, E.: \emph{Noncommutative residues and manifolds with conical singularities}. 
 J. Funct. Anal. \textbf{150} (1997), no. 1, 146--174.
 
 \bibitem[St]{St:SICRM}  Stanton, N.K.: \emph{Spectral invariants of CR manifolds}. 
 Michigan Math. J. \textbf{36} (1989), no. 2,  267--288.
 
 \bibitem[Ta]{Ta:DGSSPCM} Tanaka, N.: 
 \emph{A differential geometric study on strongly pseudo-convex manifolds}. 
Lectures in Mathematics, Department of Mathematics, Kyoto University, No. 9.  Kinokuniya Book-Store Co., Ltd., Tokyo, 1975.
 
\bibitem[Tay]{Ta:NCMA} Taylor, M.E.: 
\emph{Noncommutative microlocal analysis. I.} Mem. Amer. Math. Soc. 52 (1984), no. 313,

 \bibitem[We]{We:PHSRH} Webster, S.: \emph{Pseudo-Hermitian structures on a real hypersurface}. 
 J. Differential Geom. \textbf{13} (1978), no. 1, 25--41.
 
\bibitem[Va]{Va:PhD} Van Erp, E.: PhD thesis, Pennsylvania State University, 2005.

\bibitem[Vas]{Vas:PhD} Vassout, S.: \emph{Feuilletages et r\'esidu non commutatif longitudinal}. PhD thesis, University of Paris 7, 2001.

\bibitem[Wo1]{Wo:LISA} Wodzicki, M.: \emph{Local invariants of spectral asymmetry}. 
Invent. Math. \textbf{75} (1984), no. 1, 143--177. 

 \bibitem[Wo2]{Wo:PhD} Wodzicki, M.: \emph{Spectral asymmetry and noncommutative residue} (in Russian), 
 Habilitation Thesis, Steklov Institute, (former) Soviet Academy of Sciences, Moscow, 1984.

\bibitem[Wo3]{Wo:NCRF} Wodzicki, M.: \emph{Noncommutative residue. I. Fundamentals}. 
\emph{$K$-theory, arithmetic and geometry} (Moscow, 1984--1986), 320--399, Lecture Notes in Math., 1289, Springer, Berlin-New York, 1987. 

  \bibitem[Wo4]{Wo:LESCHAEA} Wodzicki, M.: \emph{The long exact sequence in cyclic homology associated with extensions of algebras}. 
 C. R. Acad. Sci. Paris S\'er. I Math. \textbf{306} (1988), no. 9, 399--403.
\end{thebibliography}
\end{document}